\documentclass[10pt,twocolumn,letterpaper]{article}

\usepackage{cvpr}

\usepackage{times}
\usepackage{epsfig}
\usepackage{graphicx}
\usepackage{amsmath}
\usepackage{amssymb}
\usepackage{times}  
\usepackage{helvet}  
\usepackage{courier}  
\usepackage{url}  
\usepackage{graphicx}  
\usepackage{cancel,soul,ulem}
\usepackage{enumitem}
\usepackage{natbib}
\bibpunct{[}{]}{,}{n}{}{}
\usepackage{xcolor}
\usepackage{soul}

\usepackage{amsmath}
\usepackage{amssymb}
\usepackage{url}
\usepackage{amsthm} 
\usepackage{algorithm}
\usepackage{algorithmic}
\usepackage{diagbox} 
\usepackage{bbm}
\usepackage{bm}
\usepackage{diagbox}
\usepackage{graphicx}
\usepackage{subcaption}
\usepackage{epstopdf}
\usepackage{booktabs}
\usepackage{pdfpages}

\usepackage{courier}
\usepackage{color}
\definecolor{colorpink}{RGB}{251,53,155}
\definecolor{colorblue}{RGB}{0,148,200}
\definecolor{colorgreen}{RGB}{0,150,0}

\def\beq{\begin{eqnarray}}
\def\eeq{\end{eqnarray}}
\def\noi{\noindent}
\def\nn{\nonumber}
\def\la{\langle}
\def\ra{\rangle}

\newtheorem{lemma}{Lemma}
\newtheorem{theorem}{Theorem}
\newtheorem{definition}{Definition}
\newtheorem{proposition}{Proposition}

\newcommand{\bbb}[1]{\boldsymbol{\mathbf{#1}}}
\def\S{\bar{\mathcal{S}}}
\def\Z{\bar{\mathcal{Z}}}

\newcommand{\tinytiny}{\fontsize{4}{4}\selectfont}

\def\threefigwid{0.3\textwidth}
\def\fourfigwid{0.245\textwidth}
\def\objimghei{0.16\textheight}
\def\ghs{\hspace{0.05cm}} 

\def\E{\mathbb{E}}

\usepackage{multirow}
\usepackage{rotating}

\usepackage[pagebackref=true,breaklinks=true,letterpaper=true,colorlinks,bookmarks=false]{hyperref}

\cvprfinalcopy 

\ifcvprfinal\fi
\begin{document}

\title{A Decomposition Algorithm for the Sparse Generalized Eigenvalue Problem}

\author{Ganzhao Yuan$^{1}$, Li Shen$^{2}$, Wei-Shi Zheng$^{3}$\\
$^1$ Center for Quantum Computing, Peng Cheng Laboratory, Shenzhen 518005, China\\
$^2$ Tencent AI Lab, Shenzhen, China\\
$^3$ School of Data and Computer Science, Sun Yat-sen University, China\\
{\tt\small yuanganzhao@foxmail.com,~mathshenli@gmail.com,~zhwshi@mail.sysu.edu.cn}
}

\maketitle
\begin{abstract}
The sparse generalized eigenvalue problem arises in a number of standard and modern statistical learning models, including sparse principal component analysis, sparse Fisher discriminant analysis, and sparse canonical correlation analysis. However, this problem is difficult to solve since it is NP-hard. In this paper, we consider a new decomposition method to tackle this problem. Specifically, we use random or/and swapping strategies to find a working set and perform global combinatorial search over the small subset of variables. We consider a bisection search method and a coordinate descent method for solving the quadratic fractional programming subproblem. In addition, we provide some theoretical analysis for the proposed method. Our experiments have shown that the proposed method significantly and consistently outperforms existing solutions in term of accuracy.





\end{abstract}


\section{Introduction}
In this paper, we mainly focus on the following sparse generalized eigenvalue problem (`$\triangleq$' means define):
\beq \label{eq:main}
\begin{split}
\min_{\bbb{x} \neq \bbb{0},{\bbb{x}\in\Omega}}~f(\bbb{x})\triangleq \frac{h(\bbb{x})}{g(\bbb{x})},~~\text{with}~~\Omega \triangleq \{ \bbb{x}~|~\|\bbb{x}\|_0 \leq s\},\\
h(\bbb{x}) \triangleq \tfrac{1}{2}\bbb{x}^T\bbb{A}\bbb{x},~~g(\bbb{x}) \triangleq\tfrac{1}{2} \bbb{x}^T\bbb{C}\bbb{x} .~~~~~~~~~~~~~~~
\end{split}
\eeq
\noi Here, $\bbb{x}\in \mathbb{R}^n$, and $\|\cdot \|_0$ is a function that counts the number of nonzero elements in a vector. $\bbb{A} \in \mathbb{R}^{n\times n}$ and $\bbb{C}\in \mathbb{R}^{n\times n}$ are some symmetry matrices. We assume that $\bbb{C}$ is strictly positive definite and $s \in [1,n]$ is a positive integer.


The sparse generalized eigenvalue problem in (\ref{eq:main}) describes many applications of interest in both computer vision and machine learning, including object recognition \cite{NaikalYS11}, visual tracking \cite{KwonL10}, object detection \cite{PaisitkriangkraiSZ09,PaisitkriangkraiSZ11,ShenPZ11}, pixel/part selection \cite{MoghaddamWA07}, and text summarization \cite{zhang2011large}. We notice that the objective function and sparsity constraint in Problem (\ref{eq:main}) is scale-invariant (multiplying $\bbb{x}$ with a positive constant does not change the value of the objective function and the satisfiability of the sparsity constraint). Thus, it is equivalent to the following problem: $\min_{\bbb{x}}~\bbb{x}^T\bbb{A}\bbb{x},~s.t.~\bbb{x}^T\bbb{C}\bbb{x}=1,~\|\bbb{x}\|_0\leq s$. Moreover, without the sparsity constraint, Problem (\ref{eq:main}) reduces to the minimum generalized eigenvalue problem and it has several equivalent formulations \cite{beck2010minimizing}: $\min_{\bbb{x} \neq \bbb{0}}~{(\bbb{x}^T\bbb{Ax})}/{(\bbb{x}^T\bbb{Cx})}=\min\{\bbb{x}^T\bbb{Ax}:\bbb{x}^T\bbb{Cx}=1\}=\max \{\lambda:\bbb{A} - \lambda \bbb{C} \succeq \bbb{0}\} =  \lambda_{\min}(\bbb{C}^{-1/2}\bbb{AC}^{-1/2})$, where $\lambda_{\min}(\bbb{X})$ is the smallest eigenvalue of a given matrix $\bbb{X}$, and $\bbb{X}\succeq \bbb{0}$ denotes $\bbb{X}$ is positive semidefinite.


Problem (\ref{eq:main}) is closely related to the classical matrix computation in the literature \cite{horn1990matrix,ge2016efficient,asteris2016simple}. Imposing an additional sparsity constraint on the solution reduces over-fitting and improves the interpretability of the model for high-dimensional data analysis. This has evoked great research interests in developing methods that enforce sparsity on eigenvectors. For instance,
the work of \cite{jolliffe2003modified} successively choose a sparse principle direction that maximizes the variance by enforcing a sparsity constraint using a bounded $\ell_1$ norm. The work of \cite{zou2006sparse} reformulates the principle component analysis problem as a elastic-net regularized ridge regression problem, which can be solved efficiently using least angle regression. The work of \cite{d2005direct} proposes a convex relaxation for the sparse principle component analysis problem problem based on semidefinite programming.

One difficulty of solving Problem (\ref{eq:main}) comes from the combinational nature of the cardinality constraint. A conventional way to solve this problem is to simply replace $\ell_0$ norm by its convex relaxation. Recently, non-convex approximation methods such as Schatten $\ell_p$ norm, re-weighted $\ell_1$ norm, Capped $\ell_1$ function have been proposed for acquiring better accuracy \cite{SongBP15}. However, all these approximation methods fails to directly control the sparsity of the solution. In contrast, iterative hard thresholding maintain the sparsity of the solution by iteratively setting the small elements (in magnitude) to zero in a gradient descent manner. Due to its simplicity, it has been widely used and incorporated into the truncated power method \cite{yuan2013truncated} and truncated Rayleigh flow method \cite{tan2016sparse}. 


Another difficulty of solving Problem (\ref{eq:main}) is due to the non-convexity of the objective function. One popular method to overcome this difficulty is removing the quadratic term using semidefinite programming lifting technique and reformulating (\ref{eq:main}) into the following low-rank sparse optimization problem: $\min_{\bbb{X} \neq \bbb{0}}$~${tr(\bbb{AX})}/{tr(\bbb{CX})}$,~$s.t.~$$\bbb{X}\succeq 0$,~$rank(\bbb{X})=1$,~$\|\bbb{X}\|_0\leq s^2$. We remark that the objective function is quasilinear (hence both quasiconvex and quasiconcave), and one can constrain the denominator to be a positive constant using the scale-invariant property of the problem. Recently, convex semidefinite programming method drops the rank constraint and considers $\ell_1$ relaxation for the sparsity constraint \cite{d2005direct,d2008optimal,krauthgamer2015semidefinite,zhang2012sparse}. It has been shown to achieve strong guarantee under suitable assumptions. However, such a matrix lifting technique will incur expensive computation overhead.




In summary, existing methods for solving Problem (\ref{eq:main}) suffer from the following limitations. \bbb{(i)} Semidefinite programming methods \cite{d2005direct,d2008optimal,krauthgamer2015semidefinite} are not scalable due to its high computational complexity for the eigenvalue decomposition. \bbb{(ii)} Convex/Non-convex approximation methods \cite{Sriperumbudur2007,thiao2010dc,sriperumbudur2011majorization} fail to directly control the low-rank and sparse property of the solution. \bbb{(iii)} Hard thresholding methods \cite{yuan2013truncated,journee2010generalized} only obtain weak optimality guarantee and result in poor accuracy in practice \cite{beck2013sparsity,yuan2017hybrid}.

Recently, the work of \cite{beck2016sparse} considers a new optimality criterion which is based on Coordinate-Wise Optimality (CWO) condition for sparse optimization. It is proven that CWO condition is stronger than the optimality criterion based on hard thresholding. The work of \cite{yuan2017hybrid} presents a new block-$k$ optimal condition for general discrete optimization. It is shown to be stronger than CWO condition since it includes CWO condition as a special case with $k=1$. Inspired by these works, we propose a new decomposition method for the sparse generalized eigenvalue problem, along with using a greedy method based on CWO \cite{beck2016sparse} for finding the working set.






\textbf{Contributions:} This paper makes the following contributions. \bbb{(i)} We propose a new decomposition algorithm for solving the sparse generalized eigenvalue problem (see Section \ref{sect:dec}). \bbb{(ii)} We discuss two strategies to find the working set for our decomposition algorithm (see Section \ref{sect:work}). \bbb{(iii)} We propose two methods to solve the sparse quadratic fractional programming subproblem (see Section \ref{sect:subp}). \bbb{(iv)} A convergence analysis for the decomposition method is provided (see Section \ref{sect:conv}). \bbb{(v)} Our experiments have shown that our method outperforms existing solutions in term of accuracy. (see Section \ref{sect:exp}).

\textbf{Notation:} All vectors are column vectors and superscript $T$ denotes transpose. $\bbb{X}_{i,j}$ denotes the ($i^{\text{th}}$, $j^{\text{th}}$) element of matrix $\bbb{X}$ and $\bbb{x}_i$ denotes the $i$-th element of vector $\bbb{x}$. $\bbb{e}_i$ is a unit vector with a $1$ in the $i^{th}$ entry and $0$ in all other entries. $\text{diag}(\bbb{x})$ is a diagonal matrix formed with $\bbb{x}$ as its principal diagonal. For any partition of the index vector $[1,2,...,n]$ into $[B,N]$ with $B\in \mathbb{N}^{k},~N\in \mathbb{N}^{n-k}$, we define $\bbb{U}_B \in \mathbb{R}^{n\times k},~\bbb{U}_N\in \mathbb{R}^{n\times (n-k)}$ as: {\tiny
$(\bbb{U}_B)_{j,i}=\left\{
             \begin{array}{ll}
               1, & \hbox{$B(i)=j$;} \\
               0, & \hbox{else.}
             \end{array}
           \right.,~(\bbb{U}_N)_{j,l}=\left\{
             \begin{array}{ll}
               1, & \hbox{$N(l)=j$;} \\
               0, & \hbox{else.}
             \end{array}
           \right.
$}\normalsize. Therefore, we have $\bbb{x}_B=\bbb{U}_B^T\bbb{x}$ and $\bbb{x}= \bbb{I}\bbb{x}=(\bbb{U}_B\bbb{U}_B^T+\bbb{U}_N\bbb{U}_N^T)\bbb{x}=\bbb{U}_B\bbb{x}_B+\bbb{U}_N\bbb{x}_N$. Finally, $C_n^{k}$ denotes the number of possible combinations choosing $k$ items from $n$.

\section{Generalized Eigenvalue Problems} \label{sect:app}
A number of standard and modern statistical learning models can be formulated as the sparse generalized eigenvalue problem, which we present some instances below.

$\bullet$ \textbf{Principle Component Analysis (PCA)}. Consider a data matrix $\bbb{Z} \in \mathbb{R}^{m\times d}$, where each row represents an independent sample. The covariance matrix $\bbb{\Sigma}$ is computed by $\bbb{\Sigma} = \tfrac{1}{m-1}\sum_{i=1}^m(\bbb{z}_i-\bbb{\mu})(\bbb{z}_i-\bbb{\mu})^T \in \mathbb{R}^{d\times d}$, where $\bbb{z}_i$ denotes $i^{\text{th}}$ column of $\bbb{Z}$ and $\bbb{\mu} = \sum_{i=1}^m \bbb{z}_i \in \mathbb{R}^{d}$. PCA can be cast into the following optimization problem: $\min_{\bbb{x}\neq 0}~(-\bbb{x}^T \bbb{\Sigma}\bbb{x})/(\bbb{x}^T \bbb{x})$.

$\bullet$ \textbf{Fisher Discriminant Analysis (FDA)}. Given observations with two distinct classes with $\bbb{\mu}_{(i)}$ and $\bbb{\Sigma}_{(i)}$ being the mean vector and covariance matrix of class $i$ ($i=1$~\text{or}~$2$), respectively. FDA seeks a projection vector such that the between-class variance is large relative to the within-class variance, leading to the following problem: $\min_{\bbb{x}\neq 0}~\frac{-\bbb{x}^T ((\bbb{\mu}_{(1)}-\bbb{\mu}_{(2)})(\bbb{\mu}_{(1)}-\bbb{\mu}_{(2)})^T)\bbb{x}}{\bbb{x}^T (\bbb{\Sigma}_{(1)}+\bbb{\Sigma}_{(2)})\bbb{x}}$.

$\bullet$ \textbf{Canonical Correlation Analysis (CCA)}. Given two classes of data $\bbb{X}\in\mathbb{R}^{m_1 \times d}$ and $\bbb{Y}\in\mathbb{R}^{m_2 \times d}$, the covariance matrix between samples from $\bbb{X}$ and $\bbb{Y}$ can be constructed as $\bbb{\Sigma} \triangleq {\tinytiny \begin{pmatrix}
                                                 \bbb{\Sigma}_{xx} &  \bbb{\Sigma}_{xy} \\
                                                 \bbb{\Sigma}_{yx} & \bbb{\Sigma}_{yy} \\
                                               \end{pmatrix}} \in \mathbb{R}^{(m_1+m_2)\times(m_1+m_2)}$ with $\bbb{\Sigma}_{xy}\in \mathbb{R}^{m_1\times m_2}$. CCA exploits the relation of the samples by solving the following problem: $\max_{\bbb{u}\neq 0,~\bbb{v}\neq 0}~\bbb{u}^T \bbb{\Sigma}_{xy} \bbb{v}, ~s.t.~ \bbb{u}^T\bbb{\Sigma}_{xx}\bbb{u} = \bbb{v}^T\bbb{\Sigma}_{yy}\bbb{v}=1$, where $\bbb{A} \triangleq {\tinytiny \begin{pmatrix}
                                                 \bbb{0}&  \bbb{\Sigma}_{xy} \\
                                                 \bbb{\Sigma}_{yx} & \bbb{0}\\
                                               \end{pmatrix}}$, $\bbb{C}\triangleq {\tinytiny \begin{pmatrix}
                                                 \bbb{\Sigma}_{xx} &  \bbb{0}\\
                                                 \bbb{0} & \bbb{\Sigma}_{yy} \\
                                               \end{pmatrix}}\in \mathbb{R}^{(m_1+m_2)\times(m_1+m_2)}$, and $\bbb{x}\triangleq {  [ \bbb{u}^T \bbb{v}^T]^T}$. One can rewrite CCA as the following equivalent problem: $\min_{\bbb{x}}~\frac{-\bbb{x}^T\bbb{Ax}}{\bbb{x}^T\bbb{Cx}}$.

Incorporated with the sparsity constraint, the applications listed above become special cases of the general optimization models in (\ref{eq:main}).


\section{The Proposed Decomposition Algorithm} \label{sect:dec}

This section presents our decomposition algorithm for solving (\ref{eq:main}), which is based on the following notation of block-$k$ optimality \cite{yuan2017hybrid} for general non-convex constrained optimization.

\begin{definition} \label{def:block:k}

(Block-$k$ Optimal Solution and Block-$k$ Optimality Measure ) \bbb{(i)} We denote $B \in \mathbb{N}^{k}$ as a vector containing $k$ unique integers selected from $\{1,2,...,n\}$. We define $N\triangleq \{1,2,...,n\}\setminus B,~\bbb{x} = \bbb{U}_B\bbb{x}_B + \bbb{U}_N \bbb{x}_N$ and let
\begin{align} \label{eq:welldefine}
\begin{split}
  \textstyle \mathcal{P}(B,\bbb{x}) \triangleq &\arg \textstyle\min_{\bbb{x}_{B}}~f\left(\bbb{U}_B \bbb{x}_B + \bbb{U}_N \bbb{x}_N\right), \\
& \textstyle s.t.~\left(\bbb{U}_B \bbb{x}_B + \bbb{U}_N \bbb{x}_N\right) \in\Omega.
\end{split}
\end{align}
\noi Assume that $\mathcal{P}(B,\bbb{x})$ in (\ref{eq:welldefine}) is bounded for all $\bbb{x}$ and $B$. A solution $\bar{\bbb{x}}$ is the block-$k$ optimal solution if and only if $\bar{\bbb{x}}_{B}= \mathcal{P}(B,\bar{\bbb{x}}) \text{~for all~} |B|=k$. In other words, a solution is the block-$k$ optimal solution if and only if every block coordinate of size $k$ achieves the global optimal solution. \bbb{(ii)} We define $\mathcal{M}(\bbb{x}) \triangleq \frac{1}{C_n^k} \sum_{i=1}^{C_n^k}~\| \mathcal{P}(\mathcal{B}_{(i)},~\bbb{x})-\bbb{x}_{\mathcal{B}_{(i)}}\|_2^2$ with {\scriptsize$\{\mathcal{B}_{(i)}\}_{i=1}^{C_n^k}$} being all the possible combinations of the index vectors choosing $k$ items from $n$ with $\mathcal{B}_{(i)} \in \mathbb{N}^{k}$ for all $i$. $\mathcal{M}(\bbb{x})$ is an optimality measure for Problem \ref{eq:main} in the sense that $\mathcal{M}(\bar{\bbb{x}})=0$ if and only if $\bar{\bbb{x}}$ is the block-$k$ optimal solution.




%
%
%

\end{definition}

We describe the basic idea of the decomposition method. In each iteration $t$, the indices $\{1,2,...,n\}$ of decision variable are separated to two sets $B^t$ and $N^t$, where $B^t$ is the working set and $N^t=\{1,2,...,n\} \setminus B^t$. To simplify the notation, we use $B$ instead of $B^t$. Therefore, we can rewrite $h(\cdot)$ and $g(\cdot)$ in Problem (\ref{algo:main}) as:
{\small\beq
 h(\bbb{x}_B,\bbb{x}_N) = \tfrac{1}{2} \bbb{x}_B^T\bbb{A}_{BB}\bbb{x}_B + \tfrac{1}{2}\bbb{x}_N^T\bbb{A}_{NN}\bbb{x}_N + \la \bbb{x}_B, \bbb{A}_{BN}\bbb{x}_N \ra, \nn\\
 g(\bbb{x}_B,\bbb{x}_N) = \tfrac{1}{2}\bbb{x}_B^T\bbb{C}_{BB}\bbb{x}_B + \tfrac{1}{2}\bbb{x}_N^T\bbb{C}_{NN}\bbb{x}_N + \la \bbb{x}_B, \bbb{C}_{BN}\bbb{x}_N \ra. ~\nn
\eeq}{\normalsize}\noi The vector $\bbb{x}_N$ is fixed so the objective value becomes a subproblem with the variable $\bbb{x}_B$. Our proposed algorithm iteratively solves the small-sized optimization problem with respect to the variable $\bbb{x}_B$ as in (\ref{eq:subprob}) until convergence. We summarize our method in Algorithm \ref{algo:main}.


\begin{algorithm}[!h]
\caption{ {\bf A Decomposition Algorithm for Sparse Generalized Eigenvalue Problem as in (\ref{eq:main}).} }
\begin{algorithmic}[1]
  \STATE Specify the working set parameter $k$ and the proximal term parameter $\theta$. Find {an} initial feasible solution $\bbb{x}^0$ and set $t=0$.
  \WHILE{not converge}
  \STATE (\bbb{S1}) Use some strategy to find a working set $B$ whose size is $k$. Define $N\triangleq \{1,2,...,n\}\setminus B$.
  \STATE (\bbb{S2}) Solve the following subproblem with the variable $\bbb{x}_B$ using combinatorial search:
  \beq \label{eq:subprob}
  \begin{split}
    \bbb{x}_{B}^{t+1} \Leftarrow \arg \min_{\bbb{x}_B}~\frac{h(\bbb{x}_B,\bbb{x}^t_N) + \frac{\theta}{2}\|\bbb{x}_B-\bbb{x}^t_B\|_2^2}{g(\bbb{x}_B,\bbb{x}^t_N)}\\
    s.t.~\|\bbb{x}_B\|_0  +\|\bbb{x}^t_N\|_0\leq s~~~~~~~~~~~~
  \end{split}
\eeq
  \STATE (\bbb{S3}) Increment $t$ by 1.
  \ENDWHILE
\end{algorithmic}\label{algo:main}
\end{algorithm}

\noi\bbb{Remarks.} \bbb{(i)} The concept of block-$k$ optimality has been introduced in \cite{yuan2017hybrid}. This paper extends their method for minimizing convex functions to handle general non-convex objective functions. \bbb{(ii)} Algorithm \ref{algo:main} relies on solving a small-sized quadratic fractional problem as in (\ref{eq:subprob}). However, using the specific structure of the objective function and the sparsity constraint, we can develop an efficient and practical algorithm to solve it globally. \bbb{(iii)} We propose a new proximal strategy when solving the subproblem as in (\ref{eq:subprob}). Note that the proximal strategy is only applied to the numerator instead of to the whole objective function. This is to guarantee sufficient descent condition and global convergence of Algorithm \ref{algo:main} (see Lemma \ref{lemma:suf:dec} and Theorem \ref{theorem:convergence}). \bbb{(iv)} When the dimension $n$ is small \footnote{For example, the popular pit props data set \cite{jeffers1967two,mackey2009deflation} only contains 13 dimensions.} and the parameter setting $\theta=0,~k=n$ is used, the subproblem in (\ref{eq:subprob}) is equivalent to Problem (\ref{eq:main}).


\section{Finding the Working Set} \label{sect:work}

This section shows how to find the working set (refer to Step \bbb{S1} in Algorithm \ref{algo:main}). This problem is challenging for two aspects. \bbb{(i)} Unlike convex methods that one can find the working set using the first-order optimal condition or KKT primal-dual residual \cite{joachims1998making,chang2011libsvm}, there is no general criteria to find the working set for non-convex problems. \bbb{(ii)} There are $C^k_n$ possible combinations of choice for the working set of size $k$. One cannot expect to use the cyclic strategy and alternatingly minimize over all the possible combinations (i.e., {\scriptsize $\{\mathcal{B}_{(i)}\}_{i=1}^{C_n^k}$ }) due to its high computational complexity when $k$ is large. We propose the following two strategies to find the working set:

$\bullet$ Random Strategy. We uniformly select one combination (which contains $k$ coordinates) from {\scriptsize $\{\mathcal{B}_{(i)}\}_{i=1}^{C_n^k}$ }. {In expectation}, our algorithm is still guaranteed to find the block-$k$ stationary point.

$\bullet$ Swapping Strategy. We denote $\mathcal{S}(\bbb{x})$ and $\mathcal{Z}(\bbb{x})$ as the index of non-zero elements and zero elements of $\bbb{x}$, respectively. Based on the current solution $\bbb{x}^t$, our method enumerates all the possible pairs $(i,j)$ with $i \in \mathcal{S}(\bbb{x}^t),~j\in \mathcal{Z}(\bbb{x}^t)$ that lead to the greatest descent $\bbb{D}_{i,j}$ by changing the two coordinates from zero/non-zero to non-zero/zero, as follows:
\beq \label{eq:one:dim:selection}
\textstyle \bbb{D}_{i,j} = \min_{\beta}~f(\bbb{x}^t + \beta \bbb{e}_i - \bbb{x}_j^t \bbb{e}_j) - f(\bbb{x}^t).
\eeq
\noi We then pick the top pairs of coordinates that lead to the greatest descent by measuring $\bbb{D}\in\mathbb{R}^{|\mathcal{S}(\bbb{x})|\times |\mathcal{Z}(\bbb{x})|}$. Specifically, we sort the elements in $\bbb{D}$ with $\bbb{D}_{P_1,S_1}\leq \bbb{D}_{P_2,S_2}\leq \bbb{D}_{P_3,S_3} \leq ,..., \bbb{D}_{P_n,S_n}$, where $P \in \mathbb{N}^{n}$ and $S \in \mathbb{N}^{n}$ are the index vectors. Assuming that $k$ is an even number, we simply pick the top-$(k/2)$ nonoverlapping elements of the sequence $P$ and $S$ respectively as the working set.


We now discuss how to solve (\ref{eq:one:dim:selection}) to obtain $\bbb{D}_{i,j}$. We start from the following lemma.

\begin{lemma} \label{lemma:line}
We consider the following one-dimensional optimization problem:
\beq \label{eq:one:dim}
\beta^* = \arg \min_{\beta}\psi(\beta) \triangleq \frac{\tfrac{1}{2}\bar{a} \beta^2  +\bar{b} \beta  + \bar{c}}{\tfrac{1}{2}\bar{r} \beta^2  + \bar{s} \beta + \bar{t}},~s.t.~\beta \geq \bar{L}
\eeq
\noi Assume that $\forall \beta\geq \bar{L},~\tau \triangleq \tfrac{1}{2}\bar{r}\beta^2   + \bar{s}\beta  + \bar{t} >0$ and the optimal solution is bounded. We have: \bbb{(i)} Problem (\ref{eq:one:dim}) admits a closed-form solution as: $\beta^* = \arg \min_{\beta}~f(\beta),~\beta \in\{ \Pi(\beta_1),\Pi(\beta_2)\}$, where $\beta_1={(-\vartheta - \sqrt{\vartheta^2 -2\pi \iota })}/{\pi},~\beta_2={(-\vartheta + \sqrt{\vartheta^2 -2\pi \iota })}/{\pi}$,~$\pi \triangleq \bar{a}\bar{s}-\bar{b}\bar{r},~\vartheta \triangleq \bar{a}\bar{t} - \bar{c}\bar{r},~\iota \triangleq \bar{t}\bar{b} + \bar{c}\bar{s}$, and $\Pi(\beta) \triangleq  \max(\bar{L},\beta)$. \bbb{(ii)} Problem (\ref{eq:one:dim}) contains one unique optimal solution.

\begin{proof}
\bbb{(i)} Dropping the bound constraint and setting the gradient of $\psi(\beta)$ to zero, we have $0=\psi'(\beta) =  ((\bar{a}\beta + \bar{b}) (\tfrac{1}{2}\bar{r}\beta^2   + \bar{s}\beta   + \bar{t}) - (\tfrac{1}{2}\bar{a} \beta^2   + \bar{b} \beta   + \bar{c}) (\bar{r}\beta + \bar{s}) )/\tau^2$. Noticing $\tau>0$, we obtain the following first-order optimal condition for $\psi$: $0=(\bar{a}\beta + \bar{b}) (\tfrac{1}{2}\bar{r}\beta^2 + \bar{s}\beta +\bar{t})  - (\tfrac{1}{2}\bar{a} \beta^2  +  \bar{b} \beta + \bar{c}) (\bar{r}\beta + \bar{s})$. It can be simplified as: $0=\tfrac{1}{2} \beta^2 \pi + \beta  \vartheta + \iota$. Solving this equation, we have two solutions $\beta_1$ and $\beta_2$. We select the one between $\Pi(\beta_1)$ and $\Pi(\beta_2)$ that leads to a lower objective value as the optimal solution. \bbb{(ii)} It can be proven by contradiction. We omit the one-sided bound constraint since it does not effect the uniqueness of the optimal solution. Assume that there exist two optimal solutions $x$ and $y$ to (\ref{lemma:line}) that lead to the same objective value $\vartheta$. According to the first-order and second-order optimal condition \cite{dinkelbach1967nonlinear,zhang2011celis}, we have: $(\bar{a} - \vartheta \bar{r})x = - \bar{b} + \vartheta  \bar{s},~(\bar{a} - \vartheta \bar{r})y = - \bar{b} + \vartheta  \bar{s},~(\bar{a} - \vartheta \bar{r})>0$, which leads to the following contradiction: $\frac{\vartheta  \bar{s} - \bar{b}}{\bar{a} - \vartheta \bar{r}} = x \neq y = \frac{\vartheta  \bar{s} - \bar{b}}{\bar{a} - \vartheta \bar{r}}$. Therefore, (\ref{eq:coordinate:one:dim}) contains one unique optimal solution. Please refer to Figure \ref{fig:shape}.
\end{proof}

\end{lemma}

\begin{figure} [!th]
\setcounter{figure}{0}
\label{fig:accuracy:cpu}
\centering
\includegraphics[width=0.15\textwidth,height=0.09\textheight]{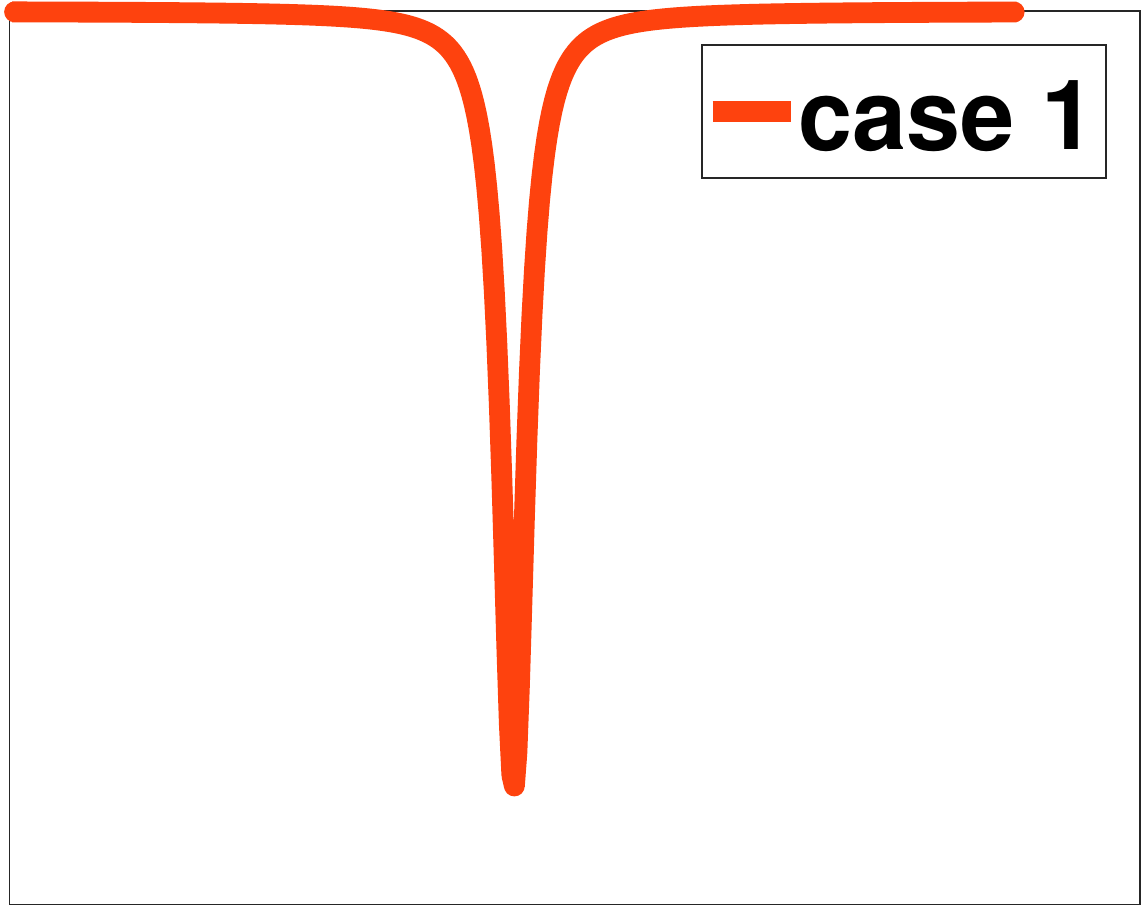}~
\includegraphics[width=0.15\textwidth,height=0.09\textheight]{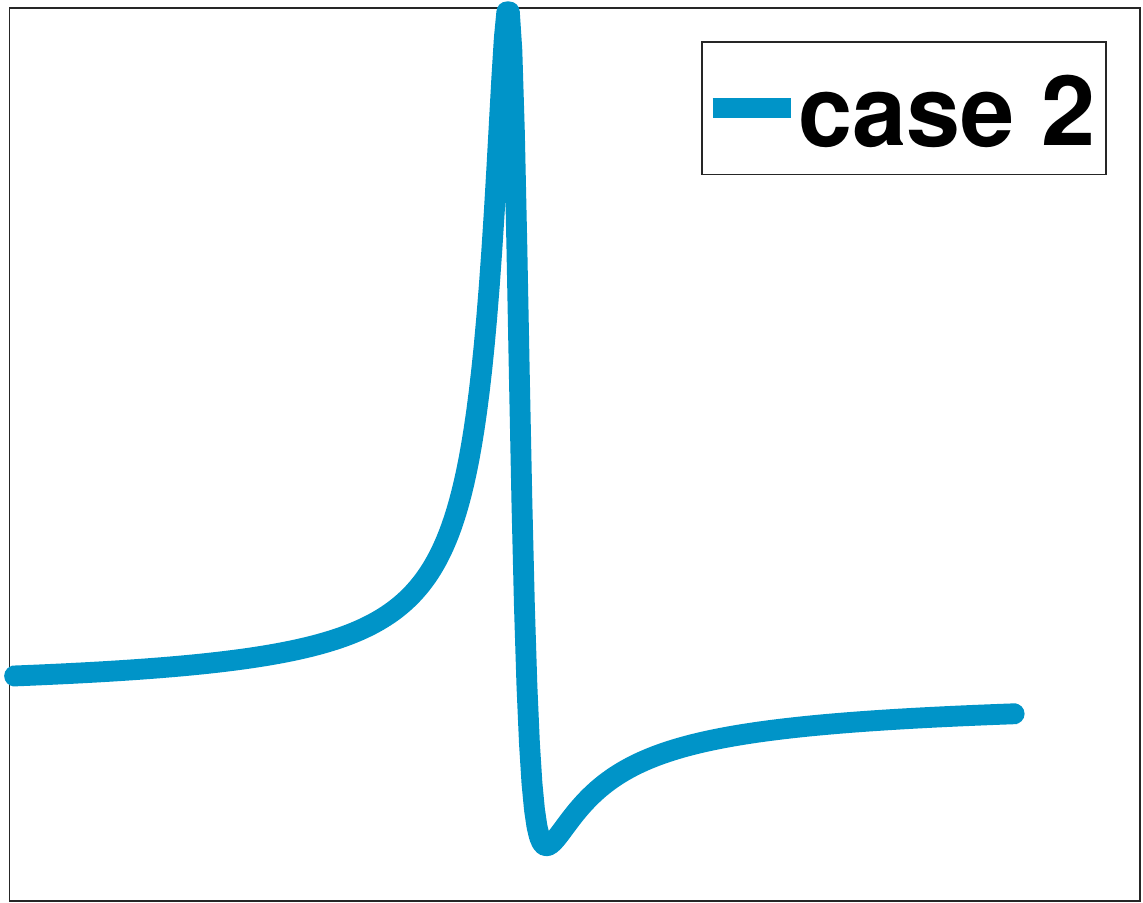}~
\includegraphics[width=0.15\textwidth,height=0.09\textheight]{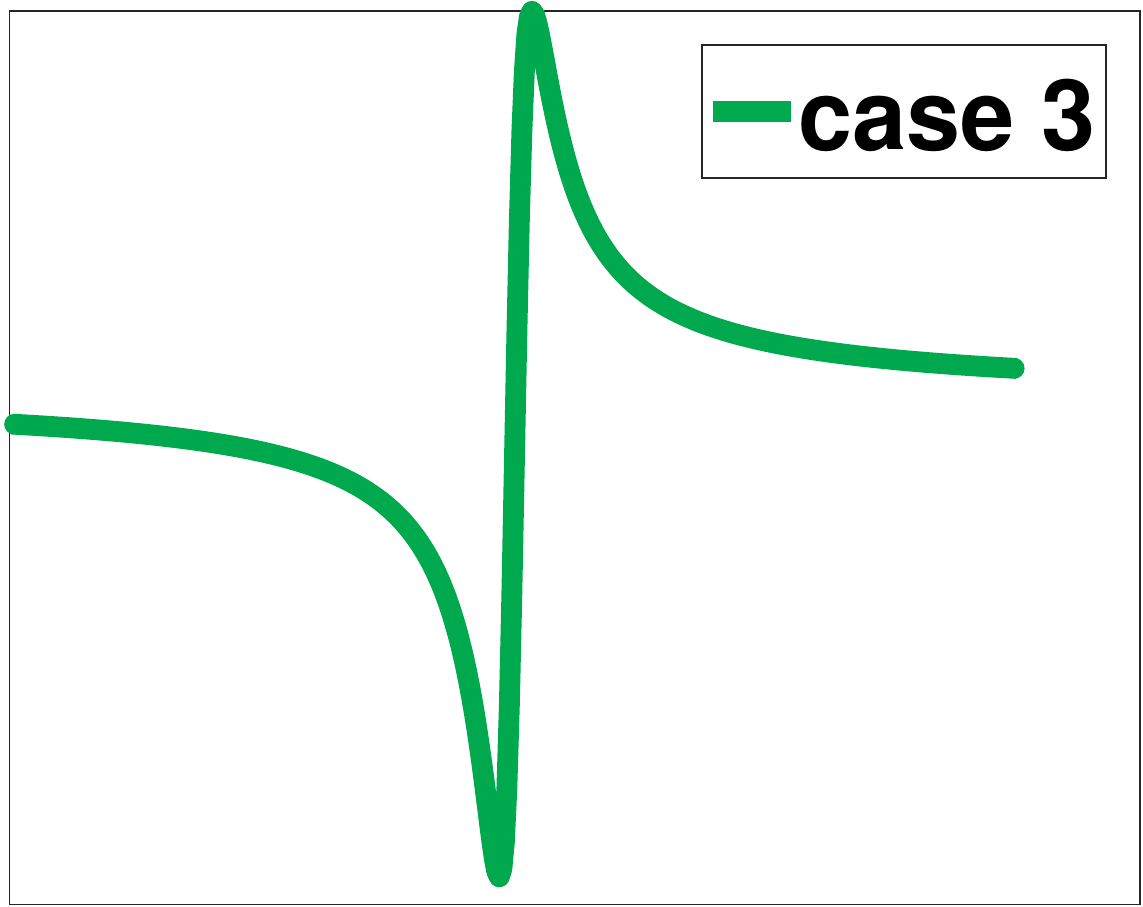}
\caption{Geometric interpretation for the one-dimensional quadratic fractional problem. Using the l'Hopital's Rule, we have $\lim_{\beta\rightarrow +\infty} \psi(\beta) = \lim_{\beta\rightarrow -\infty}~ \psi(\beta) = \frac{\bar{a}}{\bar{r}}$. Since the optimal solution is bounded and the problem at most contains two critical points, we only have the three cases above. Clearly, there exists one unique optimal solution. }\label{fig:shape}
\end{figure}


\noi Letting $\bbb{v} \triangleq \bbb{x}^t- \bbb{x}_j^t \bbb{e}_j$, we obtain: $\min_{\beta }~f(\bbb{v}+\beta \bbb{e}_i) = \min_{\beta}$~{$\tiny\frac{\tfrac{1}{2}(\bbb{v}+\beta \bbb{e}_i )^T\bbb{A}(\bbb{v}+\beta \bbb{e}_i)}{\tfrac{1}{2}(\bbb{v}+\alpha \bbb{e}_i )^T\bbb{C}(\bbb{v}+\beta \bbb{e}_i )}$}. By applying Lemma \ref{lemma:line} with $\bar{L}=-\infty,~\bar{a}=\bbb{A}_{i,i},~\bar{b}=(\bbb{Av})_i,~\bar{c}=\tfrac{1}{2}\bbb{v}^T\bbb{A}\bbb{v},~\bar{r}=\bbb{C}_{i,i},~\bar{s}=(\bbb{Cv})_i,~\bar{t}=\tfrac{1}{2}\bbb{v}^T\bbb{C}\bbb{v}$, we obtain the global optimal solution for (\ref{eq:one:dim:selection}).

\section{Solving the Subproblem} \label{sect:subp}
The subproblem (\ref{eq:subprob}) in Algorithm \ref{algo:main} reduces to the following quadratic fractional programming problem:
\beq \label{eq:subprob:k}
 \bbb{z}^* = \arg \min_{\|\bbb{z}\|_0 \leq q}~p(\bbb{z})\triangleq \frac{\tfrac{1}{2}\bbb{z}^T\bar{\bbb{Q}}\bbb{z}+\bar{\bbb{p}}^T\bbb{z}+\bar{w}}{\tfrac{1}{2}\bbb{z}^T\bar{\bbb{R}}\bbb{z}+\bar{\bbb{c}}^T\bbb{z}+\bar{v}},
\eeq
\noi where $\bbb{z}\in \mathbb{R}^k$, $\bar{\bbb{Q}} = \bbb{A}_{BB}+\theta \bbb{I}$,~$\bar{\bbb{p}}=\bbb{A}_{BN}\bbb{x}_{N}-\theta \bbb{x}^t_B$, $\bar{w}=\tfrac{1}{2} \bbb{x}_N^T \bbb{A}_{NN} \bbb{x}_N + \tfrac{\theta}{2}\|\bbb{x}_B^t\|_2^2$, $\bar{\bbb{R}} = \bbb{C}_{BB}$,~$\bar{\bbb{c}}=\bbb{C}_{BN}\bbb{x}_{N}$,~,~$\bar{v}=\frac{1}{2} \bbb{x}_N^T \bbb{C}_{NN} \bbb{x}_N,~ q=s-\|\bbb{x}_N\|_0$.

Problem (\ref{eq:subprob:k}) is equally NP-hard due to the combinatorial constraint $\|\bbb{z}\|_0\leq q$. Inspired by the work of \cite{yuan2017hybrid}, we develop an exhaustive tree/combinatorial search algorithm to solve it. Specifically, we consider to solve the following optimization problem: $\min_{\bbb{z}\in \mathbb{R}^k}~p(\bbb{z}), ~s.t.~\bbb{z}_K=0$, where $K$ has $\sum_{i=0}^q C_k^{i}$ possible choices for the coordinates. We systematically enumerate the full binary tree for $K$ to obtain all possible candidate solutions for $\bbb{z}$ and then pick the best one that leads to the lowest objective value as the optimal solution. In other words, we need to solve the following quadratic fractional programming problem with $m \triangleq k-|K|$ variables:
\beq \label{eq:quad:fraction:subproblem}
 \bbb{y}^* = \arg \min_{\bbb{y}}~\mathcal{L}(\bbb{y}) \triangleq \frac{{u}(\bbb{y})}{{q}(\bbb{y})}\triangleq\frac{\frac{1}{2}\bbb{y}^T\bbb{Q}\bbb{y}+\bbb{p}^T\bbb{y}+w}{\frac{1}{2}\bbb{y}^T\bbb{R}\bbb{y}+\bbb{c}^T\bbb{y}+v},
\eeq
\noi where $\bbb{y}\in \mathbb{R}^m$. The optimal solution of (\ref{eq:subprob:k}) can be computed as $\bbb{z}^*_K=\bbb{0}$, $\bbb{z}^*_{\bar{\bbb{K}}}=\bbb{y}^*$ with $\bar{\bbb{K}}=\{1,2,...,k\}\setminus K$. Therefore, if we find the global optimal solution of (\ref{eq:quad:fraction:subproblem}), we find the global optimal solution of (\ref{eq:subprob:k}) as well.

The non-convex problem in (\ref{eq:quad:fraction:subproblem}) is still challenging. {For solving it, we present two methods, namely a bisection search method and a coordinate descent method, which are of independent research interest.}



\subsection{A Bisection Search Method}
This subsection presents a bisection search method for finding the global optimal solution of Problem (\ref{eq:quad:fraction:subproblem}).

We now discuss the relationship between this fractional programming problem and the following parametric programming problem \cite{dinkelbach1967nonlinear}:
\beq \label{eq:parametric}
 \mathcal{J}(\alpha)=0,~\text{with}~\mathcal{J}(\alpha)  \triangleq \min_{\bbb{y}} ~u(\bbb{y})- \alpha q(\bbb{y})\nn
\eeq
\noi This is a feasibility problem with respect to $\alpha$: $\mathcal{J}(\alpha)=u(y^*(\alpha))-\alpha q(y^*(\alpha))=0$, where $y^*(\alpha)\in \mathbb{R}^m$ is defined as $y^*(\alpha) \triangleq \arg\min_{\bbb{y}}~u(\bbb{y})-\alpha q(\bbb{y})$.

The following theorem sheds some theoretic lights for the original non-convex problem in (\ref{eq:quad:fraction:subproblem}).

%

%
%
%

\begin{theorem} \label{theorem:bisection}
We have the following results. \bbb{(i)} It holds that: $\lambda_{\min}\left(  \bbb{Z}\right) \leq \min_{\bbb{y}}~\mathcal{L}(\bbb{y}) < \lambda_{\min}\left(  \bbb{O}\right)$, with $\bbb{O} \triangleq \bbb{R}^{-1/2} \bbb{Q} \bbb{R}^{-1/2}$, $\gamma \triangleq 2v-\|\bbb{R}^{-1/2}\bbb{c}\|_2^2>0$, $\bbb{g} \triangleq \bbb{R}^{-1/2}\bbb{p}-\bbb{R}^{-1/2} \bbb{Q} \bbb{R}^{-1}\bbb{c}$, $\delta \triangleq \bbb{c}^T\bbb{R}^{-1}\bbb{Q}\bbb{R}^{-1}\bbb{c}-2\bbb{c}^T\bbb{R}^{-1}\bbb{p}+2w$, and $ \bbb{Z}\triangleq{\scriptsize \begin{pmatrix}
                                                 \bbb{O} & \bbb{g}/\sqrt{\gamma} \\
                                                 \bbb{g}^T/\sqrt{\gamma} & {\delta}/{\gamma} \\
                                               \end{pmatrix}}$. \bbb{(ii)} Let $\bbb{O}=\bbb{U}diag(\bbb{d})\bbb{U}^T$ be the eigenvalue decomposition of $\bbb{O}$. The function $\mathcal{J}(\alpha)$ can be rewritten as
\beq \label{eq:linesearch}
 \mathcal{J}(\alpha)= \tfrac{1}{2}\delta - \tfrac{1}{2}\alpha \gamma -\tfrac{1}{2}\sum_{i}^m\frac{\bbb{a}^2_i}{\bbb{d}_i-\alpha},~\text{with}~\bbb{a}=\bbb{U}^T\bbb{g}
\eeq
\noi  and it is monotonically decreasing on the range $\lambda_{\min}(\bbb{Z})\leq \alpha < \lambda_{\min}(\bbb{O})$. The optimal solution can be computed as $\bbb{y}^* = \bbb{R}^{-1/2} (\bbb{u}^* - \bbb{R}^{-1/2}\bbb{c})$, with $\bbb{u}^*=- (\bbb{O}-\alpha^* \bbb{I})^{-1} \bbb{g}$ and $\alpha^*$ being the unique root of the equation $\mathcal{J}(\alpha)=0$ on the range $\lambda_{\min}(\bbb{Z})\leq \alpha < \lambda_{\min}(\bbb{O})$.

\begin{proof}
\bbb{(i)} Firstly, it is not hard to notice that Program (\ref{eq:quad:fraction:subproblem}) is equivalent to the following problem:
{
\begin{align} \label{eq:nonconvex:simple}
 &\min_{\bbb{y}}~\mathcal{L}(\bbb{y})\nn\\
 =~& \min_{\bbb{d}}~\frac{\frac{1}{2}(\bbb{R}^{-1/2}\bbb{d})^T\bbb{Q}(\bbb{R}^{-1/2}\bbb{d})+\bbb{p}^T(\bbb{R}^{-1/2}\bbb{d})+w}{\frac{1}{2}(\bbb{R}^{-1/2}\bbb{d})^T\bbb{R}(\bbb{R}^{-1/2}\bbb{d})+\bbb{c}^T(\bbb{R}^{-1/2}\bbb{d})+v} \nn\\
=~&\min_{\bbb{d}}~\frac{\frac{1}{2} \bbb{d}^T \bbb{O} \bbb{d}+\bbb{d}^T(\bbb{R}^{-1/2}\bbb{p})+w}{\frac{1}{2}\|\bbb{d}\|_2^2+\bbb{d}^T(\bbb{R}^{-1/2}\bbb{c})+v} \nn\\
=~&\min_{\bbb{d}}~\frac{\frac{1}{2} \bbb{d}^T \bbb{O} \bbb{d}+\bbb{d}^T(\bbb{R}^{-1/2}\bbb{p})+w}{\frac{1}{2}\|\bbb{d}+\bbb{R}^{-1/2}\bbb{c}\|_2^2+v-\frac{1}{2}\|\bbb{R}^{-1/2}\bbb{c}\|_2^2} \nn\\
=~&\min_{\bbb{u}}~\frac{\frac{1}{2}\bbb{u}^T\bbb{O}\bbb{u} + \bbb{u}^T \bbb{g} + \frac{1}{2}\delta }{\frac{1}{2}\|\bbb{u}\|_2^2+\frac{1}{2}\gamma},
\end{align}}{\normalsize}where the first step uses the variable substitution that $\bbb{y}=\bbb{R}^{-1/2}\bbb{d}$; the third step uses the transformation that $\bbb{u}=\bbb{d}+\bbb{R}^{-1/2}\bbb{c}$. We notice that the denominator is always strictly positive for all decision variables. Letting $\bbb{u}=\bbb{0}$ in (\ref{eq:nonconvex:simple}), we obtain $\frac{1}{2}\gamma>0$. We naturally obtain the upper bound for $\min_{\bbb{y}}~\mathcal{L}(\bbb{y})$:
{
\begin{align*}
\min_{\bbb{y}}~\mathcal{L}(\bbb{y})=~&\textstyle\min_{\bbb{u},\eta=\sqrt{\gamma}}~\frac{\frac{1}{2}\bbb{u}^T\bbb{O}\bbb{u} + \frac{1}{\sqrt{\gamma}}\bbb{u}^T \bbb{g}\eta + \frac{\delta}{2\gamma} \eta^2 }{\frac{1}{2}\|\bbb{u}\|_2^2+\frac{1}{2}\eta^2}~~~~~~~~~~~~\nn\\
\geq~&\min_{\bbb{u},\eta}~\frac{\frac{1}{2}\bbb{u}^T\bbb{O}\bbb{u} + \frac{1}{\sqrt{\gamma}}\bbb{u}^T \bbb{g}\eta + \frac{\delta}{2\gamma} \eta^2 }{\frac{1}{2}\|\bbb{u}\|_2^2+\frac{1}{2}\eta^2}~~~~~~~~~~~~~~~~\nn\\
=~&\min_{\bbb{u},\eta}~\frac{\frac{1}{2}\left[\bbb{u}^T ~|~\eta^T\right]^T\bbb{Z}~\left[\bbb{u}^T ~|~\eta^T\right]}{\frac{1}{2}\|\bbb{u}\|_2^2+\frac{1}{2}\eta^2} \nn\\
=~& \lambda_{\min}(\bbb{Z}),\nn
\end{align*}}{\normalsize}\noi where the first inequality uses the fact that $\min_{\bbb{x}} f(\bbb{x}) \leq \min_{\bbb{x}\in\Psi } f(\bbb{x})$ for all $f(\cdot)$ and $\Psi$.

We now derive the upper bound of $\min_{\bbb{y}}~\mathcal{L}(\bbb{y})$. Since the objective function $\mathcal{J}(\alpha)$ is always bounded, there must exist $\alpha$ with $\bbb{Q}-\alpha \bbb{R}\succ 0$, such that the value of $\bbb{y}$ minimizing the function $\left(u(\bbb{y})- \alpha q(\bbb{y})\right)$. Therefore, we have $\bbb{Q}-\alpha \bbb{R}\succ 0 \Rightarrow \bbb{Q}-\alpha \bbb{R}^{1/2}\bbb{I}\bbb{R}^{1/2}\succ 0 \Rightarrow \bbb{R}^{-1/2}\bbb{Q}\bbb{R}^{-1/2}-\alpha \bbb{I}\succ 0 \Rightarrow \alpha< \lambda_{\min}(\bbb{O})$.

\bbb{(ii)} Using the result of (\ref{eq:nonconvex:simple}), we can rewrite $\mathcal{J}(\alpha)$ as:
\begin{align*}
 \mathcal{J}(\alpha) &=  \min_{\bbb{u}} \tfrac{1}{2}\bbb{u}^T\bbb{O}\bbb{u} + \bbb{u}^T \bbb{g} + \tfrac{1}{2}\delta - \alpha \left(\tfrac{1}{2}\|\bbb{u}\|_2^2+\tfrac{1}{2}\gamma\right)\nn\\
&=   \min_{\bbb{u}} \tfrac{1}{2}\bbb{u}^T(\bbb{O}-\alpha I)\bbb{u} + \bbb{u}^T \bbb{g} + \tfrac{1}{2}\delta -  \tfrac{\alpha\gamma}{2}. \nn
\end{align*}
\noi Solving the quadratic optimization with respect to $\bbb{u}$ we have $\bbb{u}^*=-(\bbb{O}-\alpha \bbb{I})\bbb{g}$. Thus, we can repress $\mathcal{J}(\alpha)$ as: $\mathcal{J}(\alpha)= \textstyle -\tfrac{1}{2}\bbb{g}^T (\bbb{O}-\alpha I)^{-1}\bbb{g}+ \tfrac{1}{2}\delta -  \tfrac{\alpha\gamma}{2}$. Since it holds that $\bbb{g}^T (\bbb{O}-\alpha I)^{-1}\bbb{g} = \bbb{g}^T \bbb{U}^T \text{diag}(1\div(\bbb{d}-\alpha))\bbb{Ug}$ with $\div$ denoting the element-wise division between two vectors, we obtain (\ref{eq:linesearch}). Noticing that the first-order and second-order gradient of $\mathcal{J}(\alpha)$ with respect to $\alpha$ can be computed as: $\mathcal{J}'(\alpha) = -\tfrac{1}{2} \sum_{i}^m(\frac{\bbb{a}_i}{\bbb{d}_i-\alpha})^2 -  \tfrac{\gamma}{2} ,~\mathcal{J}''(\alpha) = - \sum_{i}^m ({\bbb{a}^2_i}/(\bbb{d}_i-\alpha)^3) $ and $\gamma>0$, we obtain $\mathcal{J}'(\alpha)<0$ and $\mathcal{J}''(\alpha)\leq 0$. Thus, the function $\mathcal{J}(\alpha)$ is concave and monotonically decreasing on the range $\lambda_{\min}(\bbb{Z})\leq \alpha < \lambda_{\min}(\bbb{O})$, and there exists a unique root of the equation $\mathcal{J}(\alpha)=0$ on the range $\lambda_{\min}(\bbb{Z})\leq \alpha < \lambda_{\min}(\bbb{O})$.

\end{proof}

\end{theorem}

Based on Theorem \ref{theorem:bisection}, we now present a bisection method for solving Problem (\ref{eq:quad:fraction:subproblem}). For notation convenience, we define $\underline{\alpha} \triangleq \lambda_{\min}(\bbb{Z})$ and $\overline{\alpha} \triangleq \lambda_{\max}(\bbb{O})-\epsilon$, where $\epsilon$ denotes the machine precision parameter which is sufficiently small. Due to the monotonically decreasing property of $\mathcal{J}(\alpha)$, we can solve (\ref{eq:linesearch}) by checking where the sign of the left-hand side changes. Specifically, we consider the following three cases for $\mathcal{J}(\alpha)$ on the range $\underline{\alpha} \leq \alpha \leq \overline{\alpha}$: \bbb{(a)} $\mathcal{J}(\underline{\alpha})\geq \mathcal{J}(\overline{\alpha})\geq 0$, \bbb{(b)} $0\geq \mathcal{J}(\underline{\alpha}) \geq \mathcal{J}(\overline{\alpha})$, and \bbb{(c)} $\mathcal{J}(\underline{\alpha}) \geq 0 \geq\mathcal{J}(\overline{\alpha})$. For case $(\bbb{a})$ and $(\bbb{b})$, we can directly return $\overline{\alpha}$ and $\underline{\alpha}$ as the optimal solution, respectively. We now consider case $(\bbb{c})$. By the Rolle mean value theorem, there always exists an $\alpha^* \in [\underline{\alpha},~\overline{\alpha}]$ such that $\mathcal{J}(\alpha^*)=0$. Thus, we can define and initialize the lower bound $lb = \underline{\alpha}$ and the upper bound $ub = \overline{\alpha}$. We then perform the following loop until the optimal solution $\alpha^*=mid$ with $\mathcal{J}(mid)\thickapprox 0$ is found: $\{mid = (lb+ub)/2,~\text{if}(\mathcal{J}(mid)>0)~lb = mid;~\text{else}~ub = mid;\}$. Such a bisection scheme is guaranteed to find the optimal solution within $\mathcal{O}(\log_2((\overline{\alpha}-\underline{\alpha})/\varepsilon))$ iterations that $ub \leq lb + \varepsilon$ \cite{boyd2004convex}.

\noi\bbb{Remarks.} \bbb{(i)} To our knowledge, this is the first algorithm for unconstrained quadratic fractional programming with global optimal guarantee. The work of \cite{GuoLYSW03} also discusses a bisection search method for the ratio of trace problem, but it can not solve our general quadratic fractional programming problem. The classical Dinkelbach's method \cite{dinkelbach1967nonlinear,Wang07traceratio} can solve our problem, but it only finds a stationary solution for the non-convex problem. Our results are based on the monotone property of the associated parametric programming problem in a restricted domain. \bbb{(ii)} The unconstrained fractional quadratic program can be solved to optimality by linear semidefinite programming and it is related to the S-lemma for the quadratically constrained quadratic program \cite{Polik2007}. In this paper, we show that it can be solved using a bisection search method. This method has the merit that it is simple and easy to implement. In addition, it is efficient and it does not require iterative eigenvalue decomposition as in the semidefinite programming lifting techniques. \bbb{(iii)} The matrix $\bbb{O}$ is a $n\times n$ principal sub-matrix of $\bbb{Z}$. Using Theorem 4.3.17 in \cite{horn1990matrix}, it always holds that $\bbb{\lambda}_1(\bbb{Z}) \leq \bbb{\lambda}_1(\bbb{O}) \leq \bbb{\lambda}_2(\bbb{Z}) \leq ... \leq \bbb{\lambda}_{n-1}(\bbb{Z}) \leq \bbb{\lambda}_{n-1}(\bbb{O}) \leq \bbb{\lambda}_{n}(\bbb{Z})$, where $\bbb{\lambda}(\bbb{X})$ denotes the eigenvalues of $\bbb{X}$ in increasing order. Thus, the bound for the $\alpha^*$ is tight.

\subsection{A Coordinate Descent Method}

This subsection presents a simple coordinate descent method \cite{Tseng2001,HsiehCLKS08,HsiehD11,VandaeleGLZD16,LeiZD16} for solving Problem (\ref{eq:quad:fraction:subproblem}). Although it can not guarantee to find the global optimal solution, it has many merits. \bbb{(i)} It is able to incorporate additional bound constraints. \bbb{(ii)} It is numerically robust and does not require additional eigenvalue solvers. \bbb{(iii)} It is guaranteed to converge to a coordinate-wise minimum point for our specific problem (see Proposition \ref{theorem:convergence:coo} below).

To illustrate the merits of the coordinate descent method, we consider incorporating the bound constraint $\bbb{x}\geq \bar{L}$ on the solution for Problem (\ref{eq:main}) \footnote{This is useful in sparse non-negative PCA \cite{AsterisPD14}.}. Our decomposition algorithm for finding the working set and strategies for handling the NP-hard $\ell_0$ norm directly follow and what one needs is to replace (\ref{eq:quad:fraction:subproblem}) and solve the following problem:
\beq \label{eq:coordinate:bound}
 \min_{\bbb{y} \in \mathbb{R}^m}~\mathcal{L}(\bbb{y})\triangleq\frac{\frac{1}{2}\bbb{y}^T\bbb{Q}\bbb{y}+\bbb{p}^T\bbb{y}+w}{\frac{1}{2}\bbb{y}^T\bbb{R}\bbb{y}+\bbb{c}^T\bbb{y}+v},~s.t.~\bbb{y}\geq \hat{L}
\eeq
 \noi for some constant $\hat{L}$.

 The coordinate descent method iteratively picks a coordinate $i\in\{1,2,...,m\}$ and solves the following one dimensional subproblem based on its current solution $\bbb{y}^j$ with $j=0,1,...\infty$:
\beq \label{eq:coordinate:one:dim}
 \beta^* = \arg \min_{ \beta}~\mathcal{L}(\bbb{y}^j+ \beta \bbb{e}_i),~s.t.~\bbb{y}^j_i+\beta\geq \hat{L}
\eeq
\noi where $j$ is the iteration counter for the coordinate descent algorithm. Problem (\ref{eq:coordinate:one:dim}) reduces to the one-dimensional subproblem as in Lemma \ref{lemma:line} with suitable parameters. In every iteration $j$, once the optimal solution $\beta^*$ in (\ref{eq:coordinate:one:dim}) is found, the intermediate solution for (\ref{eq:coordinate:bound}) is updated via $\bbb{y}_i^{j+1} \Leftarrow \bbb{y}_i^{j} + \beta^*$. There are several ways and orders to decide which coordinates to update in the literature. \bbb{(i)} Cyclic order strategy runs all coordinates in cyclic order, i.e., $1 \rightarrow 2 \rightarrow ...  \rightarrow m\rightarrow 1$. \bbb{(ii)} Random sampling strategy randomly selects one coordinate to update (sample with replacement). \bbb{(iii)} Gauss-Southwell strategy picks coordinate $i$ such that $i=\arg \max_{1\leq t \leq m} |\bar{\bbb{\nabla}} \mathcal{L}(\bbb{x}^j)|_t$, with $\bar{\bbb{\nabla}} \mathcal{L}(\bbb{x}) \in \mathbb{R}^m$ being the projected gradient of $\mathcal{L}$ at $\bbb{x}$ \cite{lin2007projected}: $\tiny
\bar{\bbb{\nabla}} \mathcal{L}(\bbb{x})_i = \left\{
  \begin{array}{ll}
    \bbb{\nabla} \mathcal{L}(\bbb{x})_i, & \hbox{$\bbb{x}_i>\hat{L}$;} \\
    \min(0,\bbb{\nabla} \mathcal{L}(\bbb{x})_i), & \hbox{$\bbb{x}_i=\hat{L}$;}
  \end{array}
\right.$ and $\bbb{\nabla} \mathcal{L}(\bbb{x})$ being the gradient of $\bbb{\nabla} \mathcal{L}$ at $\bbb{x}$. Note that $\bar{\bbb{\nabla}} \mathcal{L}(\grave{\bbb{x}})=\bbb{0}$ implies $\grave{\bbb{x}}$ is a first-order stationary point.


We now present our convergence result of the coordinate descent method for solving (\ref{eq:coordinate:bound}), which is an extension of Theorem 4.1 in \cite{Tseng2001}.  Some proofs can be found in the \textbf{Appendix}.

\begin{proposition} \label{theorem:convergence:coo}
When the cyclic order strategy is used, coordinate descent method is guaranteed to converge to a coordinate-wise minimum of Problem (\ref{eq:coordinate:bound}) that $\forall i,~\bbb{y}_i^* = \arg \min_{\alpha\geq \hat{L}} ~\mathcal{L}(\bbb{y}_i^*+\alpha \bbb{e}_i)$.

\end{proposition}
\noi\bbb{Remarks.} \bbb{(i)} Convergence of the coordinate descent method requires a unique solution in each minimization step; otherwise, it may cycle indefinitely. A simple but intriguing example is given in \cite{Powell1973}. One good feature of the non-convex problem in (\ref{eq:coordinate:bound}) is that its associated one-dimensional subproblem in (\ref{eq:coordinate:one:dim}) only contains one unique optimal solution (see part (ii) in Lemma \ref{lemma:line}). This is different from the work of \cite{VandaeleGLZD16} where their one-dimensional subproblem may have multiple optimal solutions and cause divergence. \bbb{(ii)} Coordinate descent method is guaranteed to produce a coordinate-wise stationary point which is stronger than the full gradient projection method. Note that any coordinate-wise stationary point $\bbb{x}^*$ that $\forall i,~\bbb{x}_i^* = \arg \min_{\alpha\geq \hat{L}} ~\mathcal{L}(\bbb{x}_i^*+\alpha \bbb{e}_i)$ also satisfies the first-order optimal condition with $\bar{\bbb{\nabla}} \mathcal{L}(\bbb{x}^*)=0$. However, the reverse is not true. This implies that the coordinate descent method can exploit possible higher order derivatives to escape saddle points for the non-convex problem.



\section{Convergence Analysis of Algorithm \ref{algo:main} }\label{sect:conv}

This section presents the convergence analysis of Algorithm \ref{algo:main}. We assume that $\{f(\bbb{x}^t)\}_{t=0}^{\infty}$ is generated by Algorithm \ref{algo:main} and the solution is bounded with $0<\|\bbb{x}^t\|<\infty$ for all $t$ throughout this section. We first present the following lemma.

\begin{lemma} \label{lemma:suf:dec}
(\textbf{Sufficient Decrease Condition}) It holds that: $f(\bbb{x}^{t+1}) - f(\bbb{x}^t) \leq  \frac{-\theta \|\bbb{x}^{t+1}-\bbb{x}^t\|_2^2}{(\bbb{x}^{t+1})^T\bbb{Cx}^{t+1}}.$

\end{lemma}

\noi \bbb{Remarks.} The proximal term in the numerator in (\ref{eq:subprob}) is necessary for our non-convex problem since it guarantees sufficient decrease condition which is important for convergence.

\begin{figure*} [!th]
\setcounter{figure}{0}
\label{fig:accuracy:cpu}
\captionsetup{singlelinecheck = on, format= hang, justification=justified, font=footnotesize, labelsep=space}
\centering
      \begin{subfigure}{\threefigwid}\includegraphics[width=\textwidth,height=\objimghei]{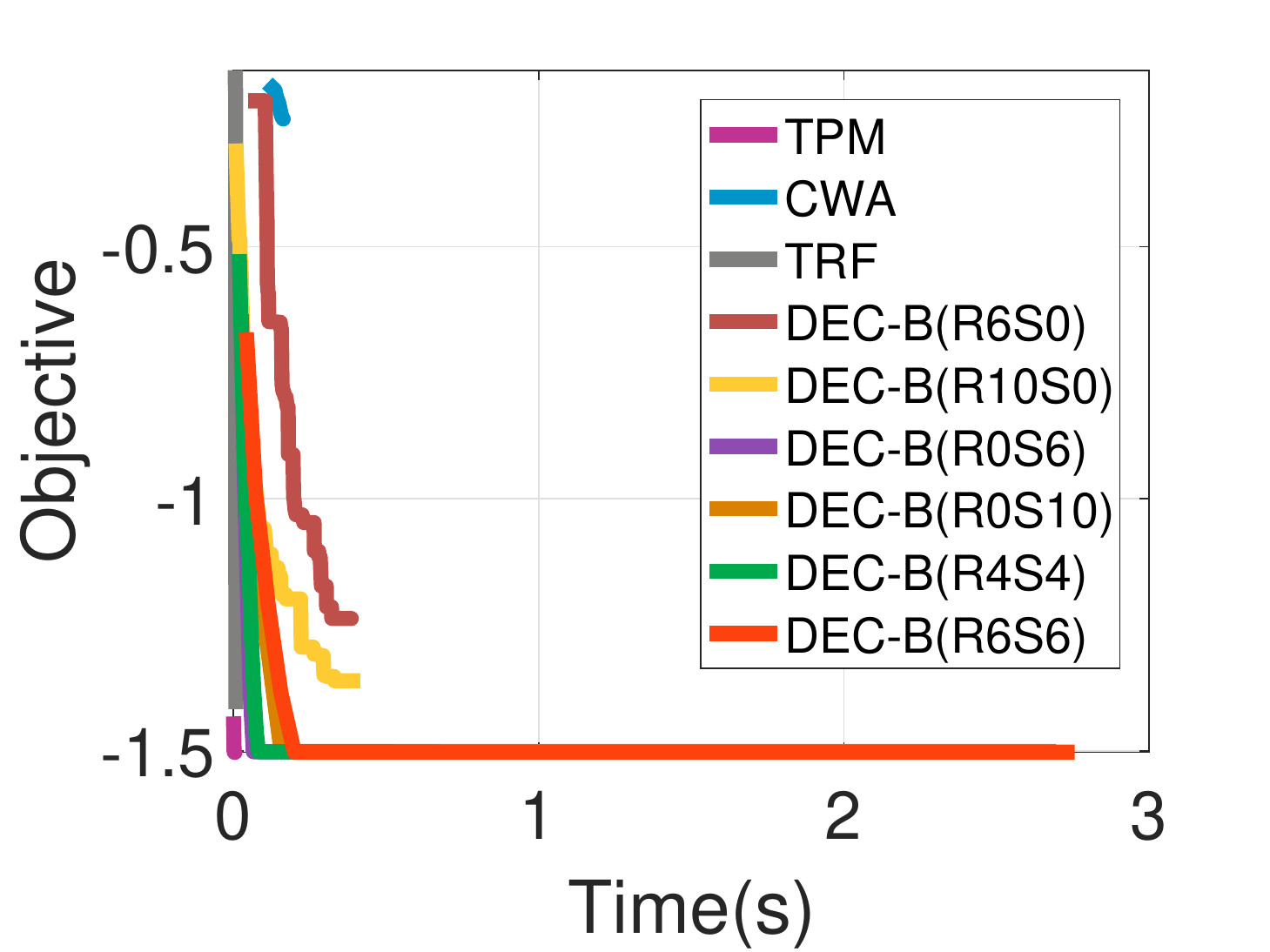}\vspace{-6pt}\caption{\footnotesize sparse PCA, `w1a', s=15}\end{subfigure}\ghs
      \begin{subfigure}{\threefigwid}\includegraphics[width=\textwidth,height=\objimghei]{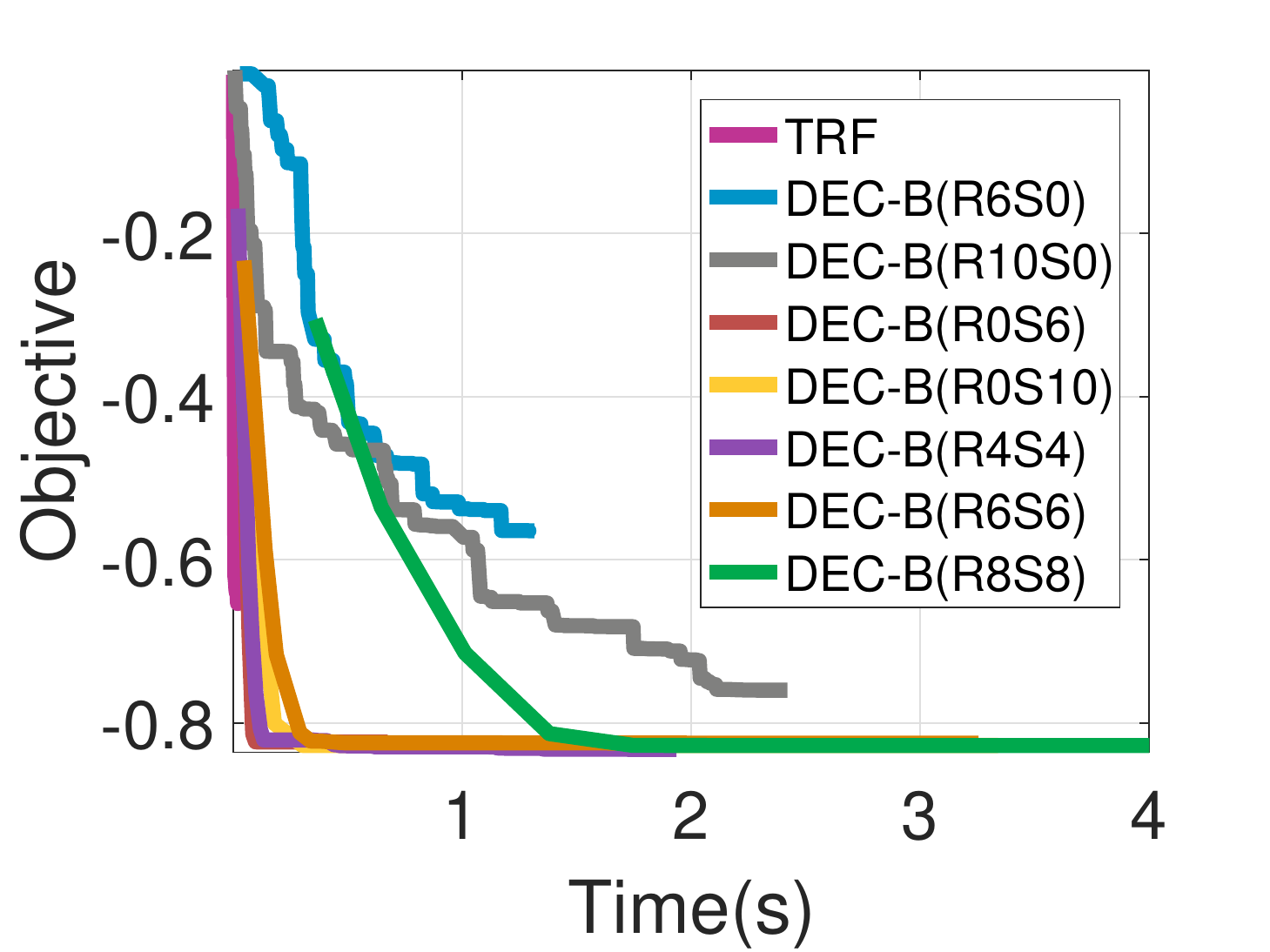}\vspace{-6pt}\caption{\footnotesize sparse FDA, `randn-500', s=15}\end{subfigure}\ghs
      \begin{subfigure}{\threefigwid}\includegraphics[width=\textwidth,height=\objimghei]{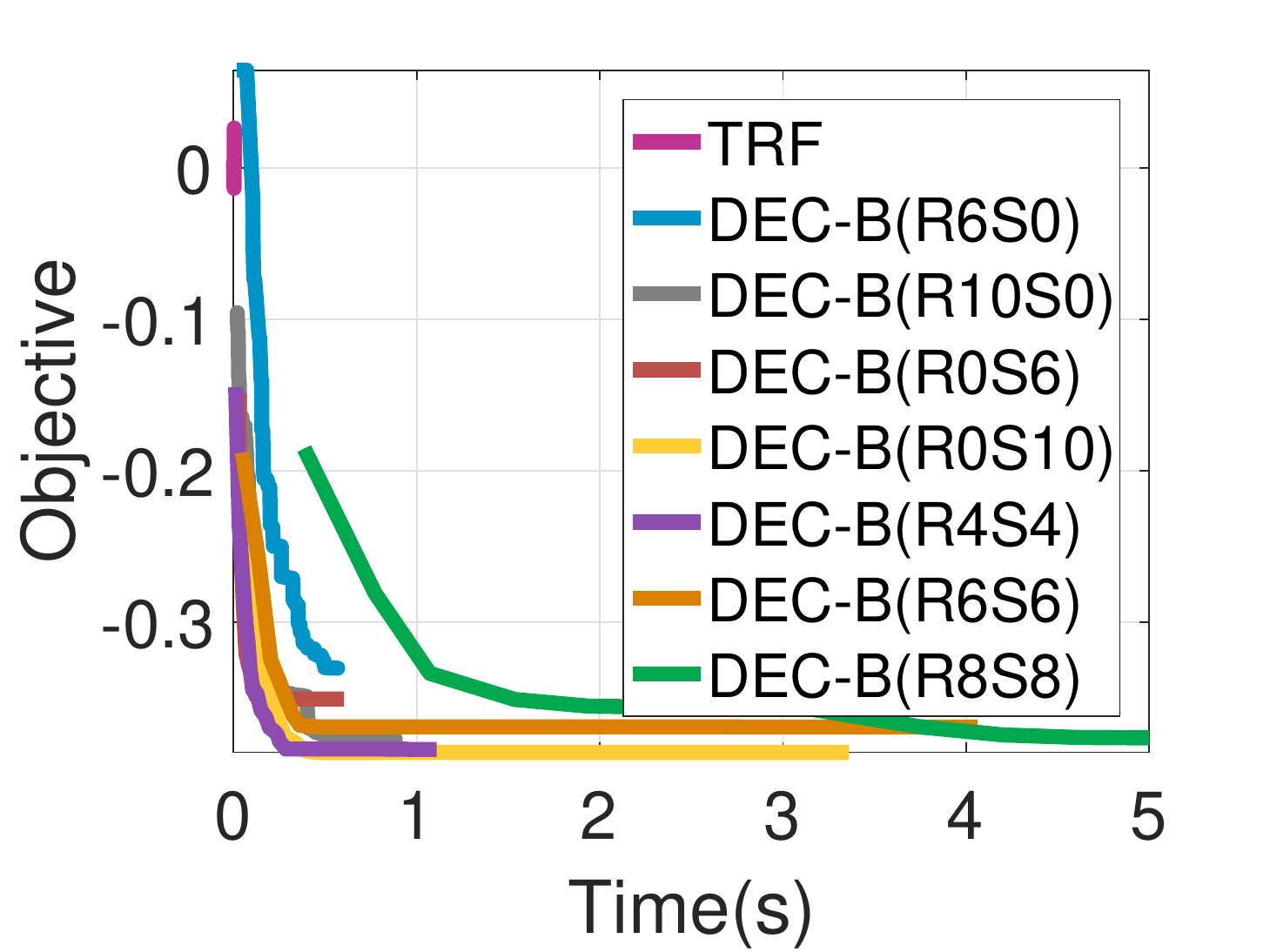}\vspace{-6pt}\caption{\footnotesize sparse CCA, `randn-500', s=15}\end{subfigure}
      \centering
      \begin{subfigure}{\threefigwid}\includegraphics[width=\textwidth,height=\objimghei]{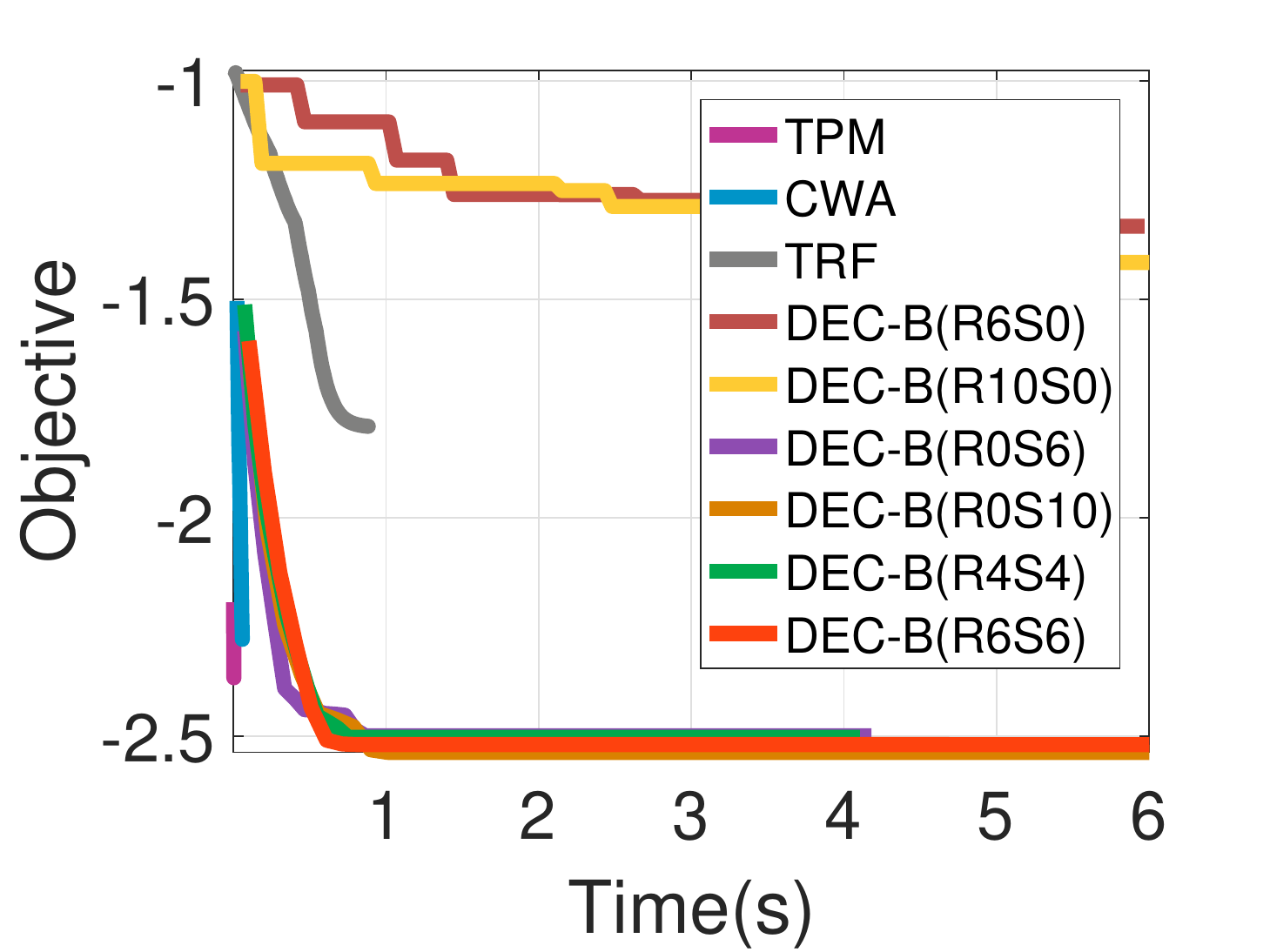}\vspace{-6pt}\caption{\footnotesize sparse PCA, `randn-2000', s=15}\end{subfigure}\ghs
      \begin{subfigure}{\threefigwid}\includegraphics[width=\textwidth,height=\objimghei]{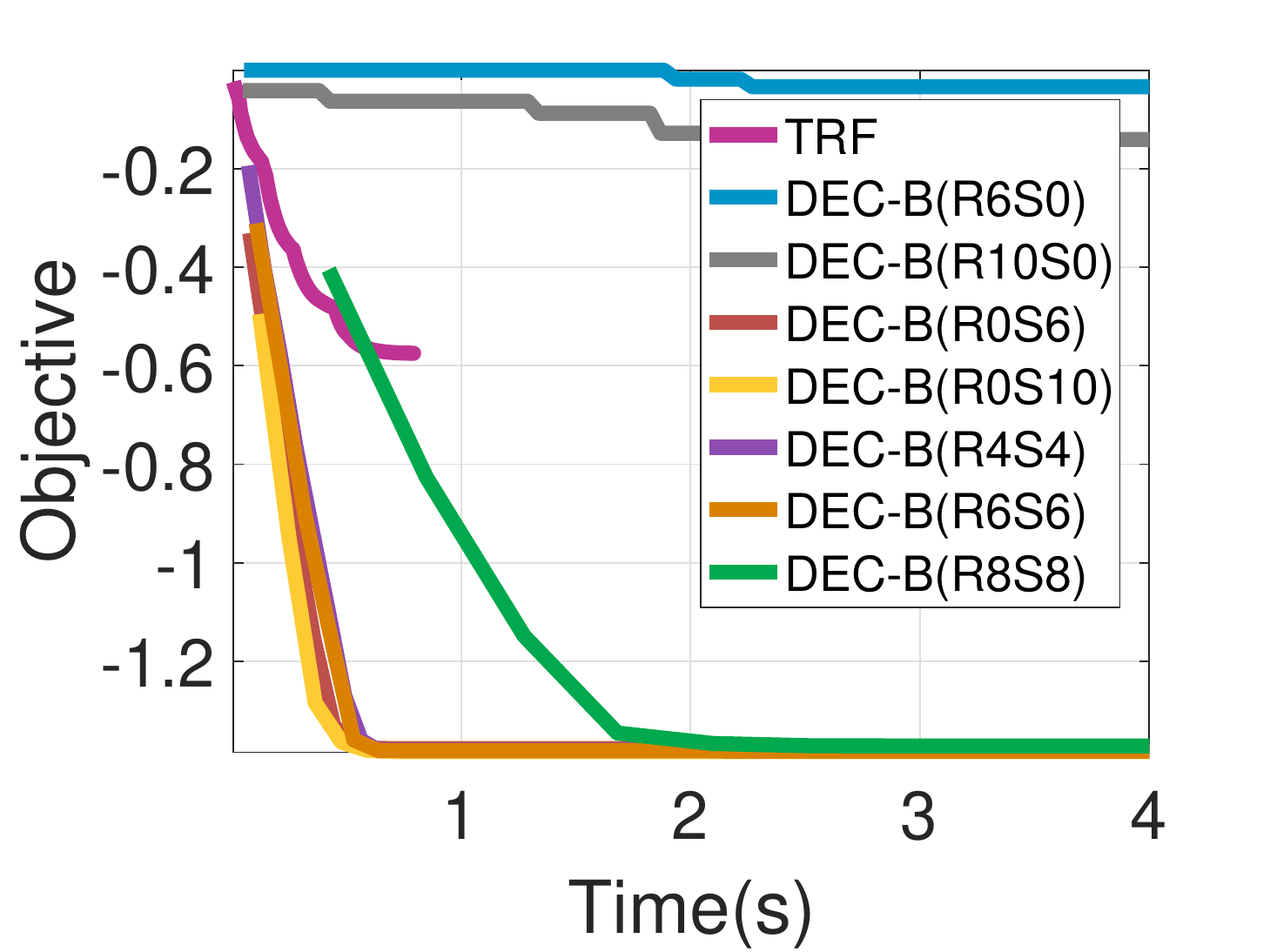}\vspace{-6pt}\caption{\footnotesize sparse FDA, `randn-2000', s=15}\end{subfigure}\ghs
      \begin{subfigure}{\threefigwid}\includegraphics[width=\textwidth,height=\objimghei]{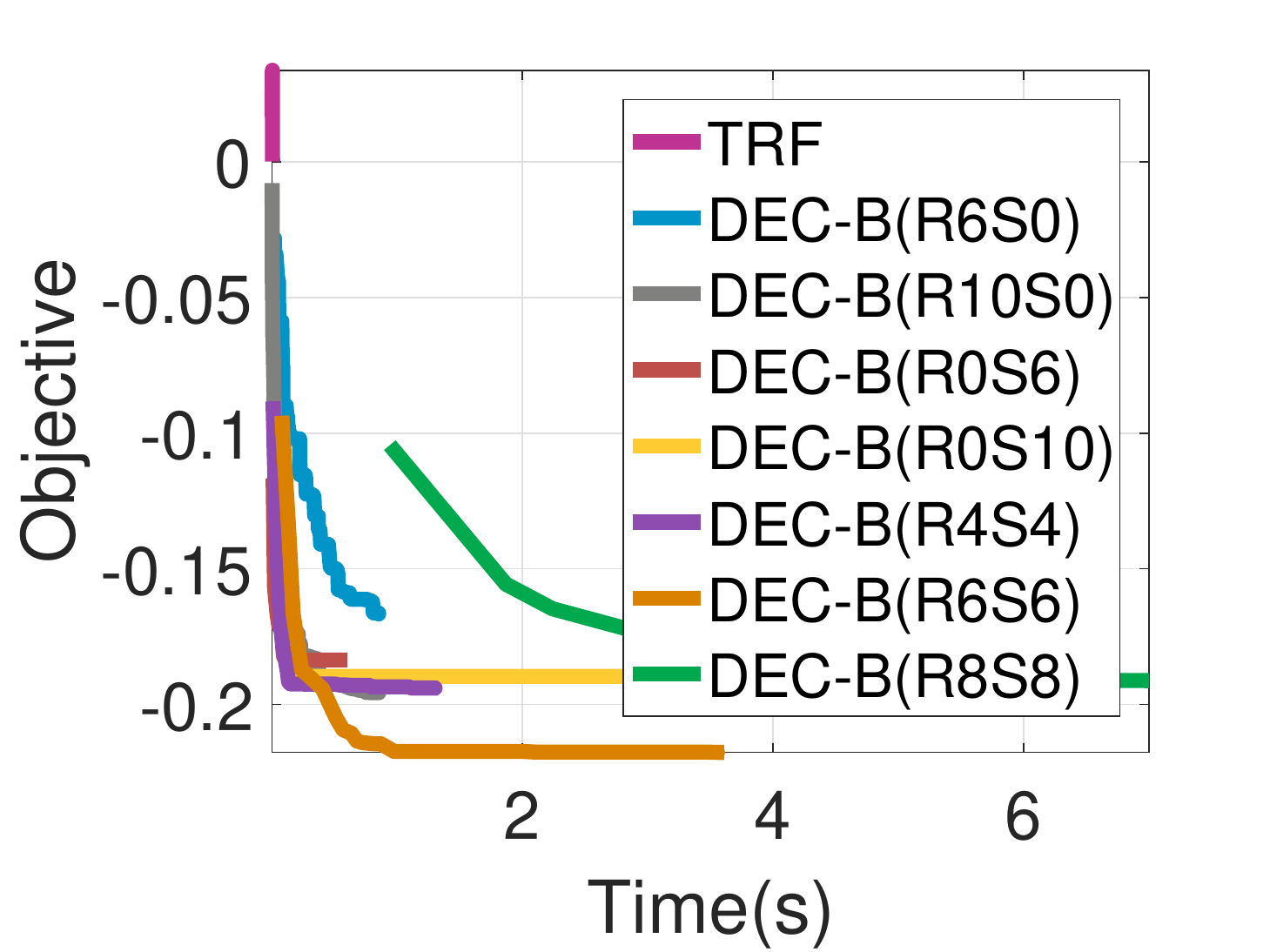}\vspace{-6pt}\caption{\footnotesize sparse CCA, `randn-2000', s=15}\end{subfigure}
\caption{Convergence behavior of different methods for sparse PCA (left column), sparse FDA (middle column), and sparse CCA (right column).}\label{fig:accuracy:cpu}
\end{figure*}

Now we present our main convergence result.

\begin{theorem} \label{theorem:convergence}

\textbf{Convergence Properties of Algorithm \ref{algo:main}.} Assume that the subproblem in (\ref{eq:subprob}) is solved globally, and there exists a constant $\sigma$ such that $\bbb{x}^{t}\bbb{C}\bbb{x}^{t}\geq \sigma >0$ for all $t$. We have the following results.

\bbb{(i)} When the random strategy is used to find the working set, we have $\lim_{t\rightarrow \infty} \E[\|\bbb{x}^{t+1} - \bbb{x}^t\|] = 0$ and Algorithm \ref{algo:main} converges to the block-$k$ stationary point in expectation.

\bbb{(ii)} When the swapping strategy is used to find the working set with $k\geq 2$, we have $\lim_{t\rightarrow \infty} \|\bbb{x}^{t+1} - \bbb{x}^t\| = 0$ and Algorithm \ref{algo:main} converges to the block-2 stationary point deterministically.

\end{theorem}
\noi \bbb{Remarks.} \bbb{(i)} Thanks to the fact that the denominator is positive and the objective function is quadratic fractional, our algorithm is still guaranteed to convergence even in the presence of non-convexity. \bbb{(ii)} We propose using a swapping strategy to find the working set which enumerates all possible swaps for all pairs of coordinates to find the greatest descent. One good feature of this strategy is that it achieves optimality guarantee which is no worse than Beck and Vaisbourd's coordinate-wise optimality condition \cite{beck2016sparse}.


%
%

\section{Experiments}\label{sect:exp}
This section demonstrates the efficacy of the proposed decomposition algorithm by considering three important applications (i.e., sparse PCA, sparse FDA, and sparse CCA) on synthetic and real-world data sets.

\bbb{$\bullet$ Data Sets}. \bbb{(i)} We consider four real-world data sets: `a1a', `w1a', `w2a', and `madelon'. We randomly select a subset of examples from the original data sets \footnote{\url{https://www.csie.ntu.edu.tw/~cjlin/libsvm/}}. The size of the data sets used in our experiments are $2000 \times 119$, $2000 \times 300$, $2000 \times 300$, $2000 \times 112$, respectively. \bbb{(ii)} We also use a similar method as in \cite{beck2016sparse} to generate synthetic Gaussian random data sets. Specifically, we produce the feature matrix $\bbb{X}\in\mathbb{R}^{m\times d}$ and the label vector $\bbb{y}\in\mathbb{R}^m$ as follows: $\bbb{X}=\text{randn}(m,d),~\bbb{y}=\text{sign}(\text{randn}(m,1))$, where $\text{randn}(m,d)$ is a function that returns a standard Gaussian random matrix of size $m \times d$ and $\text{sign}$ is a signum function. We fix $m=300$ and consider different values for $d=\{100,500,1500,2000\}$. We denote the data sets as `randn-d' and place the results in the \bbb{Appendix}.

Based on $\bbb{X}$ and $\bbb{y}$, we generate the matrices $\bbb{A}$ and $\bbb{C}$ in Problem (\ref{eq:main}) for different applications (see Section \ref{sect:app}). Note that the resulting size of the sparse generalized eigenvalue problem for sparse PCA, sparse FDA, and sparse CCA are $d$, $d$, and $m$, respectively. We vary the sparsity parameter $s\in\{4,8,12,...,40\}$ and report the objective values for Problem (\ref{eq:main}).

\begin{figure*} [!th]

\setcounter{subfigure}{0}
\centering
      \begin{subfigure}{\fourfigwid}\includegraphics[width=\textwidth,height=\objimghei]{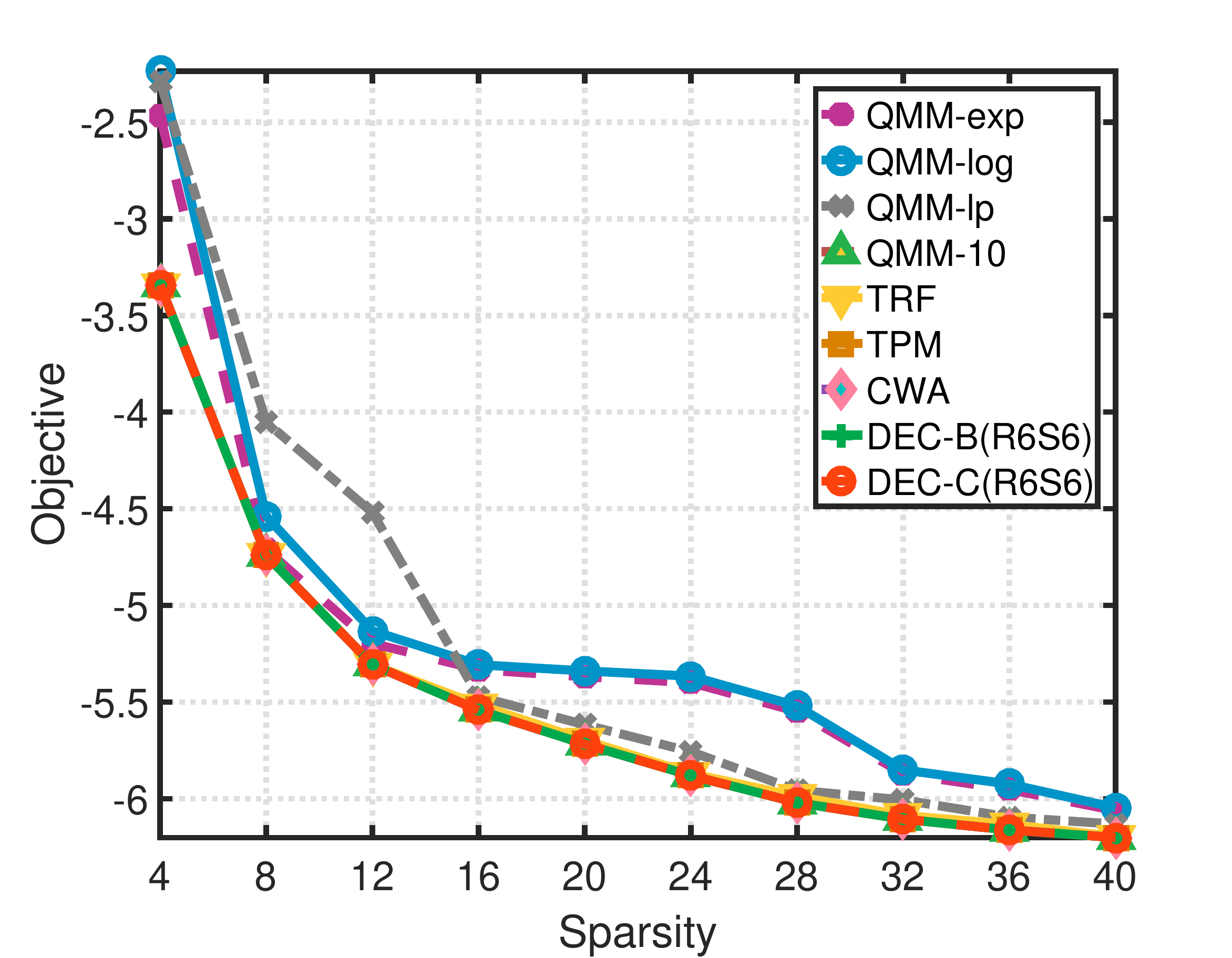}\vspace{-6pt}\caption{\footnotesize a1a}\end{subfigure}\ghs
      \begin{subfigure}{\fourfigwid}\includegraphics[width=\textwidth,height=\objimghei]{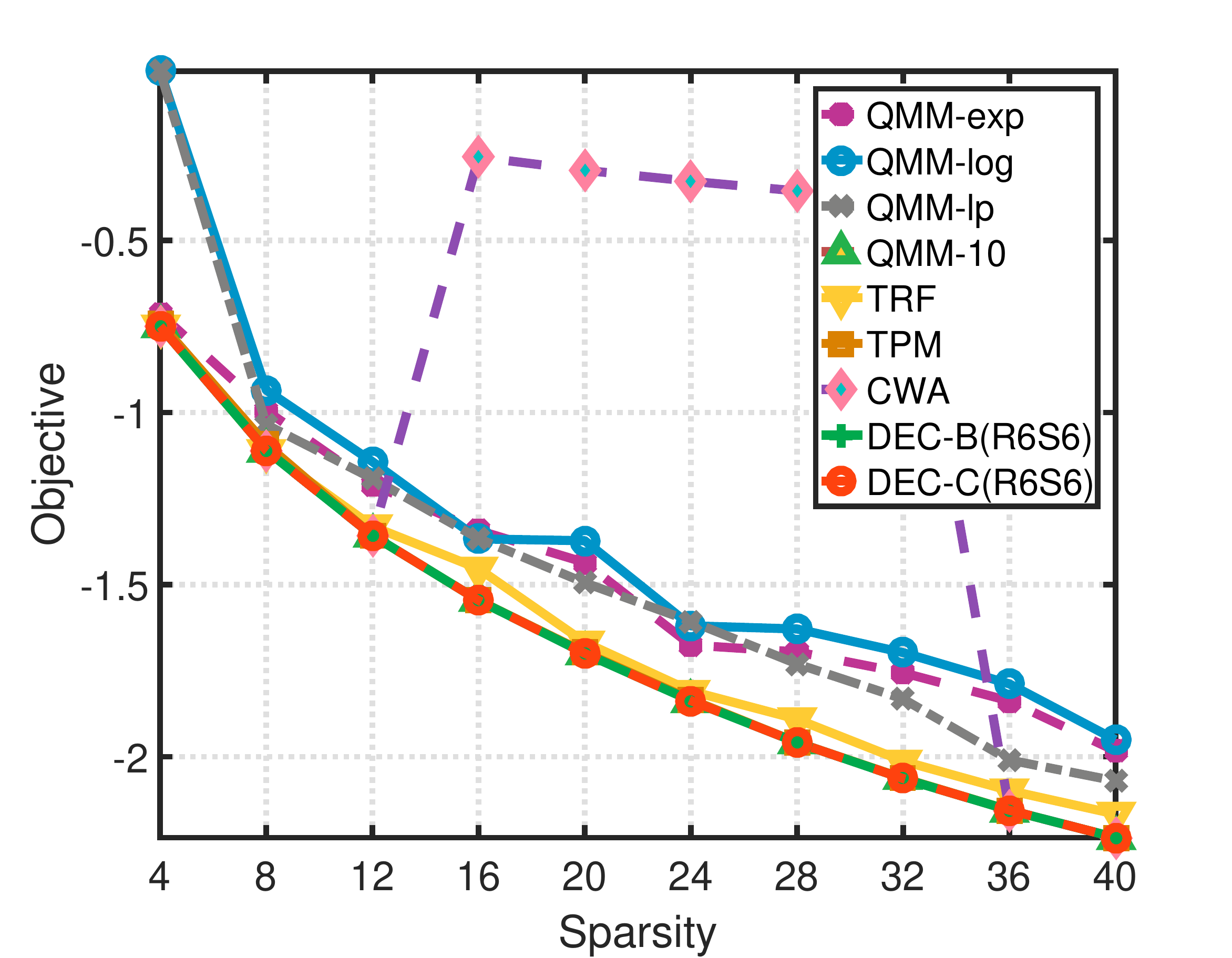}\vspace{-6pt}\caption{\footnotesize w1a}\end{subfigure}\ghs
      \begin{subfigure}{\fourfigwid}\includegraphics[width=\textwidth,height=\objimghei]{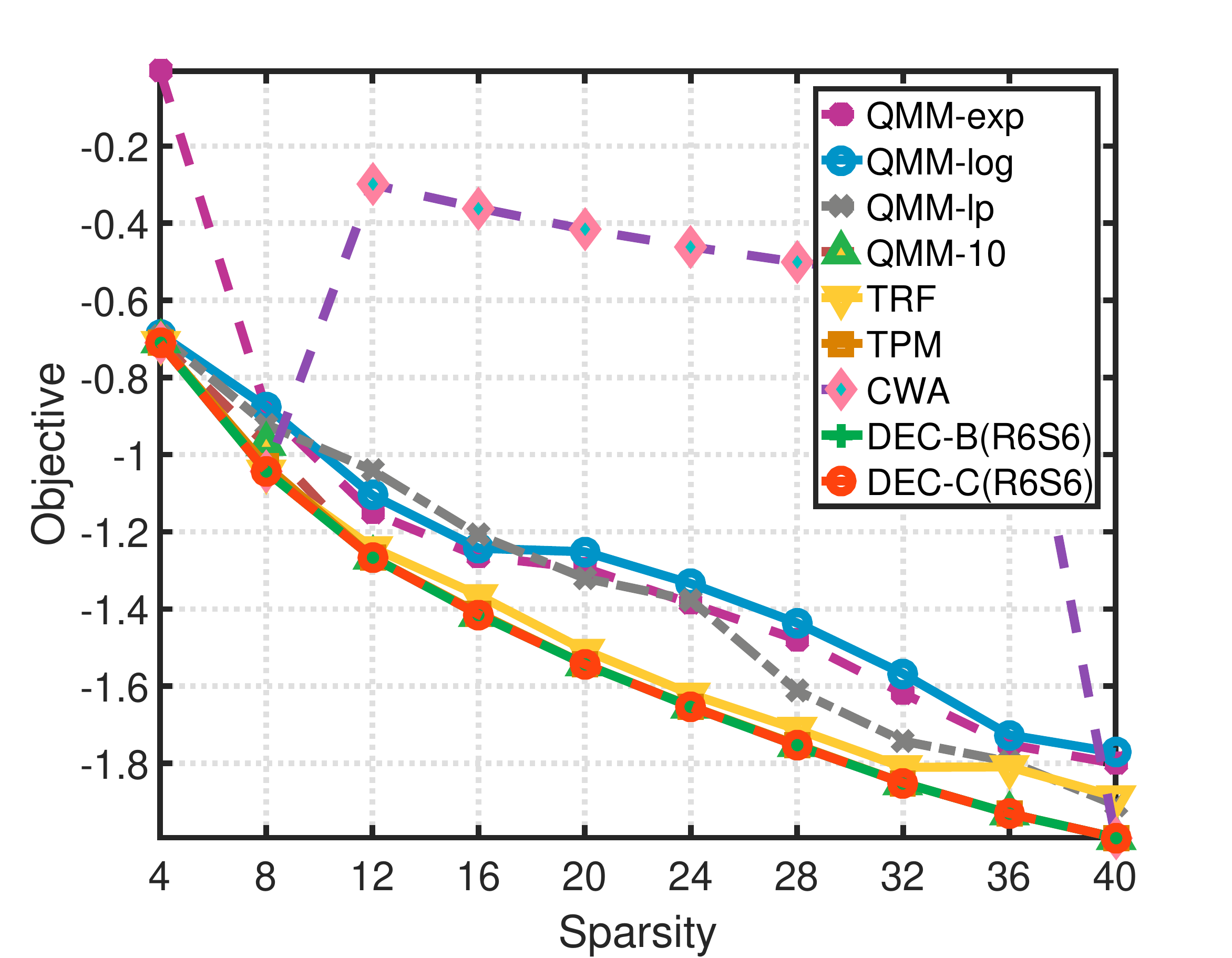}\vspace{-6pt}\caption{\footnotesize w2a}\end{subfigure}\ghs
      \begin{subfigure}{\fourfigwid}\includegraphics[width=\textwidth,height=\objimghei]{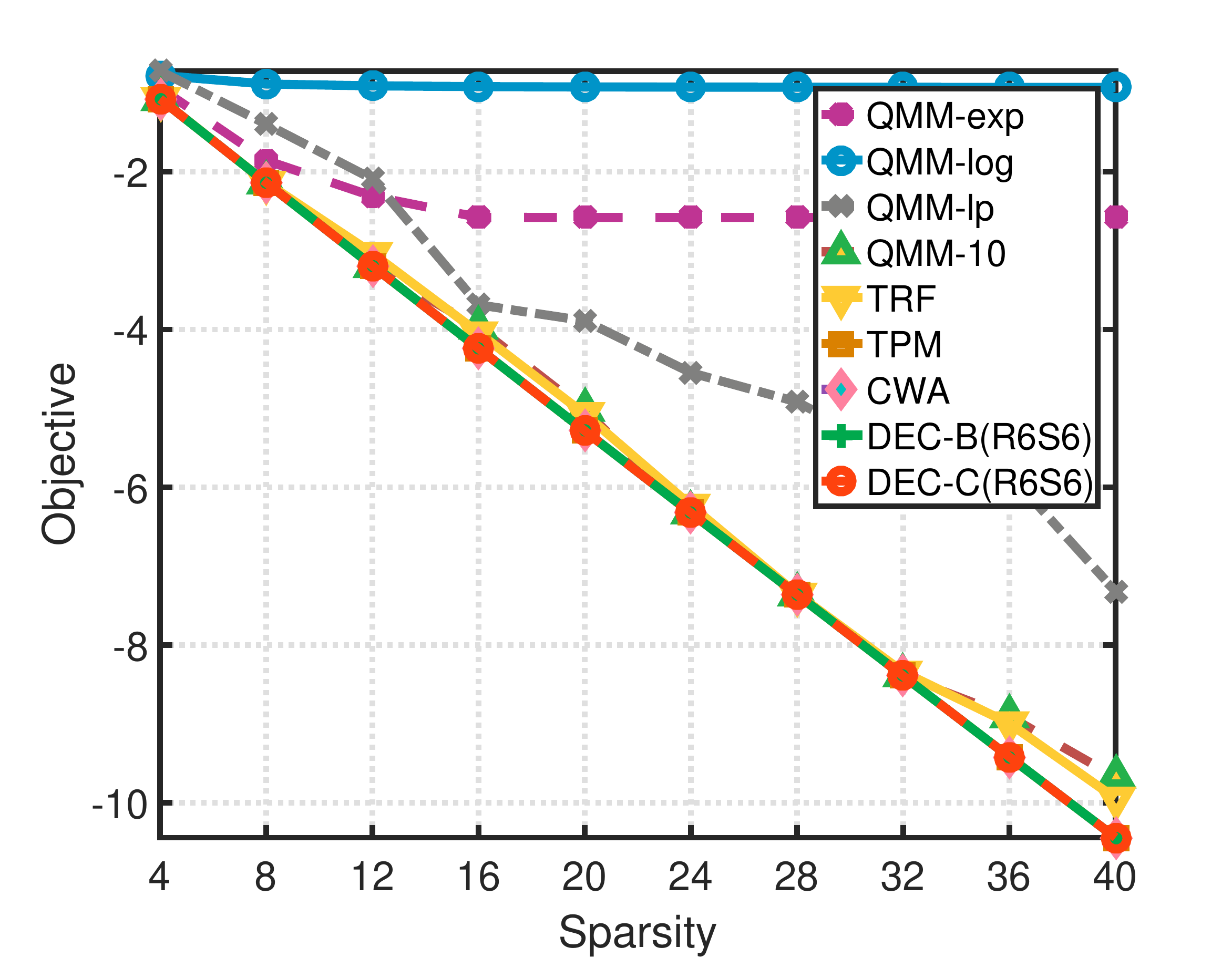}\vspace{-6pt} \caption{\footnotesize madelon}\end{subfigure}

\caption{Accuracy of different methods on different data sets for sparse PCA problem with varying the cardinalities.}\label{fig:accuracy:pca}

\setcounter{subfigure}{0}
      \centering

      \begin{subfigure}{\fourfigwid}\includegraphics[width=\textwidth,height=\objimghei]{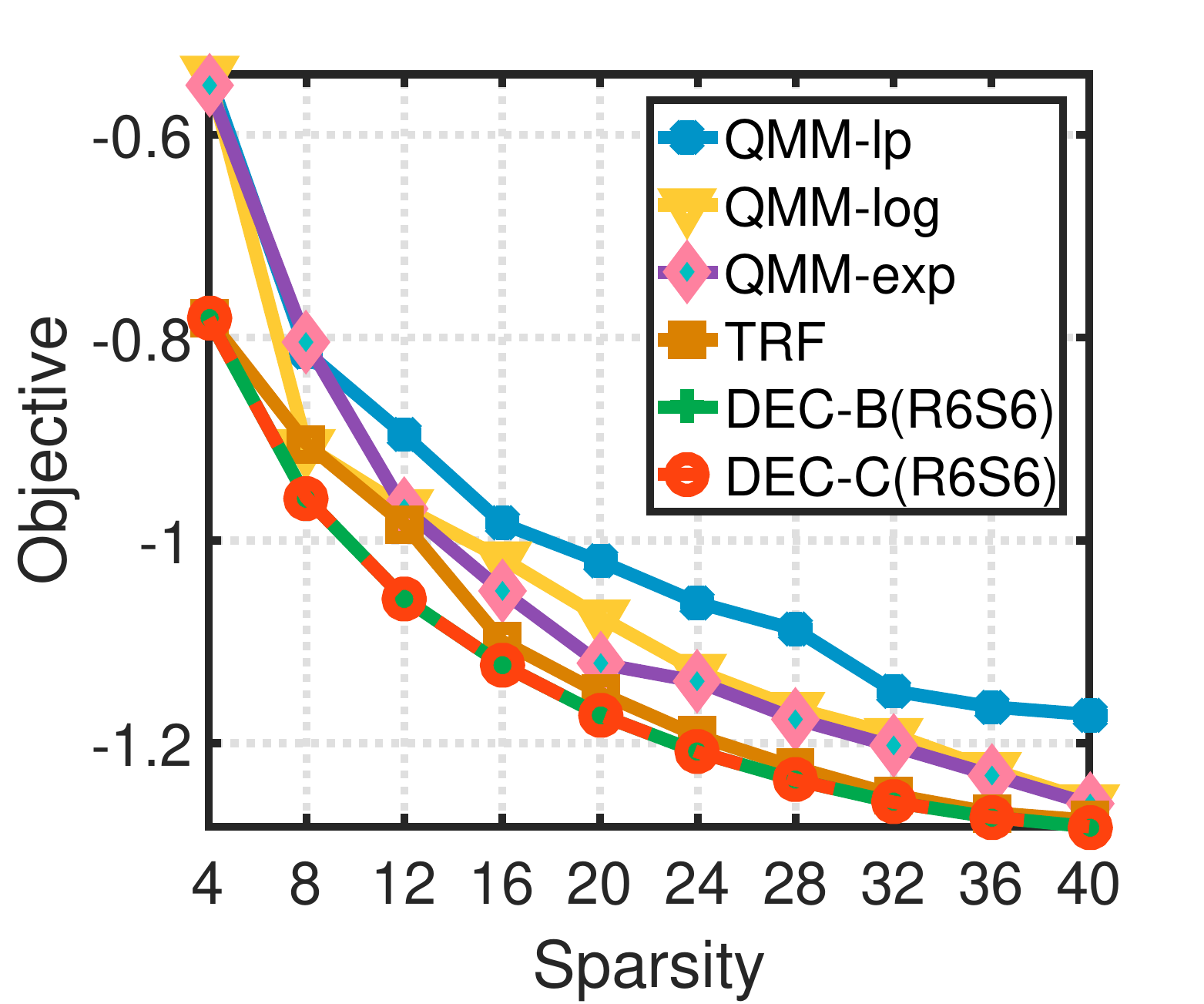}\vspace{-6pt}\caption{\footnotesize a1a}\end{subfigure}\ghs
      \begin{subfigure}{\fourfigwid}\includegraphics[width=\textwidth,height=\objimghei]{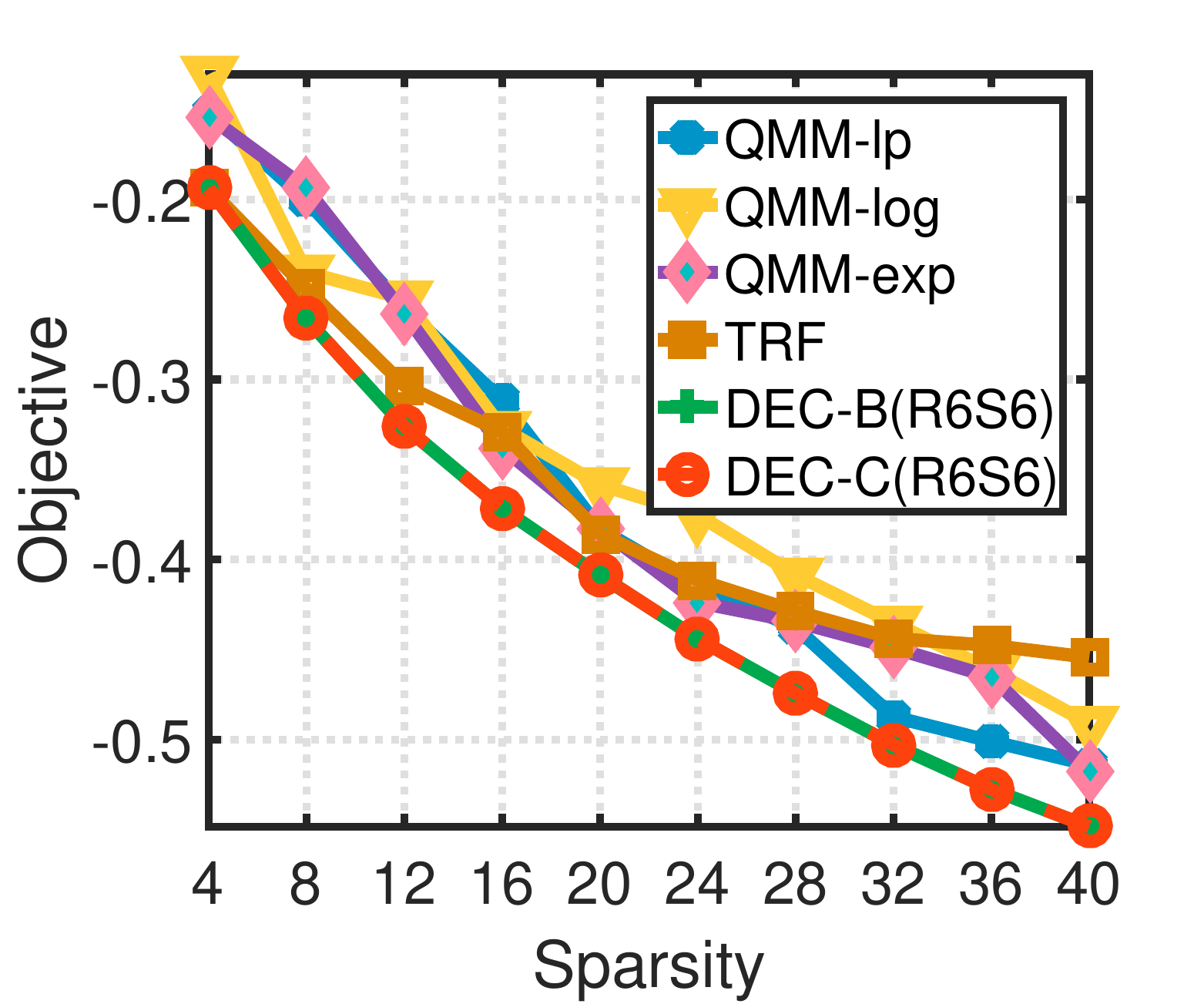}\vspace{-6pt}\caption{\footnotesize w1a}\end{subfigure}\ghs
      \begin{subfigure}{\fourfigwid}\includegraphics[width=\textwidth,height=\objimghei]{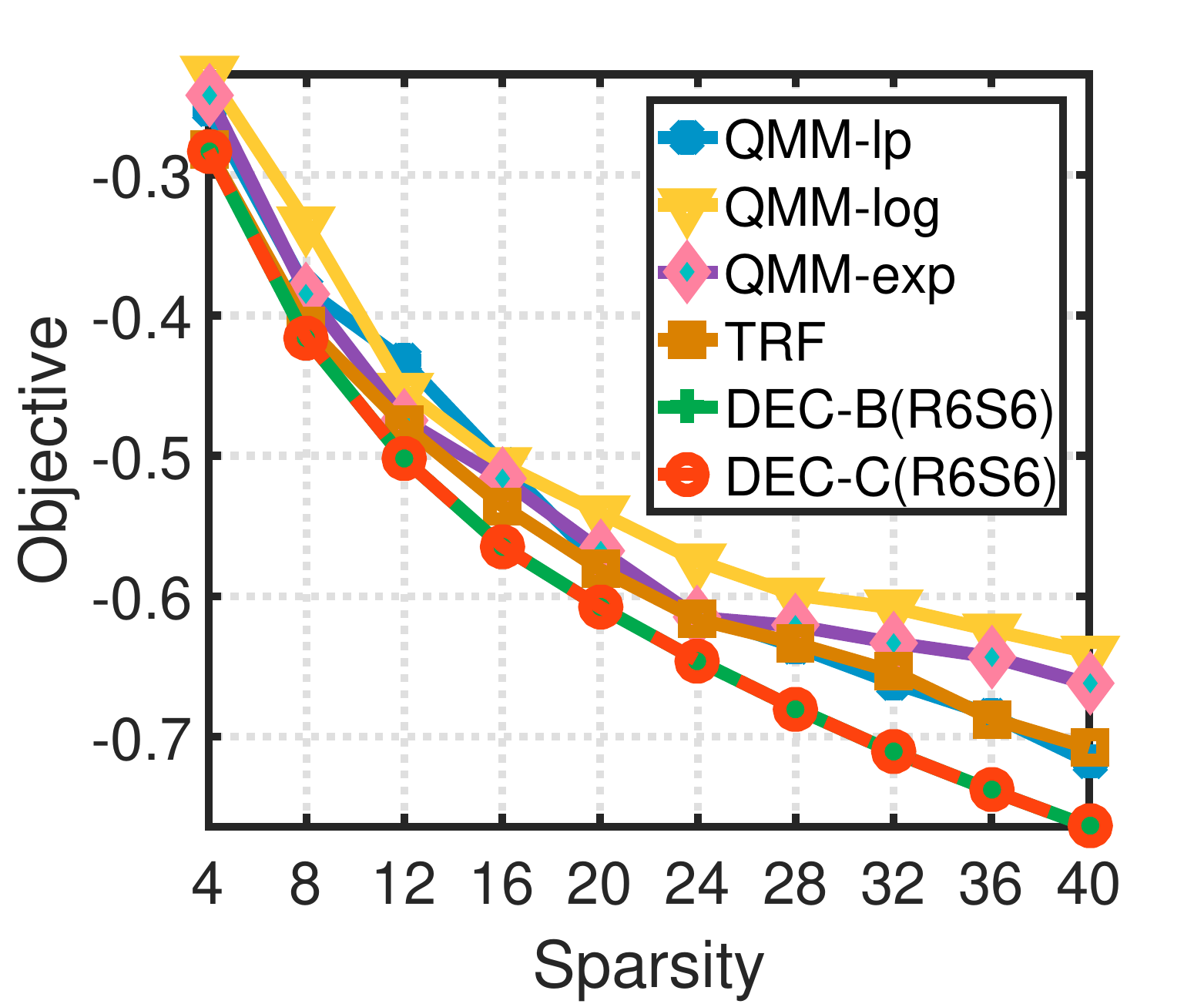}\vspace{-6pt}\caption{\footnotesize w2a}\end{subfigure}\ghs
      \begin{subfigure}{\fourfigwid}\includegraphics[width=\textwidth,height=\objimghei]{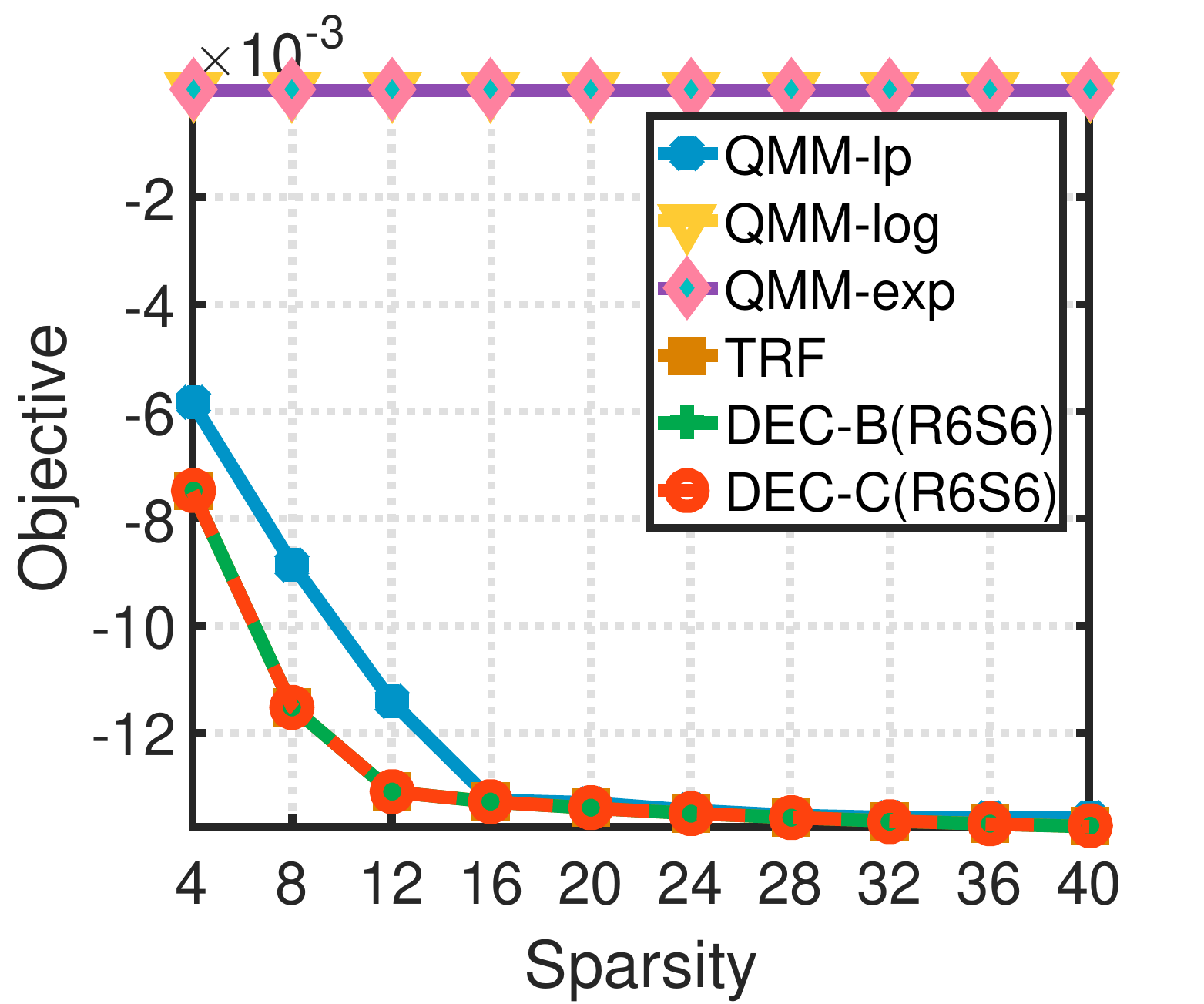}\vspace{-6pt} \caption{\footnotesize madelon}\end{subfigure}



\caption{Accuracy of different methods on different data sets for sparse FDA problem with varying cardinalities.}\label{fig:accuracy:fda}

\setcounter{subfigure}{0}
\centering
      \begin{subfigure}{\fourfigwid}\includegraphics[width=\textwidth,height=\objimghei]{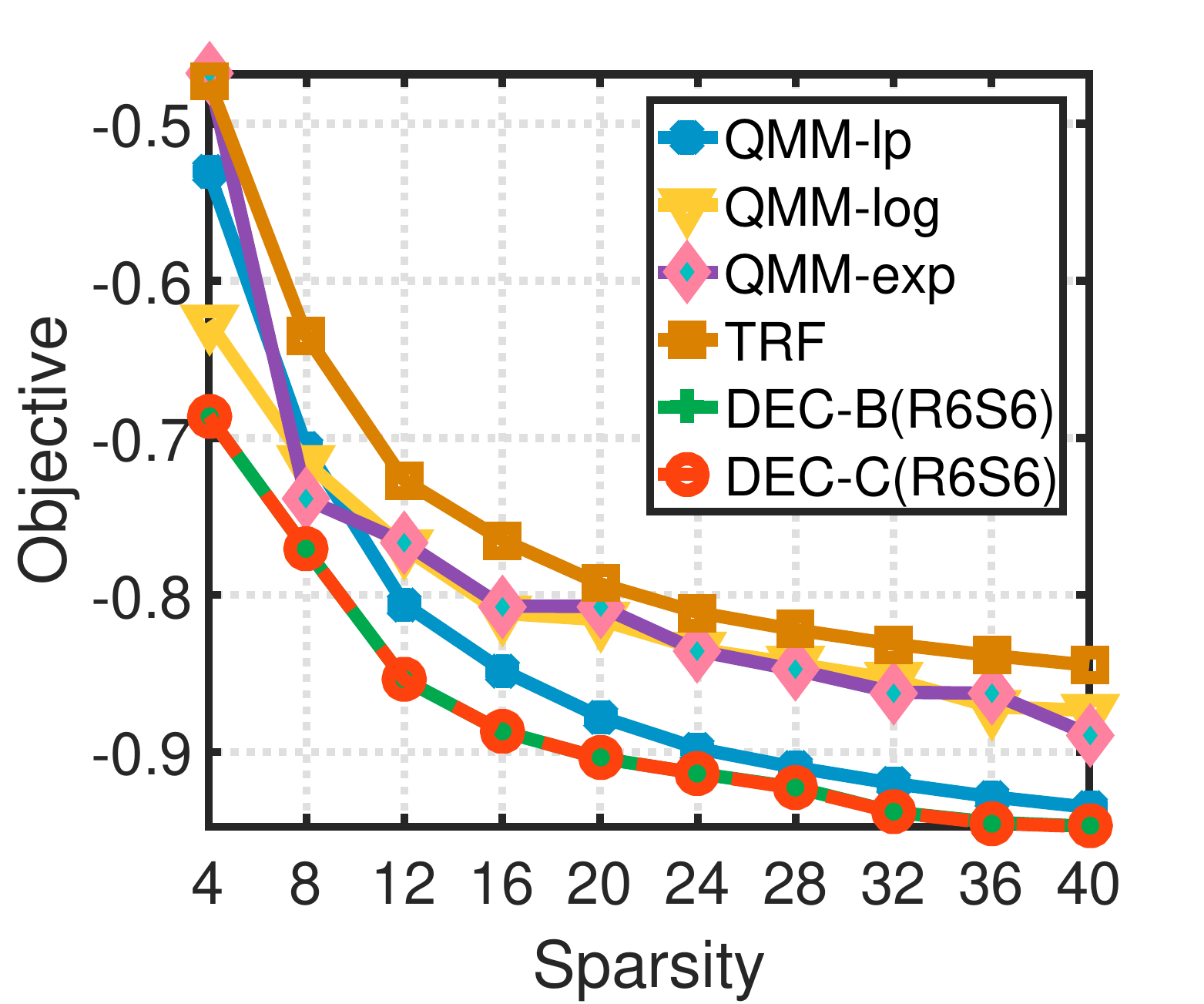}\vspace{-6pt}\caption{\footnotesize a1a}\end{subfigure}\ghs
      \begin{subfigure}{\fourfigwid}\includegraphics[width=\textwidth,height=\objimghei]{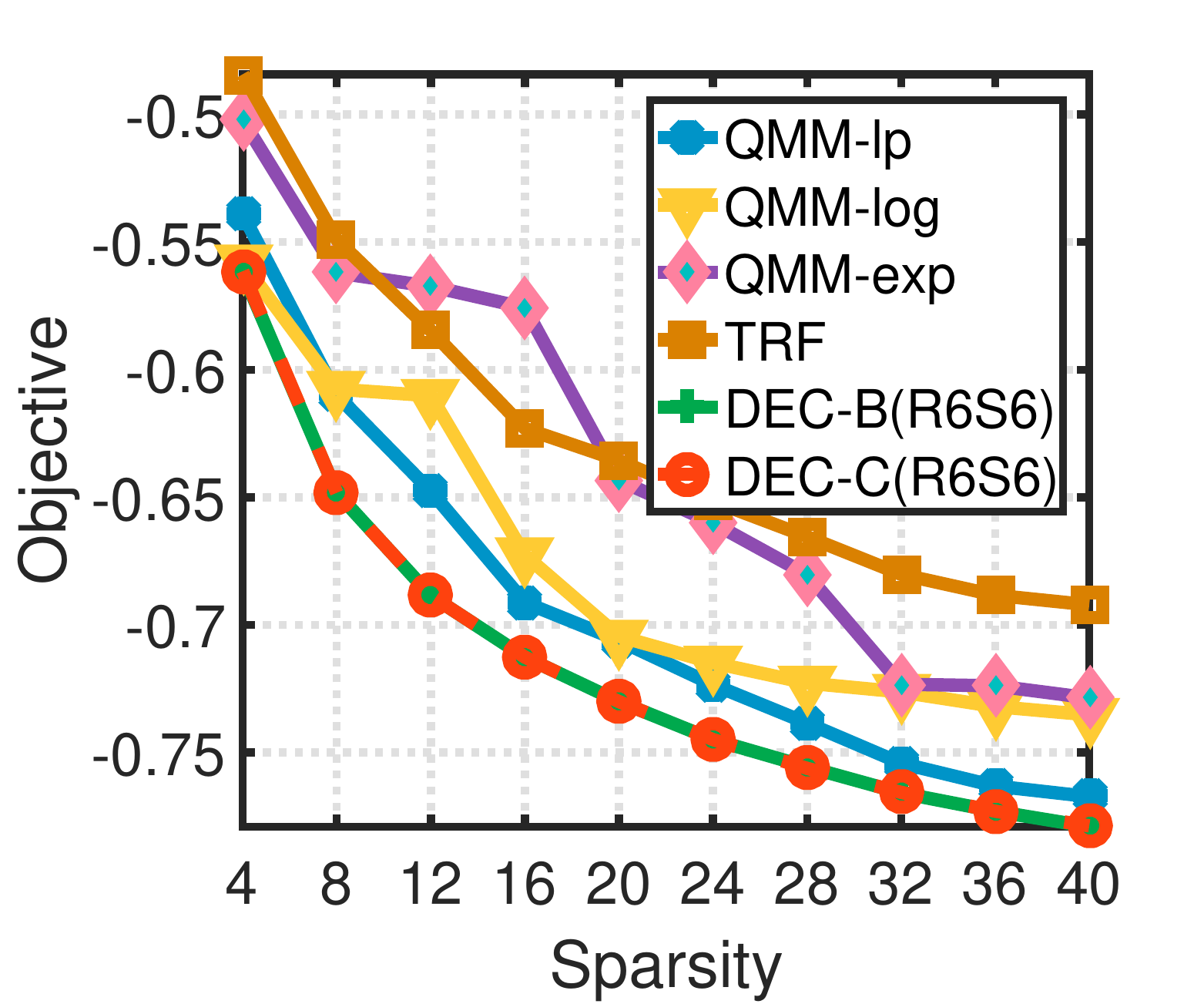}\vspace{-6pt}\caption{\footnotesize w1a}\end{subfigure}\ghs
      \begin{subfigure}{\fourfigwid}\includegraphics[width=\textwidth,height=\objimghei]{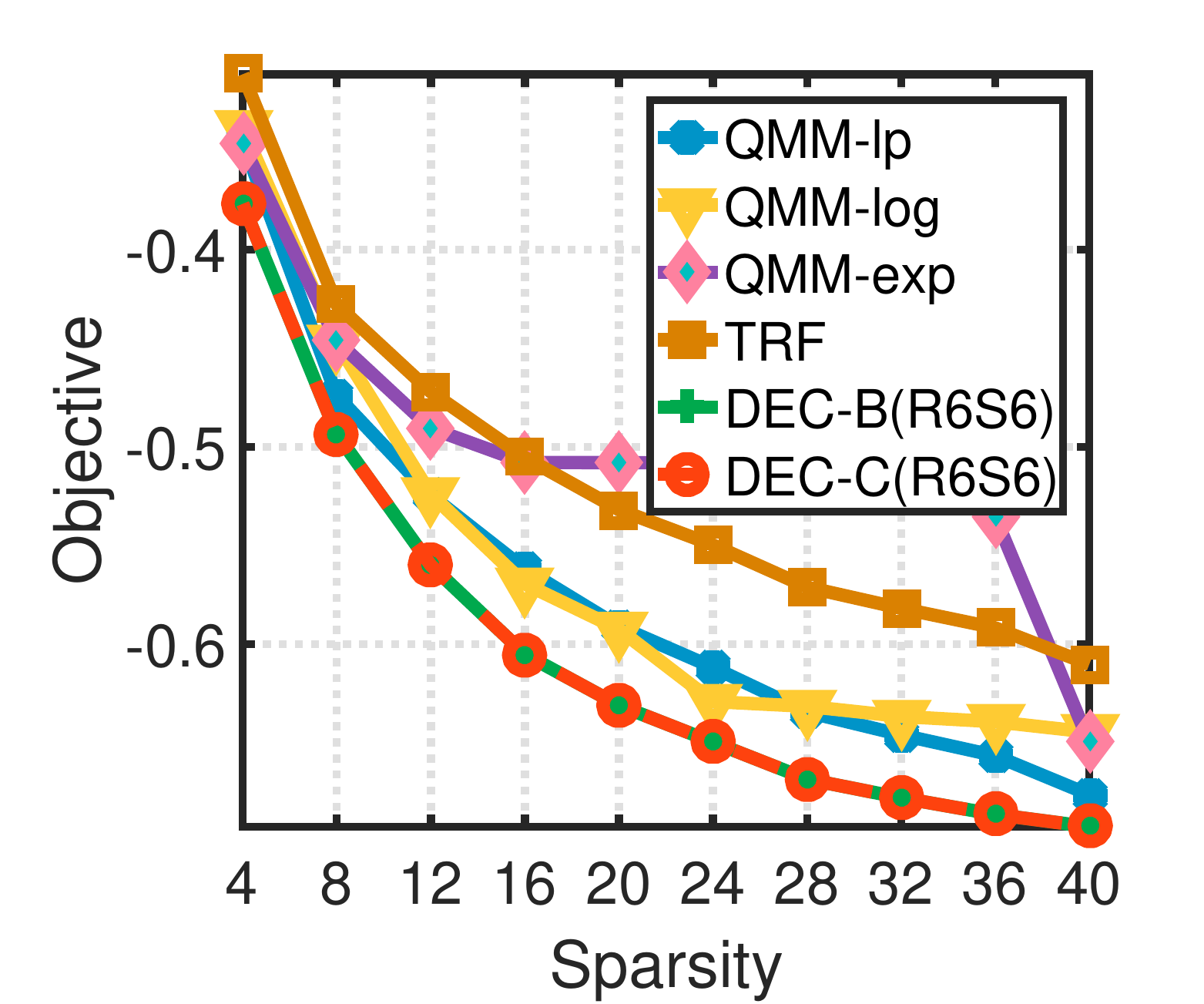}\vspace{-6pt}\caption{\footnotesize w2a}\end{subfigure}\ghs
      \begin{subfigure}{\fourfigwid}\includegraphics[width=\textwidth,height=\objimghei]{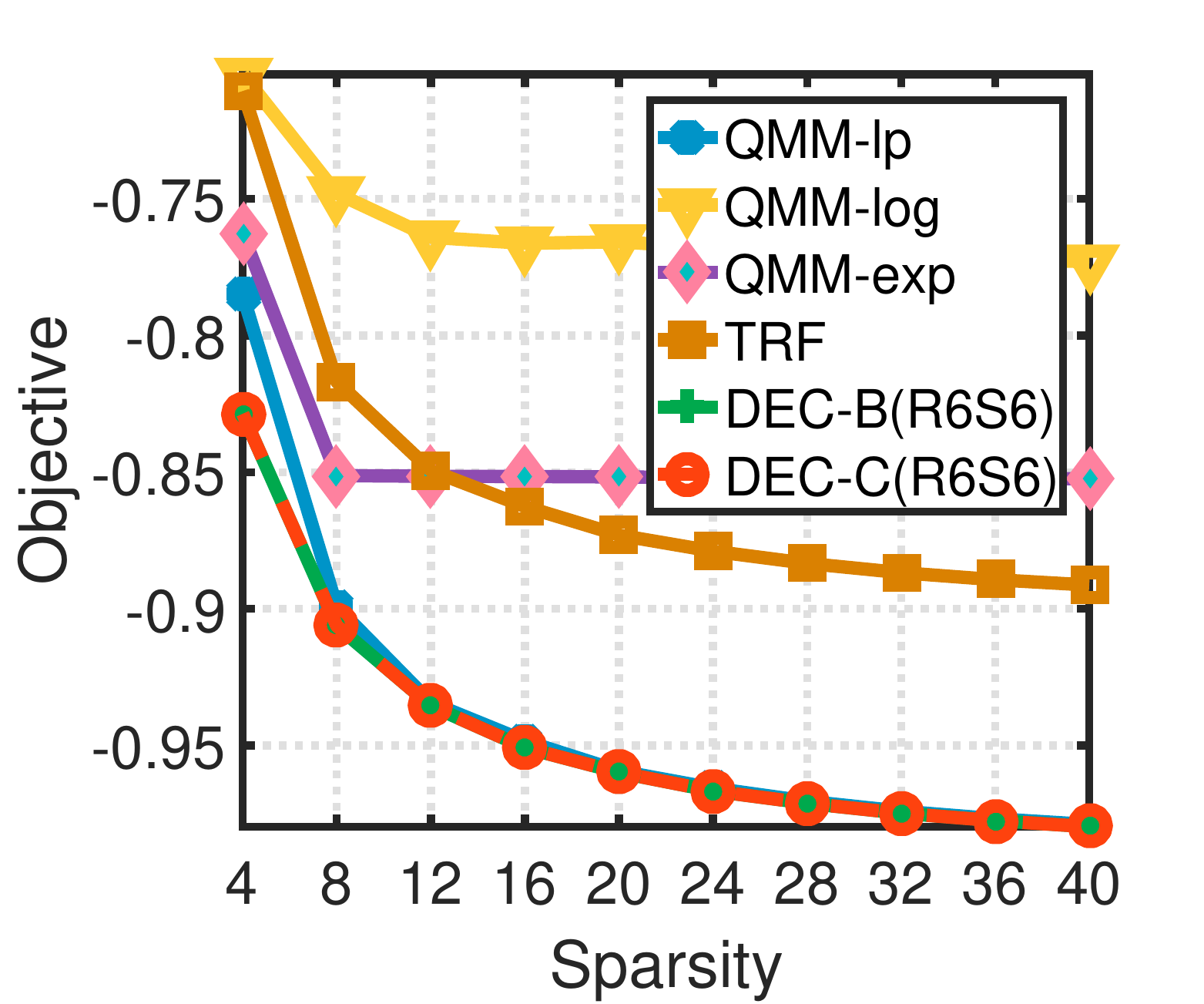}\vspace{-6pt} \caption{\footnotesize madelon}\end{subfigure}


\caption{Accuracy of different methods on different data sets for sparse CCA problem with varying cardinalities.}\label{fig:accuracy:cca}
\end{figure*}

%
%

\begin{figure*} [!t]
\captionsetup{singlelinecheck = on, format= hang, justification=justified, font=footnotesize, labelsep=space}
      \centering
      \begin{subfigure}{\threefigwid}\includegraphics[width=\textwidth,height=\objimghei]{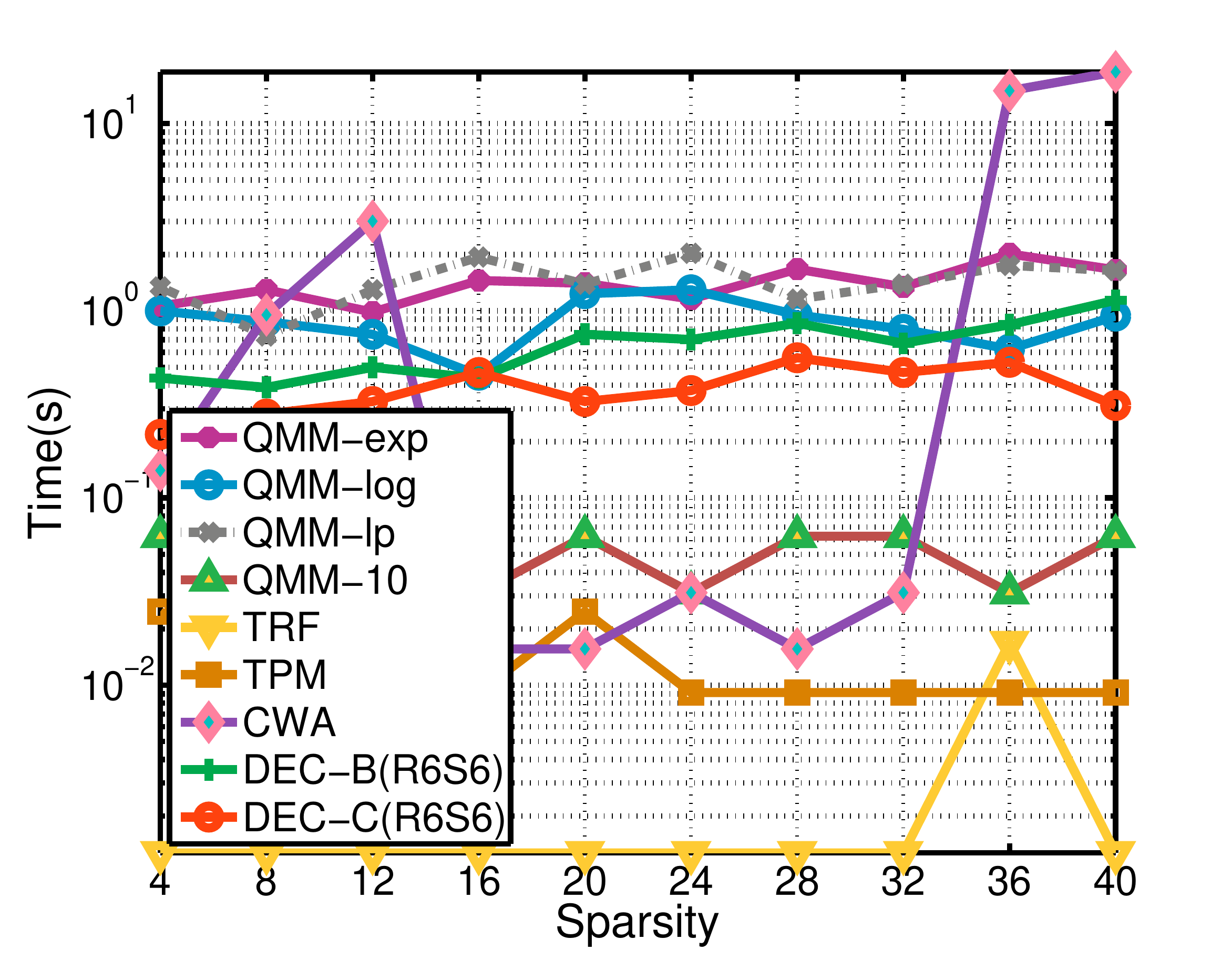}\vspace{-6pt}\caption{PCA}\end{subfigure}\ghs
      \begin{subfigure}{\threefigwid}\includegraphics[width=\textwidth,height=\objimghei]{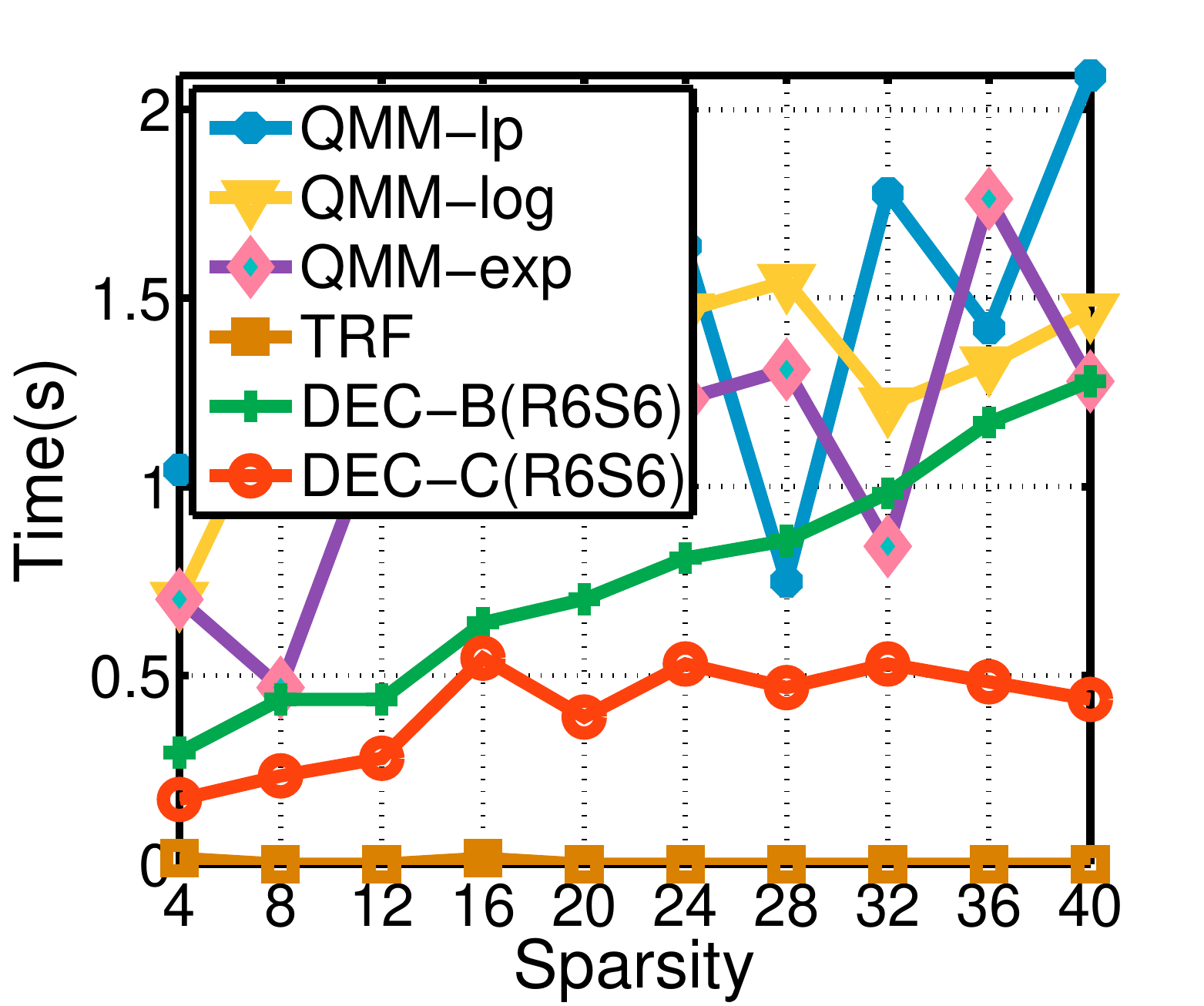}\vspace{-6pt}\caption{\footnotesize FDA}\end{subfigure}\ghs
      \begin{subfigure}{\threefigwid}\includegraphics[width=\textwidth,height=\objimghei]{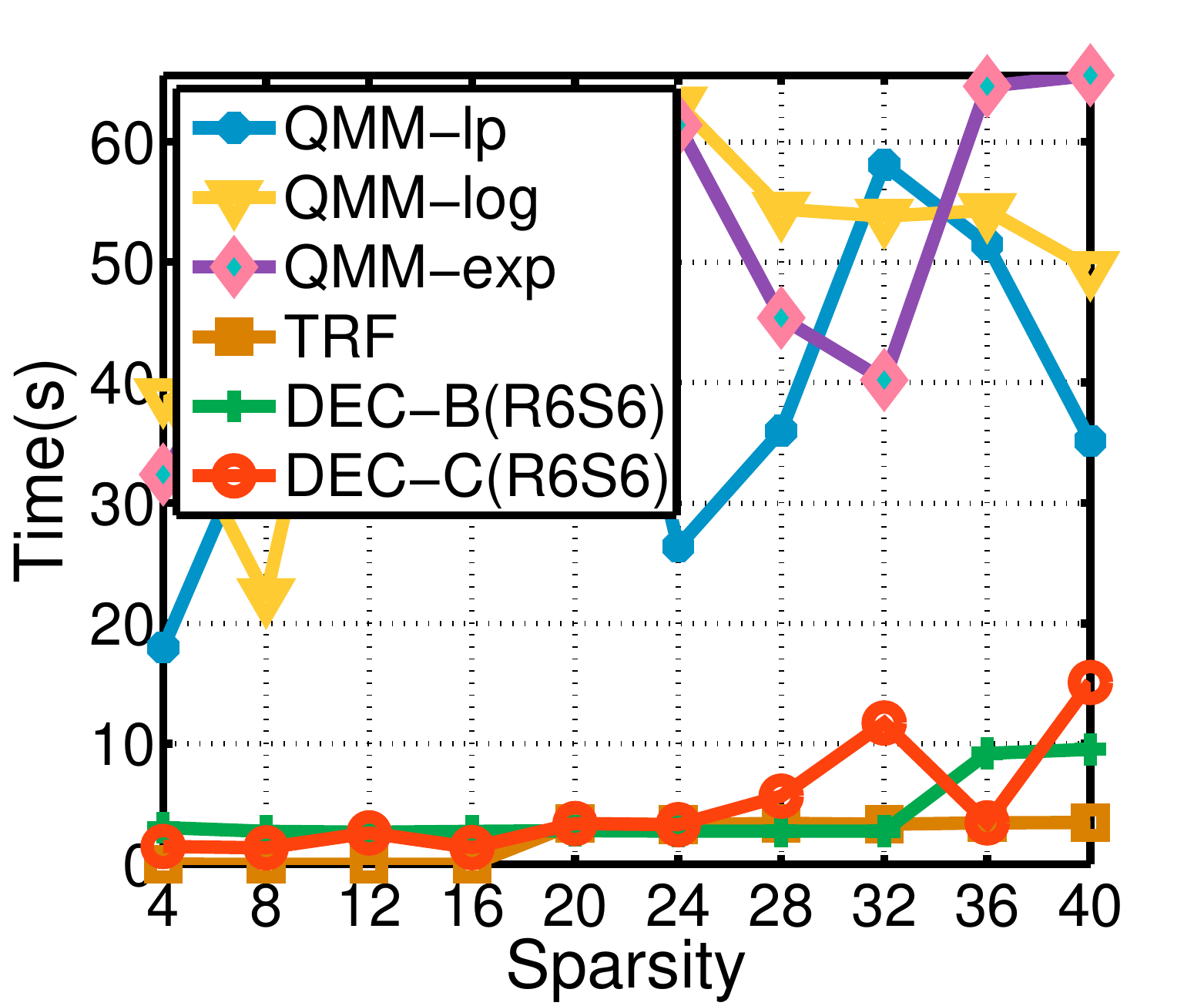}\vspace{-6pt}\caption{\footnotesize CCA}\end{subfigure}
\caption{Comparison of the computing time for different methods on `w1a' data set with varying cardinalities.}\label{fig:accuracy:cpu2}
\end{figure*}

\bbb{$\bullet$ Compared Methods}. We compare the following methods. \bbb{(i)} Truncated Power Method (TPM) \cite{yuan2013truncated} \footnote{code: \url{sites.google.com/site/xtyuan1980/}} iteratively and greedily decreases the objective while maintaining the desired sparsity for the solutions by hard thresholding truncation. \bbb{(ii)} Coordinate-Wise Algorithm (CWA) \cite{beck2013sparsity,beck2016sparse} \footnote{code: \url{sites.google.com/site/amirbeck314/}} iteratively performs an optimization step with respect to two coordinates, where the coordinates that need to be altered are chosen to be the ones that produce the maximal decrease among all possible alternatives. \bbb{(iii)} Truncated Rayleigh Flow (TRF) \cite{tan2016sparse} iteratively updates the solution using the gradient of the generalized Rayleigh quotient and performs a truncation operation to achieve sparsity. \bbb{(iv)} Quadratic Majorization Method (QMM) \cite{SongBP15} \footnote{code: \url{https://junxiaosong.github.io/}} approximates the $\ell_0$-norm by a continuous surrogate function and iteratively majorizes the surrogate function by a quadratic separable function, which at each iteration reduces to a regular generalized eigenvalue problem. Using different smooth non-convex approximation functions, they develop different versions of QMM (QMM-exp, QMM-log, QMM-$\ell_p$, QMM-$\ell_0$). Since their methods only solve a regularized problem and fail to control the sparsity of the solution, we use a simple bisection search to find the best regulation parameter and report the lowest objective value after hard thresholding. \bbb{(v)} The proposed decomposition method (denoted as DEC) is included for comparisons. We use DEC-B(Ri-Sj) and DEC-C(Ri-Sj) to denote our method based on a \bbb{B}isection search method and a \bbb{C}oordinate descent method, respectively, along with selecting $i$ coordinate using the \bbb{R}andom strategy and $j$ coordinates using the \bbb{S}wapping strategy. In each iteration, we compute $r_t = (f(\bbb{x}^t)-f(\bbb{x}^{t+1}))/f(\bbb{x}^t)$. We let Algorithm \ref{algo:main} run up to $T$ iterations and stop it at iteration $t<T$ if $\text{mean}([{r}_{t-\text{min}(t,M)+1},~{r}_{t-min(t,M)+2},...,{r}_t]) \leq \epsilon$. The default parameter $(\theta,~\epsilon,~M,~T)=(10^{-5},~10^{-5},~50,~1000)$ is used. All codes are implemented in MATLAB on an Intel 3.20GHz CPU with 8 GB RAM \footnote{For the purpose of reproducibility, we provide our code in the authors' research webpage.}. Only DEC-C is developed in C and wrapped into our MATLAB code, since it uses an elementwise loop which is inefficient in native MATLAB.

We remark that both (a) and (b) are only designed for sparse PCA with $\bbb{C}=\bbb{I}$. We do not compare against the DC programming algorithms \cite{Sriperumbudur2007,thiao2010dc} since they fail to control the sparsity of the solution and result in worse accuracy than QMM (see \cite{SongBP15}).




\bbb{$\bullet$ Convergence Behavior}. We show the convergence behavior for different methods in Figure \ref{fig:accuracy:cpu}. We do not include the results of QMM since it fails to control the sparsity of the solution. Due to space limitation, we only report the results of DEC-B in this set of experiments. We have the following observations. \bbb{(i)} The methods $\{$TPM, CWA, TRF$\}$ converge within one second and they are faster than DEC. However, they get stuck into poor local minima and result in much worse accuracy than DEC. \bbb{(ii)} The objective values of DEC stabilize after less than 5 seconds, which means it has converged, and the decrease of the objective value is negligible afterwards. This implies that one may use a looser stopping criterion without sacrificing accuracy. \bbb{(iii)} DEC-B(R6S0) and DEC-B(R8S8) converge slowly, and it seems that DEC-B(R6S6) finds a good trade-off between efficiency and effectiveness. \bbb{(iv)} DEC-B(R10S0) achieves a lower objective value than DEC-B(R6S0). This is reasonable since a larger $k$ in the block-$k$ optimality condition implies a stronger stationary point. \bbb{(v)} $\{$DEC-B(R6S0), DEC-B(R10S0)$\}$ achieve larger objective values than $\{$DEC-B(R0S6), DEC-B(R0S10)$\}$, which implies that the swapping strategy plays an indispensable role in DEC.


\bbb{$\bullet$ Experimental Results}. We show the experimental results for sparse PCA, sparse FDA, and sparse CCA in Figure \ref{fig:accuracy:pca}, \ref{fig:accuracy:fda} and \ref{fig:accuracy:cca}, respectively. Several conclusions can be drawn. \bbb{(i)} CWA generally outperforms $\{$TPM, TRF, QMM$\}$. \bbb{(ii)} CWA is not stable and generates much worse results on `w1a' and `w2a'. \bbb{(iii)} The proposed method DEC still outperforms CWA and achieves the lowest objective values. \bbb{(iv)} Both DEC-B and DEC-C perform similarly. This is because coordinate descent methods find a desirable solution for the quadratic fractional programming problem.

\bbb{$\bullet$ Computational Efficiency}. We demonstrate a comparison of the actual computing time for different methods on `w1a' data set in Figure \ref{fig:accuracy:cpu2}. Two conclusions can be drawn. \bbb{(i)} DEC takes less than 15 seconds to converge in all our instances. \bbb{(ii)} DEC is practical and it is much more efficient than QMM.

\section{Conclusions}


This paper presents an effective algorithm for solving the sparse generalized eigenvalue problem. Although the problem is NP-hard, we consider a decomposition algorithm to solve it. Our experiments on synthetic data and real-world data have shown that our method significantly and consistently outperforms existing solutions in term of accuracy.



{
\normalem
\bibliographystyle{ieee}
\bibliography{my}}

\begin{thebibliography}{10}\itemsep=-1pt

\bibitem{asteris2016simple}
M.~Asteris, A.~Kyrillidis, O.~Koyejo, and R.~Poldrack.
\newblock A simple and provable algorithm for sparse diagonal cca.
\newblock In {\em International Conference on Machine Learning (ICML)}, pages
  1148--1157, 2016.

\bibitem{AsterisPD14}
M.~Asteris, D.~S. Papailiopoulos, and A.~G. Dimakis.
\newblock Nonnegative sparse {PCA} with provable guarantees.
\newblock In {\em International Conference on Machine Learning (ICML)}, pages
  1728--1736, 2014.

\bibitem{beck2013sparsity}
A.~Beck and Y.~C. Eldar.
\newblock Sparsity constrained nonlinear optimization: Optimality conditions
  and algorithms.
\newblock {\em SIAM Journal on Optimization}, 23(3):1480--1509, 2013.

\bibitem{beck2010minimizing}
A.~Beck and M.~Teboulle.
\newblock On minimizing quadratically constrained ratio of two quadratic
  functions.
\newblock {\em Journal of Convex Analysis}, 17(3):789--804, 2010.

\bibitem{beck2016sparse}
A.~Beck and Y.~Vaisbourd.
\newblock The sparse principal component analysis problem: Optimality
  conditions and algorithms.
\newblock {\em Journal of Optimization Theory and Applications},
  170(1):119--143, 2016.

\bibitem{boyd2004convex}
S.~Boyd and L.~Vandenberghe.
\newblock {\em Convex optimization}.
\newblock Cambridge university press, 2004.

\bibitem{chang2011libsvm}
C.-C. Chang and C.-J. Lin.
\newblock Libsvm: a library for support vector machines.
\newblock {\em ACM transactions on Intelligent Systems and Technology (TIST)},
  2(3):27, 2011.

\bibitem{d2008optimal}
A.~d’Aspremont, F.~Bach, and L.~E. Ghaoui.
\newblock Optimal solutions for sparse principal component analysis.
\newblock {\em Journal of Machine Learning Research}, 9(Jul):1269--1294, 2008.

\bibitem{d2005direct}
A.~d'Aspremont, L.~E. Ghaoui, M.~I. Jordan, and G.~R. Lanckriet.
\newblock A direct formulation for sparse pca using semidefinite programming.
\newblock In {\em Advances in Neural Information Processing Systems (NIPS)},
  pages 41--48, 2005.

\bibitem{dinkelbach1967nonlinear}
W.~Dinkelbach.
\newblock On nonlinear fractional programming.
\newblock {\em Management Science}, 13(7):492--498, 1967.

\bibitem{ge2016efficient}
R.~Ge, C.~Jin, P.~Netrapalli, A.~Sidford, et~al.
\newblock Efficient algorithms for large-scale generalized eigenvector
  computation and canonical correlation analysis.
\newblock In {\em International Conference on Machine Learning (ICML)}, pages
  2741--2750, 2016.

\bibitem{GuoLYSW03}
Y.~Guo, S.~Li, J.~Yang, T.~Shu, and L.~Wu.
\newblock A generalized foley-sammon transform based on generalized fisher
  discriminant criterion and its application to face recognition.
\newblock {\em Pattern Recognition Letters}, 24(1-3):147--158, 2003.

\bibitem{horn1990matrix}
R.~A. Horn and C.~R. Johnson.
\newblock {\em Matrix analysis}.
\newblock Cambridge university press, 1990.

\bibitem{HsiehCLKS08}
C.~Hsieh, K.~Chang, C.~Lin, S.~S. Keerthi, and S.~Sundararajan.
\newblock A dual coordinate descent method for large-scale linear {SVM}.
\newblock In {\em International Conference on Machine Learning (ICML)}, pages
  408--415, 2008.

\bibitem{HsiehD11}
C.~Hsieh and I.~S. Dhillon.
\newblock Fast coordinate descent methods with variable selection for
  non-negative matrix factorization.
\newblock In {\em {ACM} International Conference on Knowledge Discovery and
  Data Mining (SIGKDD)}, pages 1064--1072, 2011.

\bibitem{jeffers1967two}
J.~Jeffers.
\newblock Two case studies in the application of principal component analysis.
\newblock {\em Applied Statistics}, pages 225--236, 1967.

\bibitem{joachims1998making}
T.~Joachims.
\newblock Making large-scale svm learning practical.
\newblock Technical report, Technical Report, SFB 475:
  Komplexit{\"a}tsreduktion in Multivariaten Datenstrukturen, Universit{\"a}t
  Dortmund, 1998.

\bibitem{jolliffe2003modified}
I.~T. Jolliffe, N.~T. Trendafilov, and M.~Uddin.
\newblock A modified principal component technique based on the lasso.
\newblock {\em Journal of computational and Graphical Statistics},
  12(3):531--547, 2003.

\bibitem{journee2010generalized}
M.~Journ{\'e}e, Y.~Nesterov, P.~Richt{\'a}rik, and R.~Sepulchre.
\newblock Generalized power method for sparse principal component analysis.
\newblock {\em Journal of Machine Learning Research}, 11(Feb):517--553, 2010.

\bibitem{krauthgamer2015semidefinite}
R.~Krauthgamer, B.~Nadler, D.~Vilenchik, et~al.
\newblock Do semidefinite relaxations solve sparse pca up to the information
  limit?
\newblock {\em The Annals of Statistics}, 43(3):1300--1322, 2015.

\bibitem{KwonL10}
J.~Kwon and K.~M. Lee.
\newblock Visual tracking decomposition.
\newblock In {\em IEEE Conference on Computer Vision and Pattern Recognition
  (CVPR)}, 2010.

\bibitem{LeiZD16}
Q.~Lei, K.~Zhong, and I.~S. Dhillon.
\newblock Coordinate-wise power method.
\newblock In {\em Advances in Neural Information Processing Systems (NIPS)},
  pages 2056--2064, 2016.

\bibitem{lin2007projected}
C.-J. Lin.
\newblock Projected gradient methods for nonnegative matrix factorization.
\newblock {\em Neural Computation}, 19(10):2756--2779, 2007.

\bibitem{mackey2009deflation}
L.~W. Mackey.
\newblock Deflation methods for sparse pca.
\newblock In {\em Advances in Neural Information Processing Systems}, pages
  1017--1024, 2009.

\bibitem{MoghaddamWA07}
B.~Moghaddam, Y.~Weiss, and S.~Avidan.
\newblock Fast pixel/part selection with sparse eigenvectors.
\newblock In {\em IEEE International Conference on Computer Vision (ICCV)},
  2007.

\bibitem{NaikalYS11}
N.~Naikal, A.~Y. Yang, and S.~Sastry.
\newblock Informative feature selection for object recognition via sparse
  {PCA}.
\newblock In {\em IEEE International Conference on Computer Vision (ICCV)},
  2011.

\bibitem{PaisitkriangkraiSZ09}
S.~Paisitkriangkrai, C.~Shen, and J.~Zhang.
\newblock Efficiently training a better visual detector with sparse
  eigenvectors.
\newblock In {\em IEEE Conference on Computer Vision and Pattern Recognition
  (CVPR)}, 2009.

\bibitem{PaisitkriangkraiSZ11}
S.~Paisitkriangkrai, C.~Shen, and J.~Zhang.
\newblock Incremental training of a detector using online sparse
  eigendecomposition.
\newblock {\em {IEEE} Transactions on Image Processing (TIP)}, 2011.

\bibitem{Polik2007}
I.~P\'{o}lik and T.~Terlaky.
\newblock A survey of the s-lemma.
\newblock {\em SIAM Review}, 49(3):371--418, July 2007.

\bibitem{Powell1973}
M.~J.~D. Powell.
\newblock On search directions for minimization algorithms.
\newblock {\em Mathematical Programming}, 4(1):193--201, Dec 1973.

\bibitem{robbins1985convergence}
H.~Robbins and D.~Siegmund.
\newblock A convergence theorem for non negative almost supermartingales and
  some applications.
\newblock In {\em Herbert Robbins Selected Papers}, pages 111--135. Springer,
  1985.

\bibitem{ShenPZ11}
C.~Shen, S.~Paisitkriangkrai, and J.~Zhang.
\newblock Efficiently learning a detection cascade with sparse eigenvectors.
\newblock {\em {IEEE} Transactions on Image Processing (TIP)}, 2011.

\bibitem{SongBP15}
J.~Song, P.~Babu, and D.~P. Palomar.
\newblock Sparse generalized eigenvalue problem via smooth optimization.
\newblock {\em IEEE Transactions on Signal Processing}, 63(7):1627--1642, 2015.

\bibitem{sriperumbudur2011majorization}
B.~K. Sriperumbudur, D.~A. Torres, and G.~R. Lanckriet.
\newblock A majorization-minimization approach to the sparse generalized
  eigenvalue problem.
\newblock {\em Machine learning}, 85(1):3--39, 2011.

\bibitem{Sriperumbudur2007}
B.~K. Sriperumbudur, D.~A. Torres, and G.~R.~G. Lanckriet.
\newblock Sparse eigen methods by d.c. programming.
\newblock In {\em International Conference on Machine Learning (ICML)}, pages
  831--838, 2007.

\bibitem{tan2016sparse}
K.~M. Tan, Z.~Wang, H.~Liu, and T.~Zhang.
\newblock Sparse generalized eigenvalue problem: Optimal statistical rates via
  truncated rayleigh flow.
\newblock {\em Journal of the Royal Statistical Society: Series B}, 2018.

\bibitem{thiao2010dc}
M.~Thiao, P.~D. Tao, and L.~An.
\newblock A dc programming approach for sparse eigenvalue problem.
\newblock In {\em International Conference on Machine Learning (ICML)}, pages
  1063--1070, 2010.

\bibitem{Tseng2001}
P.~Tseng.
\newblock Convergence of a block coordinate descent method for
  nondifferentiable minimization.
\newblock {\em Journal of Optimization Theory and Applications},
  109(3):475--494, Jun 2001.

\bibitem{VandaeleGLZD16}
A.~Vandaele, N.~Gillis, Q.~Lei, K.~Zhong, and I.~S. Dhillon.
\newblock Efficient and non-convex coordinate descent for symmetric nonnegative
  matrix factorization.
\newblock {\em IEEE Transactions on Signal Processing}, 64(21):5571--5584,
  2016.

\bibitem{Wang07traceratio}
H.~Wang, S.~Yan, D.~Xu, X.~Tang, and T.~Huang.
\newblock Trace ratio vs. ratio trace for dimensionality reduction.
\newblock In {\em IEEE Conference on Computer Vision and Pattern Recognition
  (CVPR)}, 2007.

\bibitem{yuan2017hybrid}
G.~Yuan, L.~Shen, and W.-S. Zheng.
\newblock A hybrid method of combinatorial search and coordinate descent for
  discrete optimization.
\newblock {\em arXiv preprint}, 2017.

\bibitem{yuan2013truncated}
X.-T. Yuan and T.~Zhang.
\newblock Truncated power method for sparse eigenvalue problems.
\newblock {\em Journal of Machine Learning Research}, 14(Apr):899--925, 2013.

\bibitem{zhang2011celis}
A.~Zhang and S.~Hayashi.
\newblock Celis-dennis-tapia based approach to quadratic fractional programming
  problems with two quadratic constraints.
\newblock {\em Numerical Algebra, Control and Optimization}, 1(1):83--98, 2011.

\bibitem{zhang2012sparse}
Y.~Zhang, A.~d’Aspremont, and L.~El~Ghaoui.
\newblock Sparse pca: Convex relaxations, algorithms and applications.
\newblock In {\em Handbook on Semidefinite, Conic and Polynomial Optimization},
  pages 915--940. Springer, 2012.

\bibitem{zhang2011large}
Y.~Zhang and L.~E. Ghaoui.
\newblock Large-scale sparse principal component analysis with application to
  text data.
\newblock In {\em Advances in Neural Information Processing Systems (NIPS)},
  2011.

\bibitem{zou2006sparse}
H.~Zou, T.~Hastie, and R.~Tibshirani.
\newblock Sparse principal component analysis.
\newblock {\em Journal of computational and graphical statistics},
  15(2):265--286, 2006.

\end{thebibliography}

\vspace{0.2in}
\begin{center}
 {\LARGE\bf Appendix} 
\end{center}
\appendix

\section{Proof of Proposition 1}
\setcounter{proposition}{0}

\begin{proposition} \label{theorem:convergence:coo}
When the cyclic order strategy is used, coordinate descent method is guaranteed to converge to a coordinate-wise minimum of Problem (\ref{eq:coordinate:bound}) that $\forall i,~\bbb{y}_i^* = \arg \min_{\alpha\geq \hat{L}} ~\mathcal{L}(\bbb{y}_i^*+\alpha \bbb{e}_i)$.

\begin{proof}
Note that $\mathcal{L}(\bbb{y})$ is continuous and $\{\mathcal{L}(\bbb{y}^j)\}$ converges monotonically. Assuming that it converges to $\mathcal{L}^*$ with $\lim_{j  \rightarrow \infty} \mathcal{L}(\bbb{y}^j) = \mathcal{L}^*$, we obtain that $\forall \alpha,~i=1,..,m$:
\beq \label{eq:proof:coordinate:wise1}
 \mathcal{L}^* =  \mathcal{L}(\bbb{y}^{j-1}) = \mathcal{L}(\bbb{y}^j) \leq \mathcal{L}(\bbb{y}^{j-1} + \alpha \bbb{e}_i).
\eeq
\noi Therefore, the right-handed side in (\ref{eq:proof:coordinate:wise1}) attains it minimum at both $0$ and $(\bbb{y}^j)_i - (\bbb{y}^{j-1})_i$. Combining with the fact that the subproblem only contains one unique global solution, we have $(\bbb{y}^{j-1})_i=(\bbb{y}^j)_i$. Since the coordinate $i$ is picked using cyclic order, we have: $\bbb{y}^{j-1}=\bbb{y}^j=\bbb{y}^*$ and $\bbb{y}^*$ is a coordinate-wise minimum point.
\end{proof}
\end{proposition}

\section{Proof of Lemma 2}

\setcounter{lemma}{1}
\begin{lemma} \label{lemma:suf:dec}

(\textbf{Sufficient Decrease Condition}) It holds that: $f(\bbb{x}^{t+1}) - f(\bbb{x}^t) \leq  \frac{-\theta \|\bbb{x}^{t+1}-\bbb{x}^t\|_2^2}{(\bbb{x}^{t+1})^T\bbb{Cx}^{t+1}}$.

\begin{proof}

We let $B$ be the working set in the $t$-th iteration and $N\triangleq \{1,2,...,n\}\setminus B$. Since we solve Problem (\ref{eq:subprob}) in the $t$-th iteration, we have:
\beq
&&\textstyle (h(\bbb{x}^{t+1}_B,\bbb{x}^t_N) + \frac{\theta}{2}\|\bbb{x}^{t+1}_B-\bbb{x}^t_B\|_2^2)~/~ g(\bbb{x}^{t+1}_B,\bbb{x}^t_N) \nn\\
&\leq&\textstyle (h(\bbb{z},\bbb{x}^t_N) + \frac{\theta}{2}\|\bbb{z}-\bbb{x}^t_B\|_2^2)~/~{g(\bbb{z},\bbb{x}^t_N)},~\forall \bbb{z}\in \mathbb{R}^k.\nn
\eeq
 \noi We let $\bbb{z}=\bbb{x}_B^{t}$ and combine with the fact that $\bbb{x}_N^{t+1}=\bbb{x}_N^t$, we have:
\beq
&&\textstyle (h(\bbb{x}^{t+1}_B,\bbb{x}^{t+1}_N) + \frac{\theta}{2}\|\bbb{x}^{t+1}-\bbb{x}^t\|_2^2)~/~ g(\bbb{x}^{t+1}_B,\bbb{x}^{t+1}_N) \nn\\
&\leq&\textstyle (h(\bbb{x}_B^{t},\bbb{x}^t_N) + 0)~/~{g(\bbb{x}_B^{t},\bbb{x}^t_N)}.\nn
\eeq

\noi Noticing the fact that $h(\bbb{x}^{t}_B,\bbb{x}^{t}_N) = \frac{1}{2} (\bbb{x}^{t})^T \bbb{A} \bbb{x}^{t}$ and $g(\bbb{x}^{t}_B,\bbb{x}^{t}_N) = \frac{1}{2} (\bbb{x}^{t})^T \bbb{C}\bbb{x}^{t}$, we have:
\beq
~~~\textstyle ( (\bbb{x}^{t+1})^T \bbb{A} \bbb{x}^{t+1} + \theta\|\bbb{x}^{t+1}-\bbb{x}^t\|_2^2)~/~ (  (\bbb{x}^{t+1})^T \bbb{C} \bbb{x}^{t+1}) \nn\\
\leq \textstyle (  (\bbb{x}^{t})^T \bbb{A} \bbb{x}^{t} )~/~{(  (\bbb{x}^{t})^T \bbb{C} \bbb{x}^{t} )}.~~~~~~~~~~~~~~~~~~~~~~~~~~~~~~~~~~~~~~~~~~~~\nn
\eeq
\noi Moreover, using the structure of the objective function $f(\cdot)$, we obtain: $f(\bbb{x}^{t+1})-f(\bbb{x}^{t}) = \frac{{(\bbb{x}^{t+1})}^T\bbb{Ax}^{t+1}}{{(\bbb{x}^{t+1})}^T\bbb{Cx}^{t+1}} - \frac{{(\bbb{x}^{t})}^T\bbb{Ax}^{t}}{{(\bbb{x}^{t})}^T\bbb{Cx}^{t}} \leq \frac{-\theta\|\bbb{x}^{t+1}-\bbb{x}^t\|_2^2}{{(\bbb{x}^{t+1})}^T\bbb{Cx}^{t+1}} $. Thus, we finish the proof of this lemma.

\end{proof}
\end{lemma}

\begin{figure*} [!th]
\setcounter{subfigure}{0}
\captionsetup{singlelinecheck = true, justification=justified}
\captionsetup{singlelinecheck = on, format= hang, justification=justified, font=footnotesize, labelsep=space}
\centering
       \begin{subfigure}{\fourfigwid}\includegraphics[width=\textwidth,height=\objimghei]{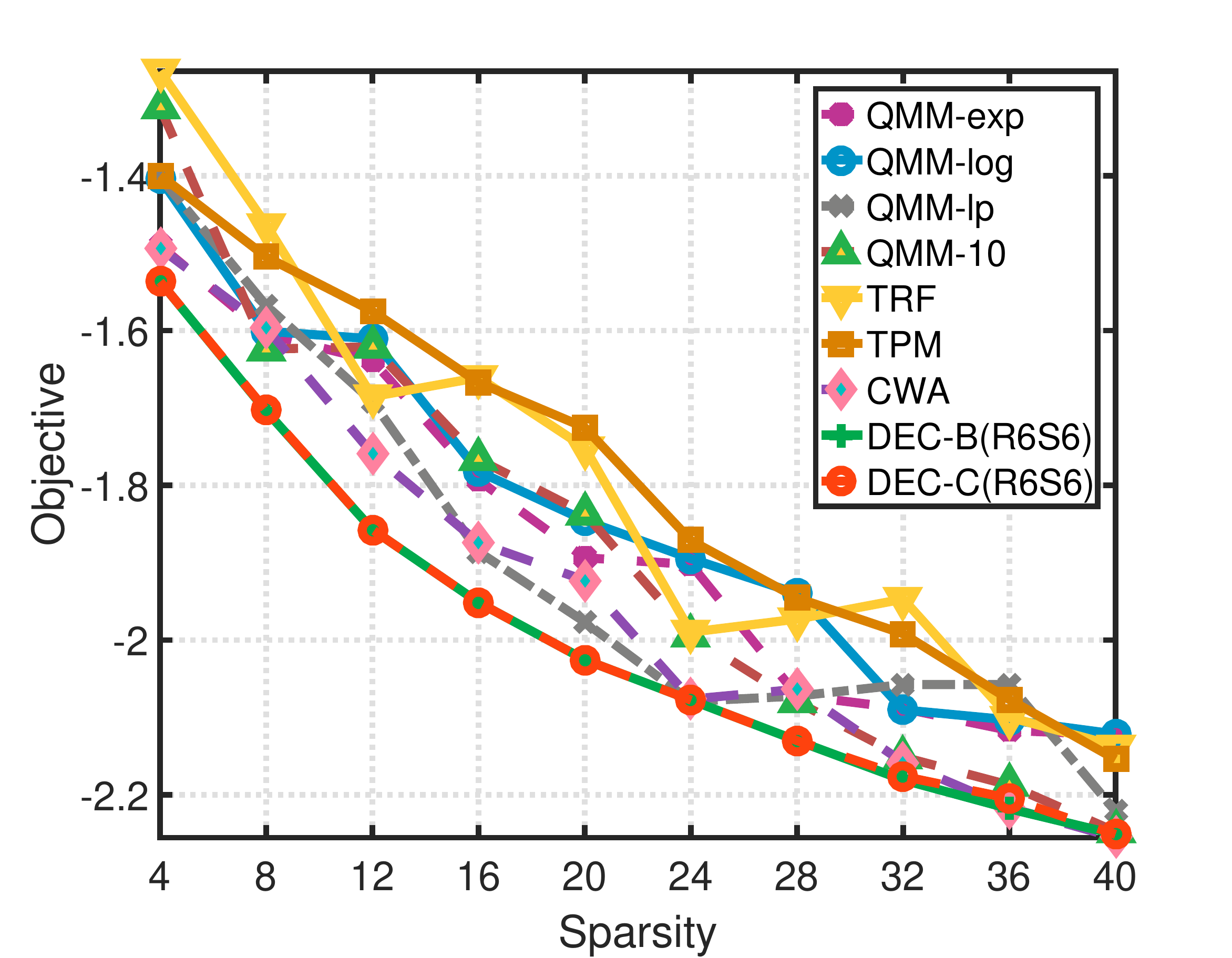}\vspace{-6pt}\caption{\footnotesize randn-100}\end{subfigure}\ghs
      \begin{subfigure}{\fourfigwid}\includegraphics[width=\textwidth,height=\objimghei]{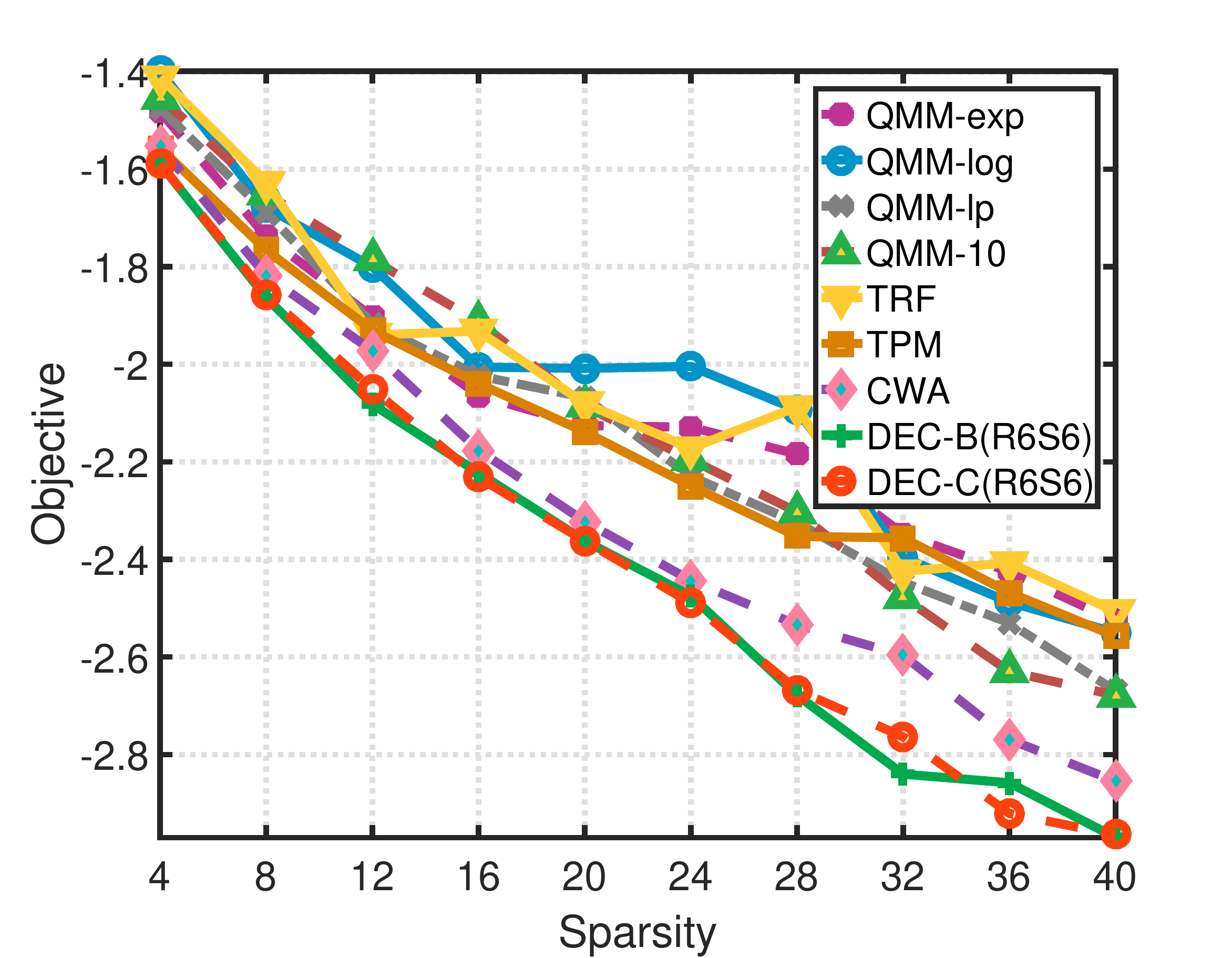}\vspace{-6pt}\caption{\footnotesize randn-500}\end{subfigure}\ghs
      \begin{subfigure}{\fourfigwid}\includegraphics[width=\textwidth,height=\objimghei]{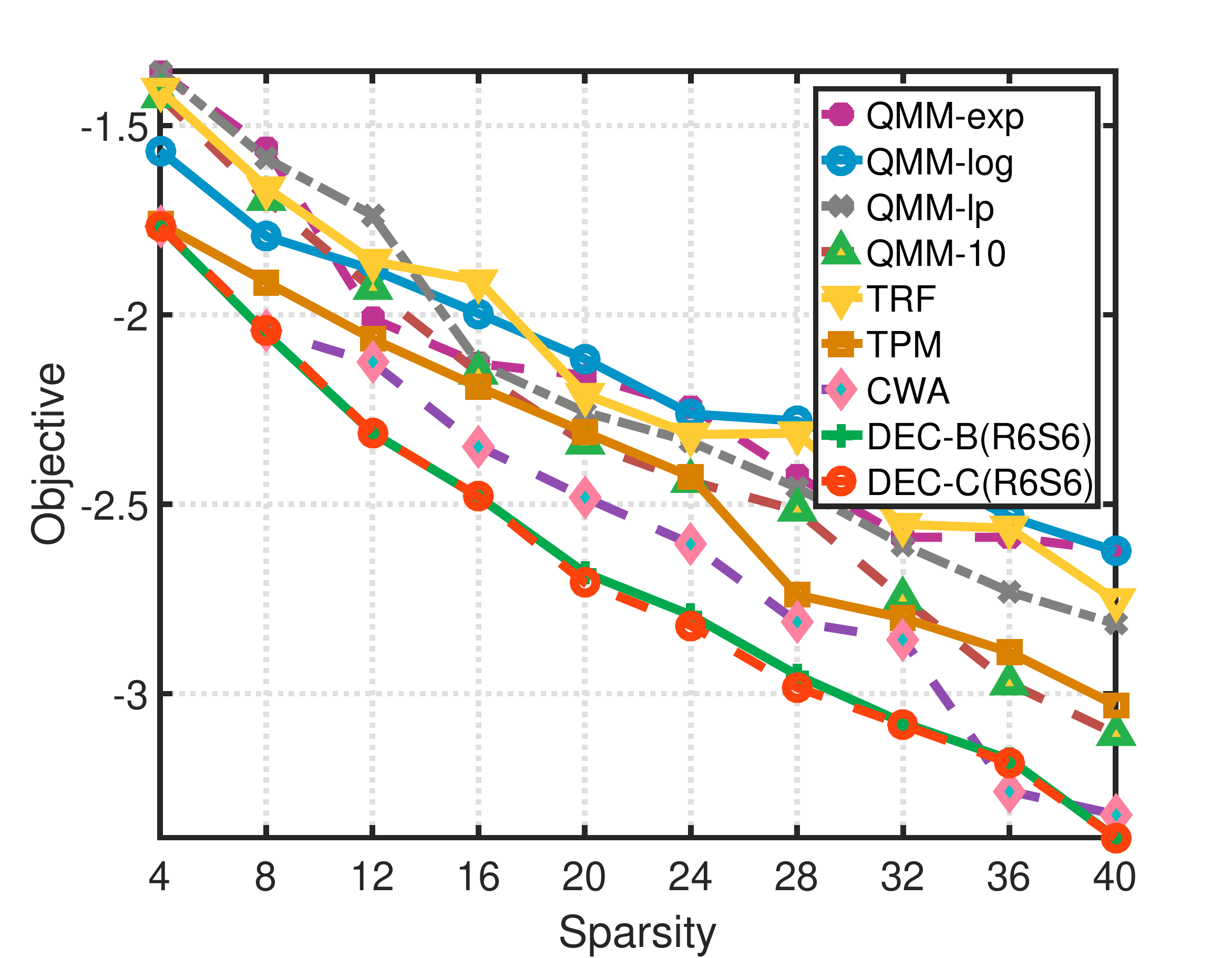}\vspace{-6pt}\caption{\footnotesize randn-1500}\end{subfigure}\ghs
      \begin{subfigure}{\fourfigwid}\includegraphics[width=\textwidth,height=\objimghei]{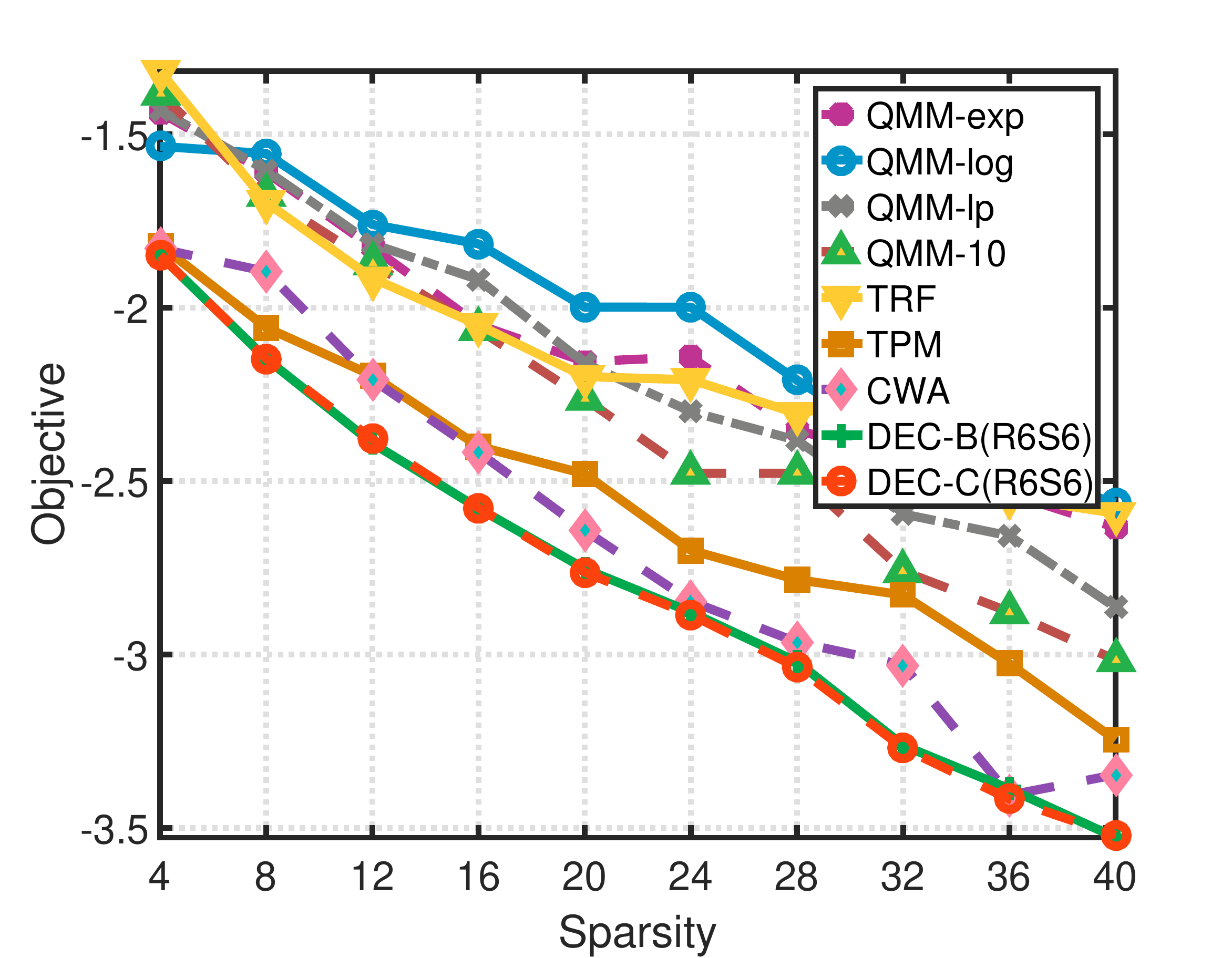}\vspace{-6pt} \caption{\footnotesize randn-2000}\end{subfigure}
\caption{Accuracy of different methods on different data sets for sparse PCA problem with varying the cardinalities.}\label{fig:accuracy:pca-random}

\setcounter{subfigure}{0}
      \centering
      \begin{subfigure}{\fourfigwid}\includegraphics[width=\textwidth,height=\objimghei]{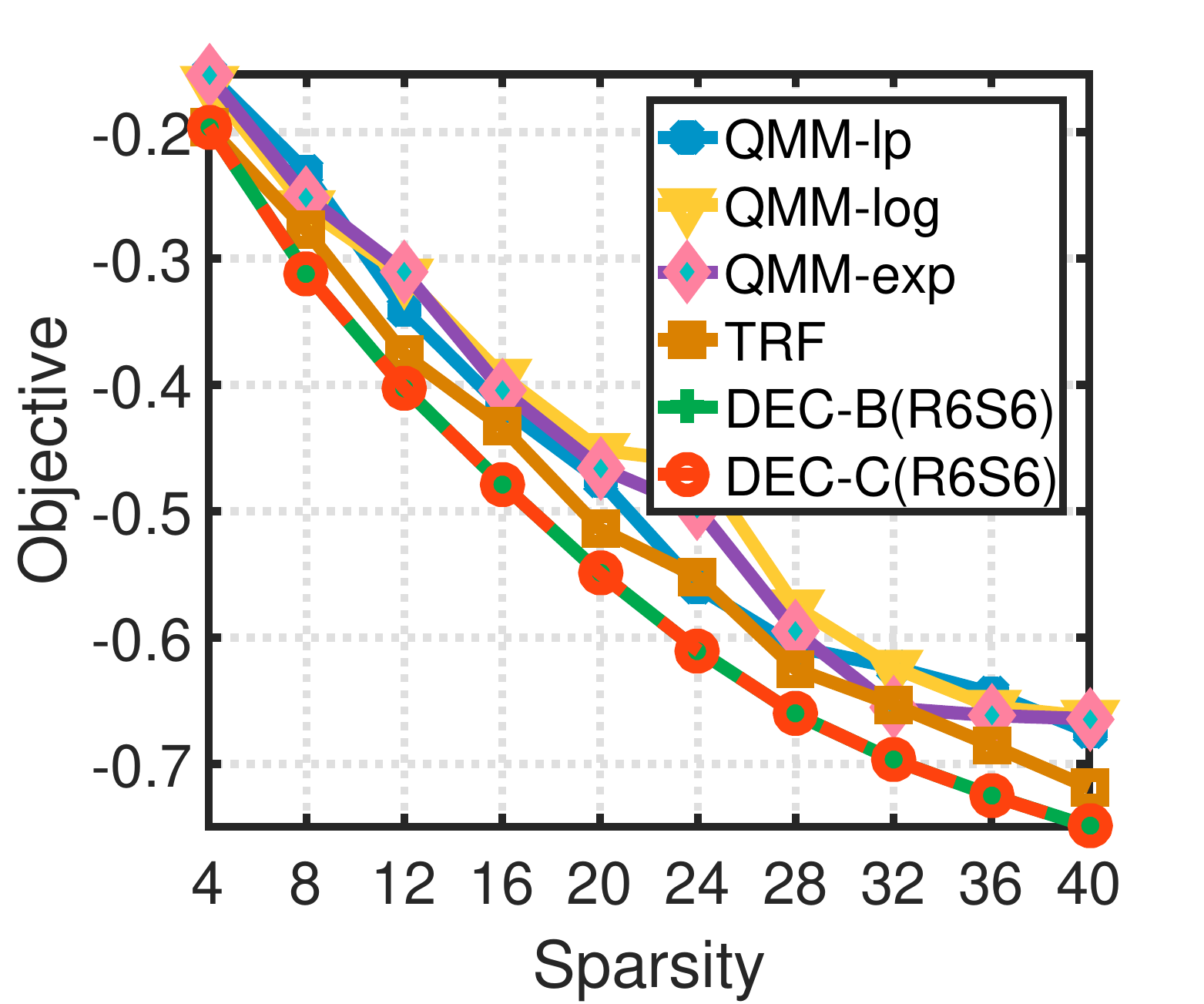}\vspace{-6pt}\caption{\footnotesize randn-100}\end{subfigure}\ghs
      \begin{subfigure}{\fourfigwid}\includegraphics[width=\textwidth,height=\objimghei]{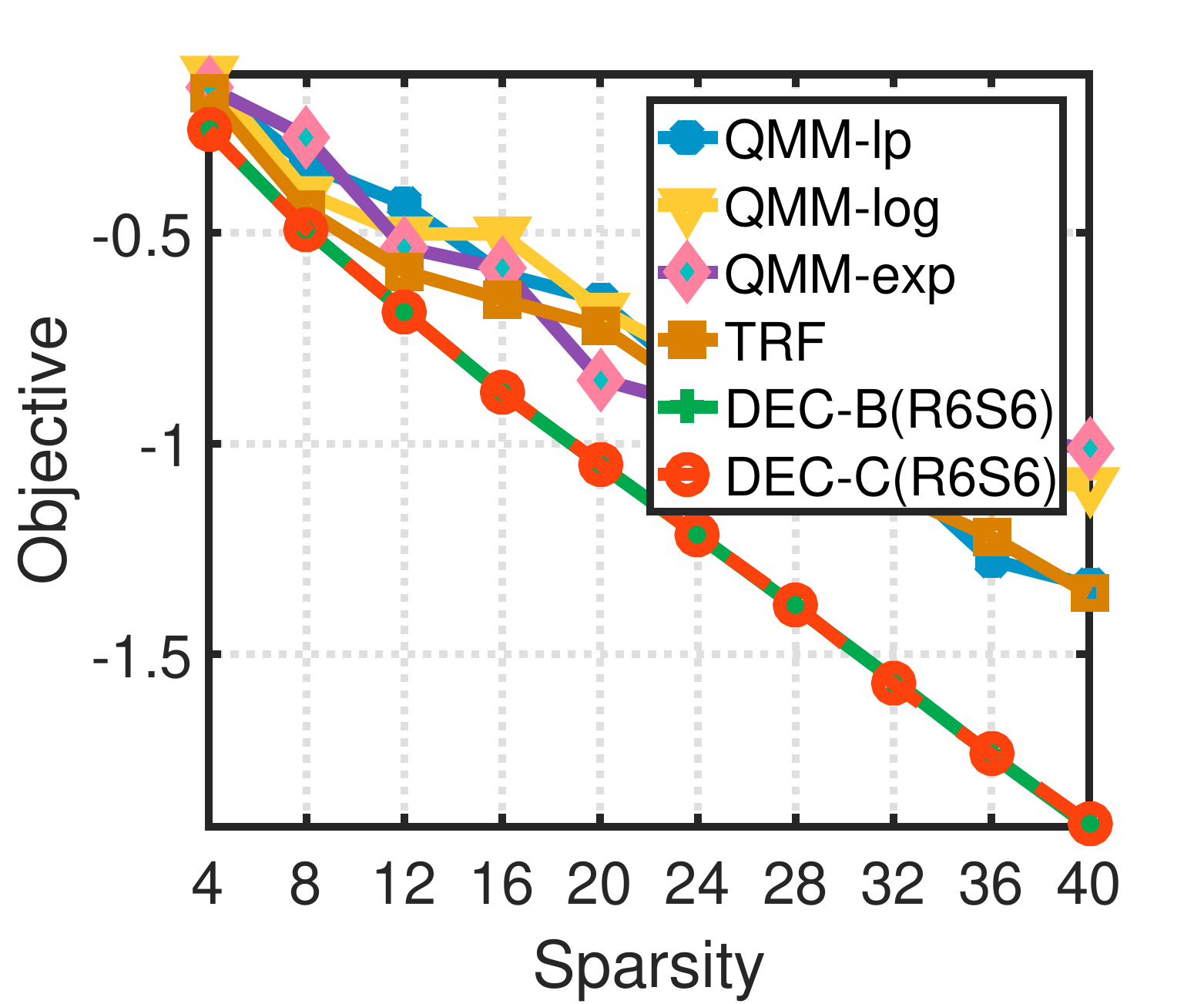}\vspace{-6pt}\caption{\footnotesize randn-500}\end{subfigure}\ghs
      \begin{subfigure}{\fourfigwid}\includegraphics[width=\textwidth,height=\objimghei]{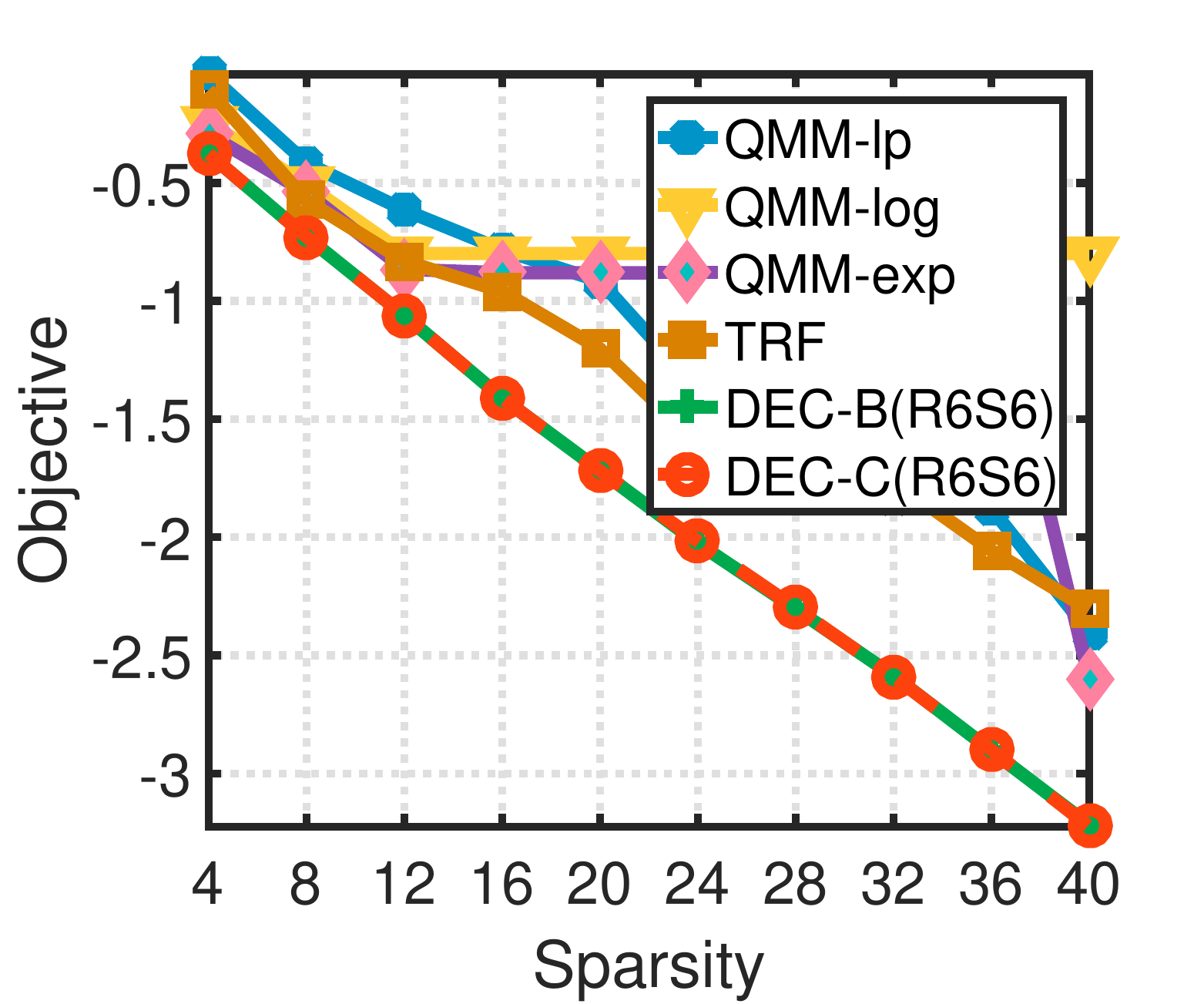}\vspace{-6pt}\caption{\footnotesize randn-1500}\end{subfigure}\ghs
      \begin{subfigure}{\fourfigwid}\includegraphics[width=\textwidth,height=\objimghei]{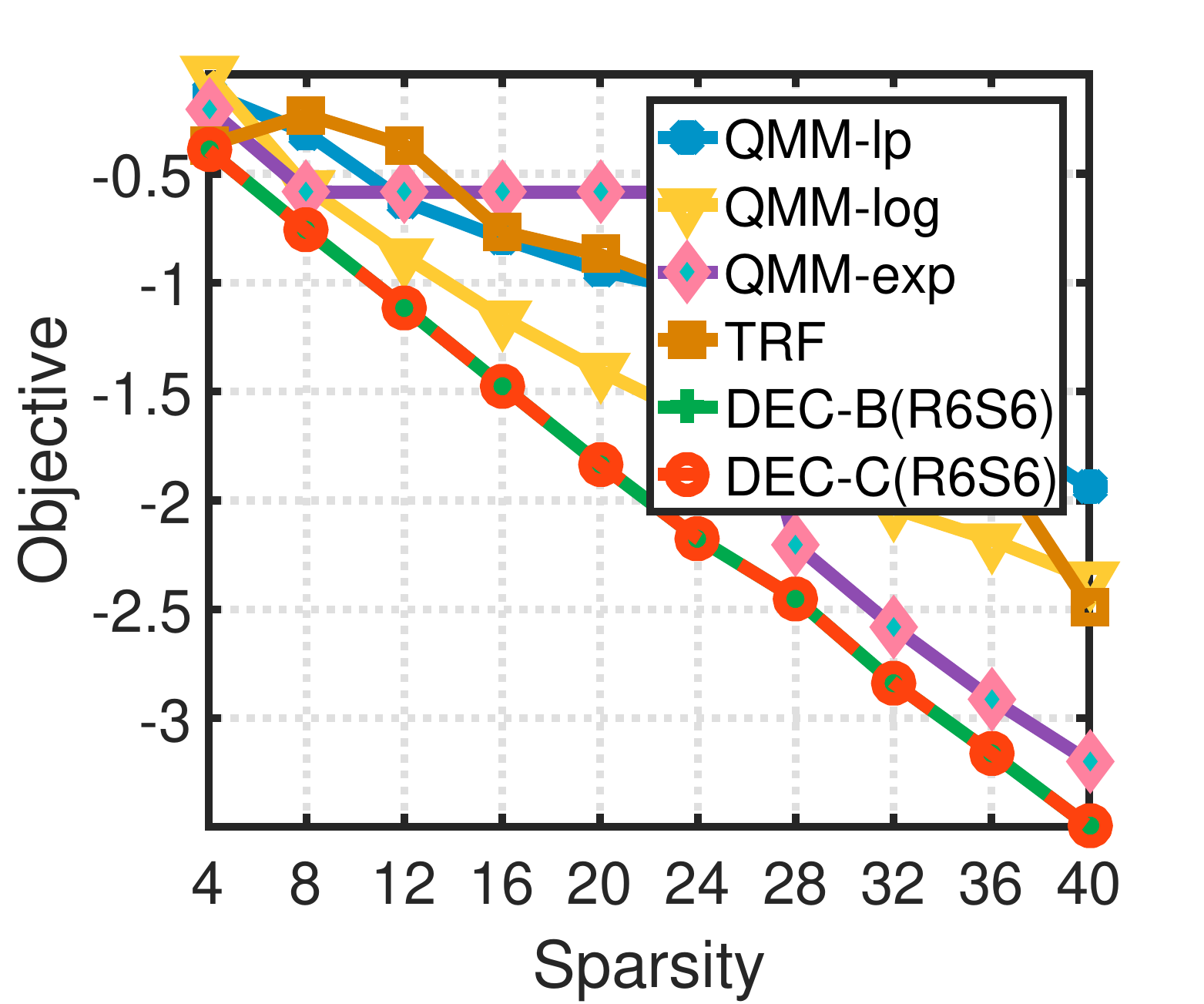}\vspace{-6pt} \caption{\footnotesize randn-2000}\end{subfigure}

\caption{Accuracy of different methods on different data sets for sparse FDA problem with varying the cardinalities.}\label{fig:accuracy:fda-random}

\setcounter{subfigure}{0}
\centering
      \begin{subfigure}{\fourfigwid}\includegraphics[width=\textwidth,height=\objimghei]{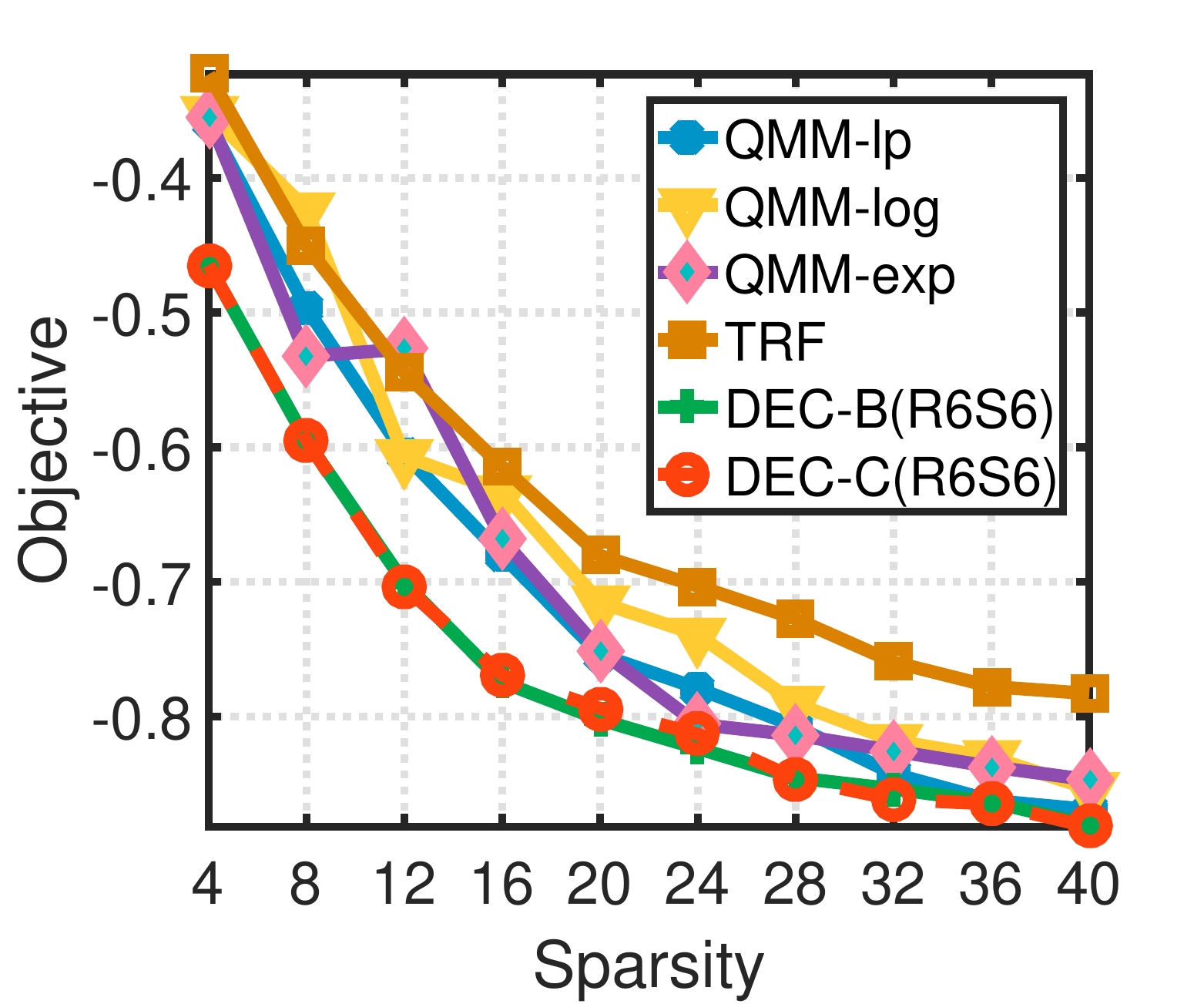}\vspace{-6pt}\caption{\footnotesize randn-100}\end{subfigure}\ghs
      \begin{subfigure}{\fourfigwid}\includegraphics[width=\textwidth,height=\objimghei]{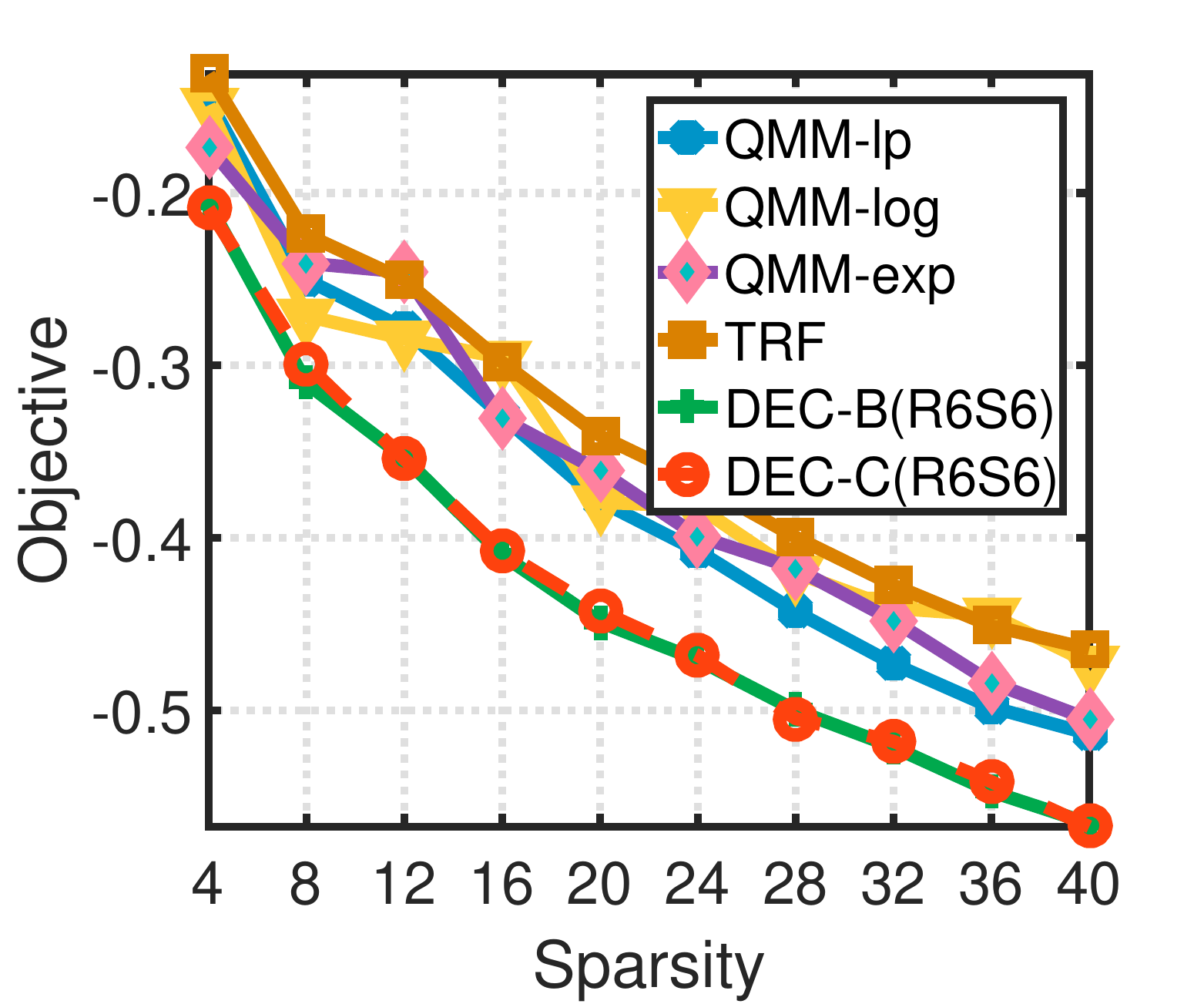}\vspace{-6pt}\caption{\footnotesize randn-500}\end{subfigure}\ghs
      \begin{subfigure}{\fourfigwid}\includegraphics[width=\textwidth,height=\objimghei]{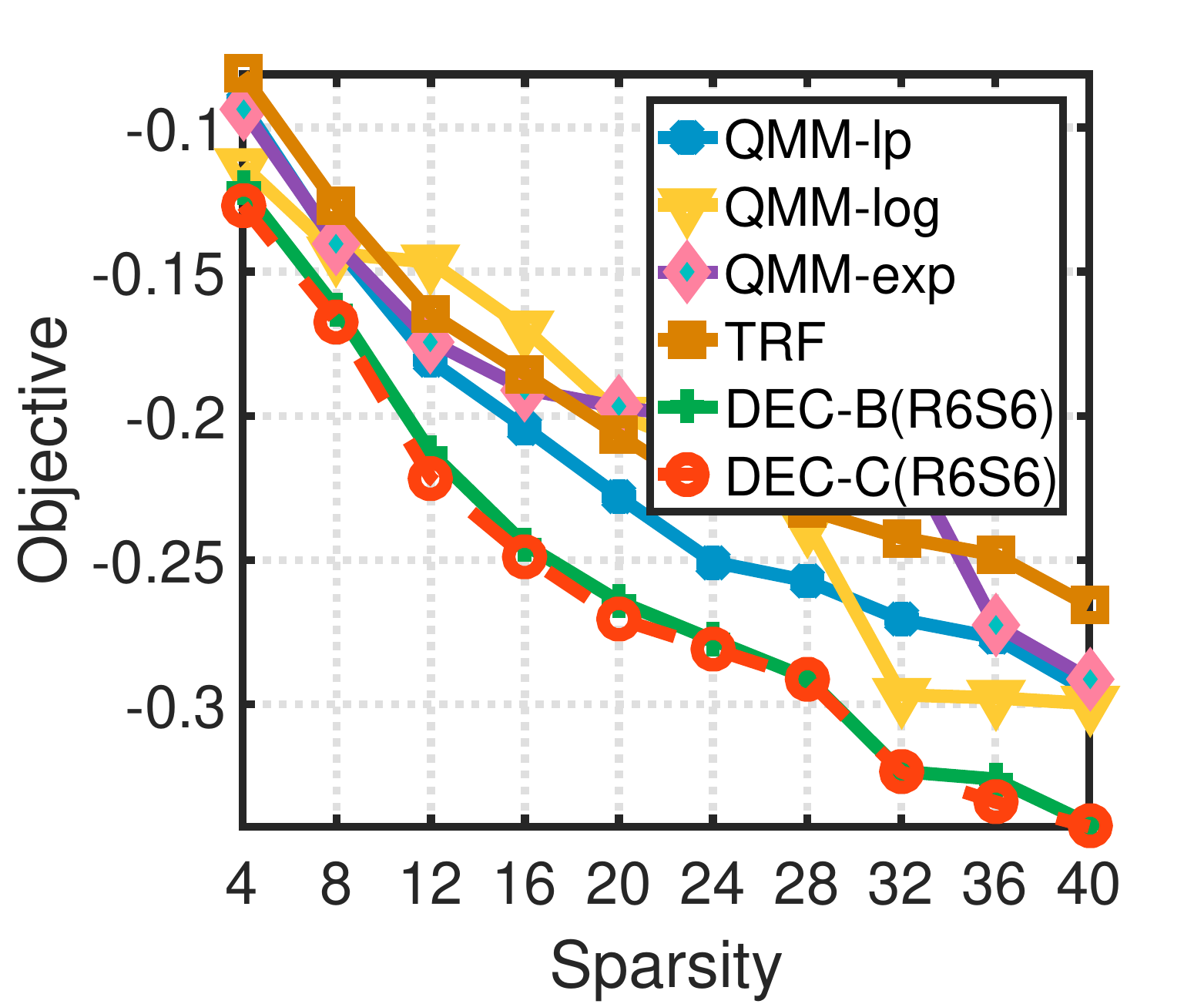}\vspace{-6pt}\caption{\footnotesize randn-1500}\end{subfigure}\ghs
      \begin{subfigure}{\fourfigwid}\includegraphics[width=\textwidth,height=\objimghei]{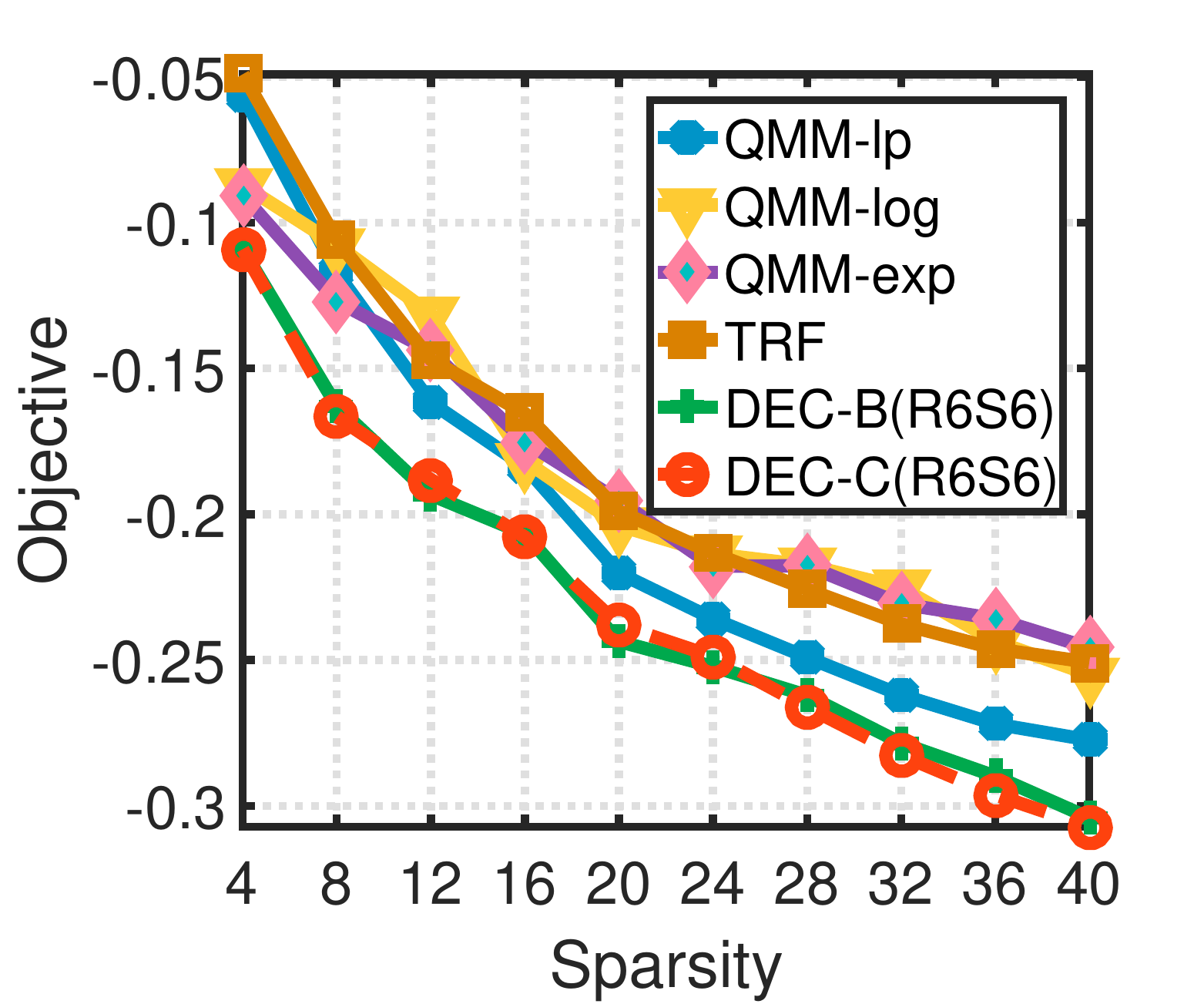}\vspace{-6pt} \caption{\footnotesize randn-2000}\end{subfigure}
\caption{Accuracy of different methods on different data sets for sparse CCA problem with varying the cardinalities.}\label{fig:accuracy:cca-random}
\end{figure*}

\section{Proof of Theorem 2}

We now prove the convergence properties of Algorithm \ref{algo:main}. The following supermartingale convergence result is useful in our analysis \cite{robbins1985convergence}.

\setcounter{lemma}{2}
\begin{lemma} \label{lemma:oooo}
\cite{robbins1985convergence} Let $\bbb{v}_t$ , $\bbb{u}_t$ and $\bbb{\alpha}_t$ be three sequences of nonnegative random variables such that
\beq
&\E[\bbb{v}_{t+1}~|~\mathcal{F}_t]\leq(1+\bbb{\alpha}_t)\bbb{v}_t-\bbb{u}_t,~\forall t\geq0 ~a.s.\nn\\
&\text{and}~\sum_{t=0}^{\infty}~\bbb{\alpha}_t < \infty ~a.s.,
\eeq
\noi where $\mathcal{F}_t$ denotes the collections $\{\bbb{v}_0,...,\bbb{v}_t$, $\bbb{u}_0,..., \bbb{u}_t$, $\bbb{\alpha}_0,...,\bbb{\alpha}_t\}$. Then, we have $\lim_{t\rightarrow \infty} \bbb{v}_t = \chi$ for a random variable $\chi\geq 0$ a.s. and $\sum_{t=0}^{\infty}\bbb{u}_t <\infty$ a.s.

\end{lemma}

%
%

We now present our main results.

\setcounter{theorem}{1}

\begin{theorem} \label{theorem:convergence}

\textbf{Convergence Properties of Algorithm \ref{algo:main}.} Assume that the subproblem in (\ref{eq:subprob}) is solved globally, and there exists a constant $\sigma$ such that $\bbb{x}^{t}\bbb{C}\bbb{x}^{t}\geq \sigma >0$ for all $t$. We have the following results.

\bbb{(i)} When the random strategy is used to find the working set, we have $\lim_{t\rightarrow \infty} \E[\|\bbb{x}^{t+1} - \bbb{x}^t\|] = 0$ and Algorithm \ref{algo:main} converges to the block-$k$ stationary point in expectation.

\bbb{(ii)} When the swapping strategy is used to find the working set with $k\geq 2$, we have $\lim_{t\rightarrow \infty} \|\bbb{x}^{t+1} - \bbb{x}^t\| = 0$ and Algorithm \ref{algo:main} converges to the block-2 stationary point deterministically.

\begin{proof}

We use $\bbb{x}^*$ and $\bar{\bbb{x}}$ to denote any optimal point and any block-$k$ stationary point of (\ref{eq:main}), respectively. We use the notation $\xi^t$ for the entire history of random index selection:
\beq
\xi^t = \{B^0,~B^1,~...,B^t\}\nn
\eeq

\bbb{(i)} We notice that $B^t$ is independent on the past $B^{t-1}$, while $\bbb{x}^t$ fully depends on $\xi^{t-1}$. Taking the expectation conditioned on $\xi^{t-1}$ for the sufficient descent inequality in Lemma \ref{lemma:suf:dec}, we obtain:
\beq\label{eq:descent:after}
&&\textstyle\E[f(\bbb{x}^{t+1})|\xi^t] - f(\bbb{x}^t) \nn\\
&\leq&  \textstyle  - \E [\frac{\theta \|\bbb{x}^{t+1}-\bbb{x}^t\|_2^2}{(\bbb{x}^{t+1})^T\bbb{Cx}^{t+1}}|\xi^t]\nn\\
&\overset{(a)}{\leq}&\textstyle    - \frac{\theta}{\sigma}\E [\|\bbb{x}^{t+1}-\bbb{x}^t\|_2^2|\xi^t]\nn\\
&=& \textstyle   - \frac{\theta}{\sigma}\frac{1}{C_n^k} \sum_{i=1}^{C_n^k}~\| \mathcal{P}(\mathcal{B}_{(i)},~\bbb{x}^t)-\bbb{x}^t\|_2^2\nn\\
&\overset{(b)}{=} &\textstyle   - \frac{\theta}{\sigma}\cdot  \mathcal{M}(\bbb{x}^t)
\eeq
\noi step (a) uses the assumption that $\bbb{x}^{t}\bbb{C}\bbb{x}^{t}\geq \sigma >0,~\forall \bbb{x}^{t}$ which clearly holds since $\bbb{C}$ is strictly positive and $\bbb{x}^{t} \neq \bbb{0}$; step (b) uses the definition of $\mathcal{M}(\bbb{x}^t)$ in Definition \ref{def:block:k}. Therefore, we have:
\beq\label{eq:expectation}
\textstyle \E[f(\bbb{x}^{t+1})~|~\xi^t] - f(\bbb{x}^{*})\leq f(\bbb{x}^t)- f(\bbb{x}^{*})  - \frac{\theta \mathcal{M}(\bbb{x}^t)}{\sigma}
\eeq
\noi Using the supermartingale convergence theorem given in Lemma \ref{lemma:oooo} with $\bbb{v}_t =\E[f(\bbb{x}^{t+1})~|~\xi^t] - f(\bbb{x}^{*}) \geq 0$ and $\bbb{u}_t = \frac{\theta \mathcal{M}(\bbb{x}^t)}{\sigma}$, we have
\beq
\textstyle \lim_{t \rightarrow \infty}~f(\bbb{x}^t) - f(\bbb{x}^*) = \chi~~a.s.\nn
\eeq
\noi for a certain random variable $\chi\geq 0$ and thus the sequence $f(\bbb{x}^t)$ converges to a random variable $\bar{F} = \chi + f(\bbb{x}^*)$. In addition, we have $\lim_{t\rightarrow \infty} f(\bbb{x}^t)-f(\bbb{x}^{t+1}) = 0$ almost surely. From (\ref{eq:descent:after}), we have
\beq
\textstyle \lim_{t \rightarrow \infty} \mathcal{M}(\bbb{x}^t) = 0,~~\lim_{t \rightarrow \infty} \|\bbb{x}^t - \bbb{x}^{t+1}\|=0.\nn
\eeq
\noi Therefore, the algorithm converges to the block-$k$ stationary point. Summing the inequality in (\ref{eq:descent:after}) over $i=0,1,...,t-1$, we have:
\beq
\textstyle\frac{\theta}{\sigma}\cdot \sum_{i=0}^t   \mathcal{M}(\bbb{x}^i) \leq f(\bbb{x}^0)- f(\bbb{x}^t).\nn
\eeq
\noi Using the fact that $f({\bbb{x}}^*)\leq f (\bbb{x}^t)$, we obtain:
\beq
&& \textstyle \frac{\theta}{\sigma} \sum_{i=0}^t \E[\| \mathcal{M}(\bbb{x}^i)~|~\xi^i] \leq f (\bbb{x}^0) - f( {\bbb{x}}^*) \nn\\
\Rightarrow && \textstyle \min_{i=1,...,t} \E[ \mathcal{M}(\bbb{x}^i)~|~\xi^i] \leq    \frac{\sigma (f (\bbb{x}^0) - f( {\bbb{x}}^*))}{t{\theta}}.\nn
\eeq
\noi We conclude that $\bbb{x}^{t}$ converges to the block-$k$ stationary point with $\min_{i=1,...,t} \E[\mathcal{M}(\bbb{x}^i)~|~\bbb{x}^i] \leq    \mathcal{O}(1/t)$.


\bbb{(ii)} We now prove the second part of this theorem. We have the following inequalities:
\beq
f(\bbb{x}^{t+1}) - f(\bbb{x}^t) &\leq& \textstyle  - \frac{\theta \|\bbb{x}^{t+1}-\bbb{x}^t\|_2^2}{(\bbb{x}^{t+1})^T\bbb{Cx}^{t+1}}  \nn\\
&\leq& \textstyle - \frac{\theta}{\sigma}  \|\bbb{x}^{t+1}-\bbb{x}^t\|_2^2\nn
\eeq
\noi Summing this inequality over $i=0,1,...,t-1$, we have:
\beq
&&\textstyle\frac{\theta}{\sigma}\cdot \sum_{i=0}^t   \|\bbb{x}^{i+1}-\bbb{x}^i\|_2^2 \leq f(\bbb{x}^0)- f(\bbb{x}^t)\nn\\
\Rightarrow && \textstyle\min_{i=1,...,t}   \|\bbb{x}^{i+1}-\bbb{x}^i\|_2^2  \leq   \frac{\sigma}{\theta} \frac{f (\bbb{x}^0) - f( {\bbb{x}}^*)}{t }.\nn
\eeq
\noi Using the fact that $f({\bbb{x}}^*)\leq f (\bbb{x}^t)$, we have $\lim_{t\rightarrow \infty} \|\bbb{x}^{t+1} - \bbb{x}^t\| = 0$. Therefore, Algorithm \ref{algo:main} is convergent when swapping strategy is used.

We now prove that Algorithm \ref{algo:main} convergence to a block-2 stationary point $\bar{\bbb{x}}$. Since Algorithm \ref{algo:main} is monotonically non-increasing and converges to a stationary point $\bar{\bbb{x}}$ such that no decrease is made, we have $\bbb{D}_{i,j}\geq 0$ for (\ref{eq:one:dim:selection}). Therefore, it holds that $\min_{\alpha}~f(\bar{\bbb{x}} + \alpha \bbb{e}_{i} - (\bar{\bbb{x}})_j \bbb{e}_j) \geq f(\bar{\bbb{x}}),~\forall i \in \S(\bar{\bbb{x}}),~j\in\Z(\bar{\bbb{x}})$. We have the following result: $f(\bar{\bbb{x}}) \leq f(\bar{\bbb{x}}+ \bbb{d}),\forall \bbb{d}~\text{with}~\|\bbb{d}-\bar{\bbb{x}}\|_0 = 2$. Therefore, $\bar{\bbb{x}}$ is a block-2 stationary point.
\end{proof}
\end{theorem}

\section{Additional Experiments}
We demonstrate the experimental results on the randomized generated data sets for sparse PCA, sparse FDA, and sparse CCA in Figure \ref{fig:accuracy:pca-random}, \ref{fig:accuracy:fda-random} and \ref{fig:accuracy:cca-random}, respectively. These results further consolidate our conclusions drawn in Section \ref{sect:exp}.


\end{document}